\documentclass{article}
\usepackage{graphicx}
\usepackage{amssymb,amsfonts,amsthm}
\usepackage{amsmath,amsthm,amssymb}
\usepackage{amsfonts}

\newtheorem{theorem}{Theorem}[section]
\newtheorem{lemma}[theorem]{Lemma}

\newtheorem{cy}[theorem]{Corollary}

\theoremstyle{definition}

\newtheorem{rk}[theorem]{Remark}

\newcounter{ppp}

\newcommand{\la}{\langle}
\newcommand{\ra}{\rangle}

\begin{document}

\renewcommand{\theequation}{\thesection.\arabic{equation}}

\title{Space functions and complexity of the word problem in semigroups.}
\author{A.Yu.Olshanskii \thanks{The
author was supported in part by the NSF grant DMS 0700811 and by the Russian Fund for Basic Research grant grant 11-01-00945}}
\maketitle

\begin{abstract} We introduce the space function $s(n)$ of a finitely presented semigroup 
$S =\langle A\mid R\rangle .$ 
%(These functions have a natural geometric analog.) 
To define $s(n)$ we consider pairs of words $w,w'$ over $A$ of length at most $n$ 
equal in $S$ and use  relations from $R$ for the transformations $w=w_0\to\dots\to w_t= w'$; $s(n)$ bounds from above the tape space (or computer memory) sufficient to implement all
such transitions $w\to\dots\to w'.$ 
One of the results obtained is the following criterion: 
A finitely generated semigroup $S$ has decidable word problem of polynomial space complexity if and only if $S$ is a subsemigroup of a finitely presented semigroup $H$ with  polynomial space function. 
\end{abstract}


\medskip


%\large

\section{Introduction}\label{intro}

Let $A$ be an alphabet, $A^*$  the set of all words in $A$, and $A^+$  the
set of non-empty words. We will use $|w|$ for the length of
a word $w,$ in particular the empty word $1$ has length $0.$
We write $S=\la A\mid R\ra$ for a semigroup (resp., monoid ) presentation
when $R\subset A^+ \times A^+ $ (resp., $R\subset A^* \times A^* $).

Let $S$ be a semigroup or monoid and $w,w' \in A^*.$ A {\it derivation}
of length $t\ge 0$ from $w$ to $w',$ where $w,w'\in A^+$ or $w,w'\in A^*,$ resp., is a sequence of words 
\begin{equation}\label{rewri}
w = w_0\to w_1\to\dots\to w_t=w',
\end{equation}
where ``$=$'' denotes the letter-for-letter equality, and for $0\le i<t,$ the
word $w_{i+1}$ results from $w_i$ after a defining relation
from $R$ is applied, i.e., $w_i  = ur'v, w_{i+1} = ur''v$ for some words $u, v,$
and $(r',r'')\in R$ or $(r'',r')\in R.$ Two words $w, w'$ represent the same element of $S$
(or they are equal in $S:$  $w=_S w'$) iff there exists a derivation $w\to\dots\to w'.$

The minimal (non-decreasing) function $d(n)\colon \mathbb{N}\to \mathbb{N}$ such that
for every two words $w,w'$ equal in $S$ and having length $\le n,$ there exists
a derivation (\ref{rewri}) with $t\le d(n),$ is called the Dehn function of the presentation 
$S=\langle A\mid R\rangle$ (\cite{Gro}, \cite{Bir}). For {\it finitely presented} $S$ (i.e., both sets $A$ and $R$ are finite),
Dehn functions are usually taken up to equivalence to get rid of the dependence  on
a finite presentation for $S$ (see \cite{MO}). To introduce this equivalence $\sim,$ we write $f\preceq g$
if there is a positive integer $c$ such that
\begin{equation}\label{prec}
f(n)\le cg(cn)+cn\;\;\; for \;\; any \;\;n\in \mathbb{N} 
\end{equation}

For example, we say that a function $f$ is {\it polynomial} if $f\preceq g$ for a polynomial $g.$
From now on, we use the following equivalence for nondecreasing functions $f$ and $g$ on $\mathbb{N}.$

\begin{equation}\label{equiv}
f\sim g\;\;\; if\;\; both\;\; f\preceq g\;\; and\;\; g\preceq f   
\end{equation}

It is not difficult to see that the Dehn function $d(n)$ of a finitely presented semigroup or monoid,
or group $S$ is recursive (or bounded from above by a recursive function) iff the word problem is algorithmically decidable for $S$ (see \cite{Ger1}, \cite{BRS}). In this case, the word
problem can be solved by a primitive algorithm that, given a pair of words $w,w'$ of length $\le n,$ just checks if
there exists a  derivation (\ref{rewri}) of length $\le d(n).$  Therefore the nondeterministic {\it time}  complexity of the word problem in $S$ is bounded from above by $d(n).$ Moreover if $Q$ is a finitely generated subsemigroup (submonoid, subgroup) of $S,$
then one can use the rewriting procedure (\ref{rewri}) for $Q,$ and so the nondeterministic time complexity of the word
problem for $Q$ is also bounded by a function equivalent to $d(n).$  

A converse statement is also true. Assume that the word problem can be solved in
a {\it finitely generated} semigroup $S$ by a nondeterministic Turing machine ($NTM$) with  time complexity $\le T(n),$ where $T(n)$ is a superadditive function (i.e. $T(m+n)\ge T(m)+T(n)$).
Then $S$ is a subsemigroup of a {\it finitely presented} semigroup $H$ with Dehn function   $O(T(n)^2).$ This is proved in \cite{Bir} while a similar statement for groups (but 
with the function $n^2T(n^2)^4$ instead of $T(n)^2$) is obtained in 
\cite{BORS}.
As the main corollary, one concludes that the word problem in a finitely generated semigroup  
(group) $H$ has time
complexity of class $NP$ (i.e., there {\it exists} a nondeterministic  algorithm of polynomial time
complexity, which solves the word problem for $H$) iff $H$ is a subsemigroup (resp., subgroup) of a finitely presented semigroup (resp., group)
with polynomial Dehn function.  

Hence the notion of Dehn function is the (semi)group-theoretical counterpart
of the concept of time complexity for algorithms. It turns out that the
filling length functions introduced earlier in \cite{Gro}, \cite{GR}, \cite{Bir1} (or briefly, space functions) of finitely presented
groups are counterparts of the concept of space complexity of algorithms.
The main theorem of \cite{O} says that  for a finitely generated group $G$ such that the word problem in $G$ is
decidable by a deterministic Turing machine ($DTM$) with space complexity $f(n)$, there  is an embedding of $G$ in a a finitely
presented group $H$ with space function equivalent to $f(n).$ 
In particular the following criterion is obtained: 
A finitely generated group $H$ has decidable word problem of polynomial space complexity if and only if $H$ is a subgroup of a finitely presented group $G$ with a polynomial space function. 

Thus, on the one hand, theorems from \cite{Bir} and \cite{BORS} provide a logical
connection between Dehn functions of semigroups and groups and the time complexity
of their word problems; and on the other hand, similar interrelation of space functions
of groups and the space complexity is obtained in \cite{O}. So it is
natural to fill a gap regarding  space functions of semigroups and the space complexity
of the algorithmic word problem in semigroups. 

In the present paper, we say that the derivation (\ref{rewri}) has space $\max_{i=1}^t |w_i|.$
For two words $w$ and $w'$ equal in $S=\la A| R\ra,$ we denote by $space_S(w,w')$ the minimum of spaces 
of the derivations connecting $w$ and $w',$ and define the value of the {\it space function} $s(n)$
to be equal to $\max (space(w,w'))$ over all pairs $(w,w')$ of equal in $S$ words with $|w|, |w'|\le n.$
An accurate definition of the space complexity (function) $f(n)$ for a Turing machine ($TM$) will be recalled in Subsection \ref{definitions}.
Now we just note that  
the space complexities of machines are taken here up to the same equivalence \ref{equiv} as  the space functions of semigroups. Up to this equivalence, the time and space complexities of the word problem
do not depend on the choice of a finite generator set, see \cite{Bir}, Prop. 2.1.

\begin{theorem} \label{main} Let $S$ be a finitely generated semigroup (monoid) such that the word problem in $S$ is
decidable by a $DTM$ with space complexity $f(n)$. Then $S$ is a subsemigroup (resp., submonoid) of a finitely
presented monoid $P$ with space function equivalent to $ f(n).$ 
\end{theorem}

\begin{rk}
It follows from \cite{MO}, \cite{CMO} that even if $S$ is finitely presented, one
cannot define $P=S$ in Theorem \ref{main}.
%the literal converse statement fails. 
Baumslag's \cite{B}
$1$-relator group  $G=\langle a,b\mid (aba^{-1})b(aba^{-1})^{-1}=b^2\rangle$
is a particular counter-example 
because the space function of $G$ is not bounded from above by any multi-exponential function
(see  \cite{Ger} and  \cite{Pl}) while the space (and time) complexity of the
word problem in $G$ is 
polynomial \cite{MUD}.
\end{rk}

\begin{cy} \label{pspace} The word problem in a {\em finitely generated} semigroup (monoid) $S$  is polynomial space decidable  if and only if
$S$ is a subsemigroup (resp., submonoid) of a {\em finitely presented} monoid $H$ with polynomial
space function.
\end{cy}

We apply our approach to the realization problem: Which functions
$f(n)\colon \mathbb {N}\to\mathbb {N}$ are, up to equivalence, the space functions of finitely presented
semigroups? It is not difficult to find examples of (semi)groups with  linear and exponential  space functions, but it is not easy even to specify a (semi)group with
space function $n^2.$

\begin{cy} \label{realiz} The space complexity $f(n)$ of arbitrary $DTM$ $M$ is equivalent
to the space function of some finitely presented semigroup (or monoid) $P$. 
\end{cy}

This corollary reveals an extensive class of space functions of semigroups, including functions equivalent to 
$[\exp{\sqrt[3]{n}}],$ $[n^k]$ ($k\in \mathbb{N}$), $[n^k\log^l n],$  $[n^k\log^l (\log\log n)^m],$ etc.  Note that we do not assume in the formulation of Theorem \ref{realiz} that the function
$f(n)$ is superadditive (i.e., $f(m+n)\ge f(m)+f(n)$) or grows sufficiently fast. (Compare with theorems in \cite{SBR} and \cite{Bir} on Dehn functions of groups and semigroups.) It follows, in particular, that there exists a finitely presented
semigroup whose space function is not equivalent to any superadditive function. Recall that 
it is unknown if the Dehn function of arbitrary finitely presented group
 is equivalent to a superadditive function; see  \cite{GS}. 
 
 \begin{cy} \label{complete} There is a finitely presented semigroup (and monoid) $P$ with  polynomial space complete word problem and  with polynomial space function.
\end{cy}

We also describe the functions $n^{\alpha}$ which are (up to equivalence)  space functions of semigroups
(and monoids).
As in \cite{O}, our approach is based on a  modification of the proof of Savitch's theorem from \cite{DK} and the proof from \cite{SBR}, where the similar problem was
considered for Dehn functions if $\alpha \ge 4,$ and close necessary and sufficient
conditions were obtained. (See also a dense series of examples with $\alpha\ge 2$ presented in \cite{BB}.)  Now we have $\alpha\ge 1$ in Corollary \ref{alpha} below. 
Also it is worth to note that for space functions, the necessary and sufficient
conditions just coincide. 

To formulate the criterion, 
we  call a real number $\alpha$ {\it computable with space} $\le f(m)$ if there exists a
$DTM$ which, given a natural number $m,$ computes a binary rational approximation of $\alpha$
with an error  $O(2^{-m}),$ and the space of this computation $\le f(m).$

\begin{cy} \label{alpha} For a real number $\alpha\ge 1,$ 
the function $[n^{\alpha}]$  is equivalent to the space function of a finitely
presented semigroup (or monoid) iff $\alpha$ is  computable with space $\le 2^{2^m}.$ 

\end{cy}

It follows that functions  $[n^{\alpha}]$ with any algebraic $\alpha\ge 1$ are  all  space functions of finitely presented semigroups (and monoids), as well as $[n^e],$ $[n^{\sqrt \pi}],$ etc.

\begin{rk} Some of the above statements  
sound similar to propositions from \cite{O}, but they cannot be
deduced from \cite{O} since there, the set of admissible transformations
in derivations is larger than here. (After Bridson and Riley \cite{BR},
we also allowed cyclic permutations and fragmentations of words in \cite{O}.) It is
an open question if Theorem \ref{main} and its corollaries valid  for groups
as well, provided the derivations are based only on the
applications of defining relations, as this is accepted in the present paper. 
\end{rk}
 
Recall that Higman proved in \cite{H} that every recursively presented group
is embeddable in a finitely presented one. A semigroup analog of this theorem
was proved by Murskii in \cite{Mu}. The  approaches to groups and to semigroups are different,
since in the group case one can use conjugations and HNN extensions to synchronize
applications of all defining relations corresponding to one machine command (see \cite{R}).
For semigroup embeddings, Murskii \cite{Mu} and Birget \cite{Bir} use one tape symmetric
input-output machines. We follow this line but one should look after the space
complexity of the constructed machines. Moreover the {\it generalized} space complexity (see
Subsection \ref{definitions} for the definitions) of the latest modification $M_5$ must be equal to the space
complexity of the initial machine $M_0.$  The trick used in Subsection \ref{M2} for
this purpose, works if the machine $M_0$ is deterministic. (Therefore Theorem \ref{main} is
formulated in terms of deterministic space complexity while  \cite{BORS} and
\cite{Bir} consider only nondeterministic time complexity.) 

As in \cite{Mu} and \cite{Bir}, our embedding of $S$ in a finitely presented semigroup $P$ in Theorem \ref{main} is based 
on the commands of the constructed machine. Unlike \cite{Mu} and \cite{Bir}, now we should control the space
function of $P$. Some additional technical difficulties appear because we want to obtain
a monoid embedding (i.e., $1\to 1$) if the semigroup $S$ is a monoid. Monoid relation $v=1$ is
less convenient since the word $v$ can be inserted before/after any letter of any word $w.$ (Note
that only semigroup embedding are under consideration in \cite{Bir}. Whether the embedding from \cite{Mu}
is a monoid embedding or not, if $S$ is a monoid, is also left inexplicit.)

In the remaining part of the proof we introduce derivation trapezia. They 
visualize derivations and make possible to use geometric images, e.g., bands, lenses, cups and caps.
Furthermore, one can remove unnecessary parts of trapezia (e.g., see Lemmas \ref{unlab}--\ref{qa}, \ref{capcup}).  Such parts may not correspond to
subderivation, and they can hardly be defined in a different language.

\section{Machines}\label{ma}

\subsection{Definitions}\label{definitions}
We will use a model of {\it recognizing} $TM$ which is close to the model from \cite{SBR}.

Recall that a {\it (multi-tape) $TM$ with $k$ tapes and $k$
heads} is  a tuple $$M= \langle A, Y, Q,
\Theta, \vec s_1, \vec s_0 \rangle$$ where $A$ is the input alphabet, 
$Y=\sqcup_{i=1}^k Y_i$ is the tape alphabet, $Y_1 \supset A,$
$Q=\sqcup_{i=1}^k Q_i$ is the set of states of the heads of the
machine, $\Theta$ is a set of transitions (commands), $\vec s_1$ is
the $k$-vector of start states, $\vec s_0$ is the $k$-vector of
accept states. ($\sqcup$ denotes the disjoint union.) The sets $Y, Q, \Theta$
are finite.

We
assume that the machine normally starts working with
states of the heads forming the vector $\vec s_1$,  with the head
placed at the right end of each tape, and accepts if it reaches the
state vector $\vec s_0$.  In general, the machine can be turned on
in any configuration and turned off at any time.

A {\em configuration} of tape number $i$ of a $TM$ is a word
$u q v$ where $q\in Q_i$ is the current state of the head,
$u$ is the word to the left of the head, and $v$ is the word to the
right of the head, $u,v\in Y_i^*.$ 
A tape is {\em empty} if $u$, $v$ are empty words.

A {\em configuration} $U$ of the machine $M$ is a word 
$$\alpha_1U_1\omega_1\alpha_2U_2\omega_2... \alpha_kU_k\omega_k$$
where $U_i$ is the configuration of tape $i$, and the endmarkers $\alpha_i, \omega_i$ 
of the $i$-th tape are special separating symbols.  

An {\it input configuration} $w(u)$ is a configuration, where all tapes,
except for the first one, are empty, the configuration of the first
tape (let us call it the {\it input tape}) is of the form $uq$,
$q\in Q_1$, $u$ is a word in the alphabet $A$, 
and the states form the start vector $\vec s_1$. The
{\em accept configuration} is the configuration where the state
vector is $\vec s_0$, the accept vector of the machine, and all
tapes are empty.  (The requirement that the tapes must be
empty is often removed for auxiliary machines which are used 
in construction of bigger machines.)  

To every $\theta\in \Theta,$ there corresponds a command (marked by the
same letter $\theta$), i.e., a pair of sequences of words
 $[V_1,...,V_k]$ and $[V'_1,...,V'_k]$ 
 such that for each $j\le k,$ either both $V_j=uqv$ and $V'_j=u'q'v'$ are configurations of the tape number $j,$ or 
 $V_j =\alpha_j qv$ and $V'_j = \alpha_j q'v',$
 or $V_j = uq\omega_j$ and $V'_j = u'q' \omega_j ,$ or $V_j = \alpha_j q\omega_j$ and $V'_j = \alpha_j q' \omega_j $ ($q,q'\in Q_j$ ).

In order to execute this command, the machine checks if $V_i$
is a subword of the current configuration of the machine,
and if this condition holds the machine replaces $V_i$ by $V'_i$
for all $i=1,\dots,k.$ Therefore we also use the notation:
$\theta: [V_1\to V'_1,\dots, V_k\to V'_k],$  where $V_j \to V'_j$ is
called the $j$-th {\it part} of the command $\theta.$

 Suppose we have a sequence of configurations $w_0,...,w_t$ and a word
 $h= \theta_1\dots\theta_{t}$ in the alphabet $\Theta,$
such that for every $i=1,..., t$ the machine passes from $w_{i-1}$ to
$w_i$ by applying the command $\theta_i$. Then the sequence
$(w_0\to w_1\to\dots\to w_t)$ is said to be a {\it computation with
history} $h.$ In this case we shall write $w_0\cdot h=w_t.$
 The number $t$ will be called the {\em time} or {\em length} of the computation.

A configuration $w$ is called {\em accepted} by a machine $M$ if
there exists at least one computation which starts with $w$ and
ends with the accept configuration. We do not only consider 
deterministic $TM$s, for example, we allow several 
transitions with the same left side. 

A word $u$ in the input alphabet $A$ is said to be {\em accepted} by the machine if the
corresponding input configuration is accepted. (A configuration with
the vector of states $\overrightarrow s_1$
is never accepted if it is not an input configuration.)  The set of all
accepted words over the alphabet $A$ is called the {\em language  ${\cal L}_M$
recognized by the machine $M$}.

If a DTM $M$  halts on an input word $w\in A^*$  at a non-accepting state $\overrightarrow s$
with all tapes empty, 
%\footnote{Eto nado predpolozhit,. t.k. inache ne budet ner-va $S(n)\le S'(n)$.} , 
then one says that $M$ {\it rejects} $w.$ Speaking on deterministic $TM,$ we will
assume that every input configuration is either accepted or rejected, i.e., $M$
may not operate forever being switched on at an input configuration. In other words,
we consider DTM-s $M$ with recursive languages ${\cal L}_M .$

  Let $|w_i|_a$ ($i=0,...,t$) be the number of tape letters (or tape squares) in 
  the configuration $w_i$. (As in \cite{SBR}, the tape letters are called $a$-{\it letters}.) 
  Then the maximum of all $|w_i|_a$ will be called the {\em space
of computation} $C: w_0\to w_1\to\dots\to w_t$
and will be denoted by $space_M(C)$. 
If $u\in A^*$ then, by definition, $space_M(u)$ is the minimal space of the computation
that accepts or rejects the corresponding input configuration $w=w(u).$

The number $S(n)=S_M(n)$ is the maximum of the numbers $space(u)$ over all
words $u\in A^*$ with $|u|\le n.$
The function $S(n)$ will be called the {\it space complexity} of the DTM $M$. 

The definition of the {\it generalized space complexity}
$S'(n)=S'_M(n) $ is  similar to the definition of  space complexity
but we consider arbitrary pair $w_0, w_t$ of configurations which
can be connected by a computation $w_0\to\dots\to w_t$ (not just input
configurations as in the definition of  $S(n)$).
We define $space_M(w_0,w_t)$ to be the minimal
space of computations connecting $w_0$ and $w_t,$ and $S'(n)$ is the minimal function
that bounds from above all numbers $space_M(w_0,w_t)$ under the condition $|w_0|_a, |w_t|\le n.$  It is clear that
$S(n)\leq S'(n).$

\subsection{Input-output machine} \label{inout}

Assume that $S$ is a semigroup (or monoid) generated by a finite set $A,$ and the
word problem in $S$ is decidable by a DTM  $M_0$ with space function $S_0(n).$ We
define this more exactly as follows. The set of input words of $M_0$ consists of
the words $uv',$ where $u$ is a word in $A$ and $v'$ is a word in a disjoint alphabet $A'$
which is a copy of $A.$
Let $u$ and $v$ be two words over $A,$ and
$|u|+|v|\le n.$ Then (1) for a copy $v'$ of $v$ in  $A',$ the
word $uv'$ belongs to the language ${\cal L}_0$ of $M_0$ iff $u=_S v,$ (2) every
input word of $M_0$ of length $\le n$ is accepted or rejected with space
$\le S_0(n),$ and (3) $S_0(n)$ is the minimal function with Property (2).

However to obtain a Higman embedding of $S$
into a finitely presented monoid we are not able to simulate the work of
the recognizing machine $M_0$ by semigroup relations but following Murskii \cite{Mu} and Birget \cite{Bir} (and preserving the space complexity),
we first transform it into an input-output machine $M_1.$ 
We will see that if $u$ is an input word and $v$ is an output word for a computation of $M_1,$ then
$u=_S v.$

The definition of an input-output TM is similar to the definition of a recognizing TM,
but the first tape is an ``input-output tape'' that holds the initial input and the
final output. 
To define the space
complexity $S(n)$ of an input-output machine, one consider input-output computations,
where both the input word $u$ and the output word $v$ are of length at most $n.$
The definition of the generalized space complexity $S'(n)$ is
similar to the definition of  space complexity
but we consider arbitrary 
computations $w_0\to\dots\to w_t$ with $\max\{|w_0|_a, |w_t|_a\}\le n$, not just input-output
computations as in the definition of  $S(n)$.  

We will assume that the words in $A^*$ are ShortLex ordered.

\begin{lemma}\label{M1} There exists an input-output DTM $M_1$ such that

(a) for every input word $u\in A^*$, there is an input-output computation $C$ of
$M_1$ with input $u,$ and the output is the least word $v$ equal to $u$ in $S;$

(b) the space complexity $S_1(n)$ of $M_1$ is equivalent to the space complexity $S_0(n)$ of $M_0;$ 

(c) depending on the state, any configuration $w$ of the computation $C$  contains
either (i) a copy of the word $u$ on one of the tapes or (ii) the output $v$ on the input-output tape, and to obtain the output configuration $w(v)$ in Case (ii), it remains to erase all other tapes and accept; also  we have $space(C)=S_1(|u|)$ in Case (ii);

(d) a configuration with the start vector of states $\overrightarrow s_1$ cannot be reached after an application of a command of $M_1$ to any configuration.  

(e) if a configuration with vector of states $\overrightarrow s_0$ results  
after an application of a  command of $M_1,$ then this command is the unique accepting command. 
 
\end{lemma}

\proof The machine $M_1$ has two tapes more than $M_0$. At first it writes a copy of $u$ on an
extra-tape $T.$ Then it writes a current word $v\le u$ (starting with the least $v;$ 
$v=1$ if $S$ is a monoid) on another extra-tape $T',$ and writes down the word $uv',$
where $v'$ is a copy of $v$ in a disjoint alphabet, 
on the input tape of $M_0$ (which is also the input-output tape of $M_1$). Then
$M_0$ starts working to check whether $u=_Sv$ or not. If ``yes'', then $M_1$ rewrites $v$ onto
the output tape, cleans up all other tapes, and accepts. Otherwise $M_1$ cleans up the tapes
of $M_0,$ replaces the
word $v$ by the next word $v_+\le u$ on $T',$ and repeats the cycle with $v_+.$

Since $u=_S u,$ sooner or later the machine $M_1$ accepts $u$ with Property (a). The first part of (c)
follows from the above algorithm as well. Since the current word $v$ is not longer than $u$
and $|uv'|\le 2|u|,$ we have $S_{M_1}(n)\le 3S_{M_0}(2n),$ and so $S_1(n)\preceq S_0(n).$
Now Property (b) and the second part of (c) follow from the inequality $S_0(n)\preceq S_1(n)$ which can be easily
provided if one forces the machine $M_1$ to check {\it every} pair $(\bar u,\bar v)$ with $|\bar u |, |\bar v|\le |u|$
(even the shortest $v_0$  with $v_0=_S u$ is already found). 
To obtain Property (d), it suffices to add special states for input configuration:
the first command changes these states, and the state letters from $\overrightarrow s_1$
do not occur in other commands. Similarly, one obtains Property (e).
\endproof

\subsection{Machine with equal space complexity and generalized space complexity}\label{M2}
  
In this subsection, we construct an NTM $M_2$ which inherits 
the basic
characteristics of the DTM $M_1$ and has  equivalent generalized space complexity
and space complexity.    
For this goal we adapt the approach from \cite{O}
to input-output machines. 

Assume that $M_1$ has $k$ tapes, and let its first tape  be the input-output tape. Then 
we add a tape numbered $k+1,$
which is empty for input/output configurations, and we
organize the work of the 3-stage machine $M_2$ as a sequential work
of the following machines $M_{21}$, $M_{22},$ and $M_{23}.$ 

The machine $M_{21}$ uses only one
command $\theta_*$ that does not change states and adds one square
with an auxiliary letter  $*$ to the $(k+1)$-st tape, i.e., the command $\theta_*$ has the form
$$[q_1\omega_1\to q_1\omega_1, \alpha_2q_2\omega_2\to \alpha_2q_2\omega_2,\dots, \alpha_kq_k\omega_k\to \alpha_kq_k\omega_k, q_{k+1}\omega_{k+1}\to *\;q_{k+1}\omega_{k+1}]$$ The machine $M_{21}$ can execute this 
command arbitrarily many times while the tapes numbered $1,\dots,k$ keep 
the copy of an input configuration of $M_1$ unchanged.  Then  a connecting rule 
$\theta_{12}: [q_1\to q'_1\omega_1, \dots, \alpha_kq_k\omega_k\to \alpha_kq'_k\omega_k, q_{k+1}\omega_{k+1}\to q'_{k+1}\omega_{k+1}]$ 
changes all states of the heads and switches on the machine $M_{22}.$ Here $(q'_1,\dots, q'_k)=\overrightarrow s_1$ is the vector of start states for $M_1.$

The work of $M_{22}$ on the tapes with numbers $1,\dots,k$ copies 
the work of $M_1.$ But the extension $\theta'$ of  every command $\theta$ of $M_1$ to the $(k+1)$-st tape is defined 
so that its application does not change the current 
space. More precisely,  if a command $\theta$ inserts $m_1$ tape squares and deletes $m_2$
tape squares on the first $k$ tapes, then $\theta'$ inserts $m_2-m_1$ (deletes $m_1-m_2$) squares with 
letter $*$ on the $(k+1)$-st tape if $m_1 - m_2\le 0$ (if $m_1 - m_2\ge 0$).
That is the $(k+1)$-st component of $\theta'$ has the form $q_{k+1}\omega_{k+1}\to *^{m_2-m_1}q_{k+1}\omega_{k+1}$
(resp., $*^{m_1-m_2}q_{k+1}\omega_{k+1}\to q_{k+1}\omega_{k+1}$).
Note that one cannot apply $\theta'$ if $m_1-m_2$ exceeds the current number of
squares on the tape numbered $k+1.$ 

The connecting command $\theta_{23}$ is applicable when $M_{22}$ reaches 
the output configuration on the first $k$ tapes. It changes the states and 
switches on the machine $M_{23}$ erasing all squares on the $(k+1)$-st tape (one by one).

Let $w$ be a configuration of the machine $M_2$ such that $w\cdot \theta_*$ is defined, 
or such that $w$ is obtained after an application of the connecting command $\theta_{12}.$ Then we have
an input configuration on the tapes with numbers $1,\dots,k$ (plus several $*$-s on the $(k+1)$-st
tape). We will denote by $u(w)$ the input word $u$ written on the first tape.  It is an 
input word for the machine $M_1$ as well, and 
the expression  $space_{M_1} u(w)$ makes sense. 

The connecting commands $\theta_{12}$ and $\theta_{23}$ are not invertible in $M_2$
by definition. Therefore every non-empty computation of $M_2$ has history of the
form $h_1h_2h_3$ or $h_1h_2$, or $h_1,$ or $h_2h_3$, or $h_3$, where $h_l$ is the history for $M_{2l},$
($l=1,2,3$).
(To simplify notation
we attribute the command $\theta_{12}$ (the command $\theta_{23}$) to $h_2$ (to $h_3$).)

\begin{lemma} \label{M1M2}

 (a) For every input word $u,$ the machine $M_1$ and $M_2$  give out the same output $v$. 
(b) The   space complexity $S_2(n)$
and the generalized space complexity $S'_2(n)$ of $M_2$ 
are both equivalent to $S_1(n)$. 

\end{lemma}

\proof Assume that $u$ is converted to the output word $v$ by $M_1.$ 
Then $u$ can be converted to $v$ by $M_2$ as well 
because the machine $M_{21}$ can 
insert sufficiently many squares 
(equal to  $space_{M_1}(u)-|u|$) so that the input-output computation of $M_1$ can be simulated 
by $M_{22}.$ Also it is clear from the definition of $M_2$, that every accepting 
computation  for $M_2$ having a history $h_1h_2h_3$ as above, 
simulates, at stage 2, an accepting computation of $M_1$ with history $h_2.$ This proves
Statement (a) 
and equality $ S_{1}(n) = S_{2}(n).$

Assume now that $C: w=w_0\to\dots\to w_t=w'$ is 
a computation of $M_2$ with minimal space for given $w$ and $w',$
and  $h\equiv h_1h_2h_3$ is the history with the above factorization 
(some of the factors $h_i$ can be empty here).  If the word $h_1$ is empty, then $|w_0|\ge\dots\ge |w_n|$ by
the definition of the machines $M_{22}$ and $M_{32}$. Hence the space of this computation
is equal to $|w|_a.$ Similarly, it is $|w'|_a$ if $h_3$ is empty. 
Then let both $h_1$ and $h_3$ be non-empty. It follows that the machine $M_2$
starts (ends) working with a copy of an input (resp., output) configuration of the machine $M_1,$ i.e., the input-output
tape of this configuration contains an input word $u=u(w)$ (output word $v=v(w')$) and the additional $(k+1)$-st
tape has $m$ squares (resp., $m'$ squares) for some $m\ge 0.$ 
We consider two cases.

{\bf Case 1.} Suppose $m \ge space_{M_1}(u)-|u|.$ This inequality says that the additional tape
has enough squares to enable $M_{22}$  to simulate the  computation
of $M_1$ with the input word $u.$ Hence there is an $M_2$-computation  $w_0\to\dots\to w_{n'}$
with history of the form $h'_2h'_3$, and so its space, as well as the space of our original computation, is  $|w|_a.$

{\bf Case 2.} Suppose $m < space_{M_1}(u)-|u|.$ Then there is a computation $w_0\to\dots\to w_{n'}$
such that the commands of its $M_{21}$-stage insert squares until the total number of 
squares of the $(k+1)$-st tape becomes
equal to $space_{M_1}(u)-|u|,$ and then the machines $M_{22}$ and $M_{23}$ work in their standard manner.
The space of this (and the original) computation is $space_{M_1}(u).$

The estimates obtained in cases 1 and 2 
show that $S'_{2}(n)\le\max( S_{1}(n), n)).$ Hence
$$S_{1}(n)=S_{2}(n) \le S'_{2}(n)\le\max( S_{1}(n), n))\sim S_{1}(n), $$ 
and statement (b) is completely proved too.
\endproof

 \subsection{Symmetric machine $M_3$}\label{symm}

For every command $\theta$ of a $TM$, given by a vector $[V_1\to V'_1,\dots,V_k\to V'_k]$,
the vector $[V'_1\to V_1,\dots,V'_k\to V_k]$ also gives  a command
of some $TM$. These two commands $\theta$ and $\theta^{-1}$ are called 
{\em mutually inverse}.

Since the machine $M_1$ is deterministic, the machine $M_2$ has no invertible commands at all.
The definition of the {\it symmetric} machine $M_3= M_2^{sym}$ is the following. Suppose $M_2= \langle X, Y, Q,
\Theta, \vec s_1, \vec s_0 \rangle.$ 
Then by definition, $M_2^{sym} = \langle X, Y, Q,
\Theta^{sym}, \vec s_1, \vec s_0 \rangle,$ where $\Theta^{sym}$
is the minimal {\it symmetric} set containing $\Theta,$ that is,
with every command $[V_1\to V'_1,\dots, V_{k+1}\to V'_{k+1}]$ it contains the inverse command 
$[V'_1\to V_1,\dots, V'_{k+1}\to V_{k+1}]$; in other words,
$\Theta^{sym}=\Theta^+\sqcup\Theta^-,$ where 
$\Theta^+=\Theta$ (the set of positive commands) and $\Theta^-=\{\theta^{-1}\mid \theta\in \Theta\}$
(the set of negative commands).

A computation $w_0\to\dots\to w_t$ of $M_3$ (or other machine) is called {\it reduced} if its history 
is a reduced word. If the history $h=\theta_1\dots\theta_t$ contains a subword $\theta_{i}\theta_{i+1},$ where the commands
$\theta_{i}$ and $\theta_{i+1}$ are mutually inverse, then obviously there is a shorter computation
$w_0\to\dots \to w_{i-1}= w_{i+1}\to\dots\to w_t$ whose space does not exceed the space of the original
one.

\begin{lemma}\label{dvatau} Let $C: w_0\to\dots\to w_t$ be a reduced computation of $M_3$
with history $h=\tau h'\tau',$ where $\tau, \tau'\in\{\theta_{12}^{\pm 1}, \theta_{23}^{\pm 1}\}$
and every command from $h'$ is a command of $M_{22}$ or its inverse. Then the words
$u=u(w_0)$ and $v=u(w_t)$ are equal in $S$ (recall that $u(w)$ is the subword of $w$ written on the input-output tape),
and $space(C)\ge S_1(\max(|u|,|v|)).$ Furthermore, if $\tau'=\theta_{23},$ then $\tau=\theta_{12}$ and
the word $h'$ is positive.

\end{lemma}

\proof Note that $h'$ has no 2-letter subwords $\theta^{-1}\theta',$ where
both $\theta$ and $\theta'$ are positive commands of $M_{22}$ since then different
commands $\theta$ and $\theta'$ would be applicable to the same configuration $w_i\cdot \theta^{-1}=w_0\cdot (\dots\theta^{-1}),$ and so the corresponding commands of $M_1$ would be also
applicable to the same configuration contrary to the determinism of $M_1.$
Hence $h' = g_1g_2^{-1},$ where both $g_1$ and $g_2$ are (positive) histories for $M_{22}.$

Since one may replace $C$ by the inverse computation, it suffices to consider
three cases: (a) $\tau=\theta_{12}, \tau'=\theta_{23},$ (b) $\tau=\theta_{12}, \tau'=\theta_{12}^{-1},$ and (c) $\tau=\theta_{23}^{-1}, \tau'=\theta_{23}.$

{\bf Case (a)}. In this case $g_2$ is empty since $\theta_{23}$ can be applied
only after the unique command of $M_{22}$ corresponding to the accepting command of $M_1$ (see Lemma \ref{M1} (e)). Thus $h'$ is the history of an $M_{22}$-computation,
and by Lemma \ref{M1M2}(a), the corresponding $M_1$-computation converts $u$ into $v.$
So $u=_S v$ by Lemma \ref{M1} (a).
We also have $space(C)\ge S_1(\max(|u|,|v|)$ by Lemma \ref{M1}(c), since $v\le u$
in this case by the definition of $M_1.$

{\bf Case (c)}. The same argument shows now that both $g_1$ and $g_2$ are empty, 
a contradiction. Case (c) is impossible.

{\bf Case (b)}. 
Note that both $g_1$ and $g_2$ are non-empty since a (positive) $M_{22}$-command 
cannot follow by $\theta_{12}^{-1}$ by the Property (d) from Lemma \ref{M1}
since the machine $M_{22}$ copies $M_1.$  If $u=v,$ then two $M_1$-computations $C_1$ and $C_2$ corresponding to
the $M_{22}$-computations $w_0\to\dots\to w_0\cdot g_1$ and $w_t\to\dots w_t\cdot g_2$
have equal the first and the last configurations. Since $M_1$ is deterministic 
it follows that $C_1=C_2.$ But $g_1$ and $g_2$ are completely determined by 
their $M_1$-parts $C_1$ and $C_2.$ Hence we have $g_1=g_2,$ a contradiction. Thus $u\ne v.$

Now by the inequality $u\ne v$ and Lemma \ref{M1} (c), the configuration $w_0\cdot g_1= w_t\cdot g_2$ must contain
(the same) output word on the same input-output tape for the input words  $u$ and $v$ of $M_1.$ We also have
$u=_S v$ by Lemmas \ref{M1}(a) and \ref{M1M2}. Furthermore, by Lemma \ref{M1} (c),
the spaces of $C_1$ and $C_2$ are at least $S_1(|u|)$ and $S_1(|v|),$ respectively.
Therefore $space(C)\ge S_1(\max(|u|,|v|).$ 

The claims of the lemma are proved.

\endproof

We say that a computation $w_0\to\dots\to w_t$ is an input-input 
computation of the machine $M_3$  if both $w_0$ and $w_t$ are input configurations
of $M_2$ (and of $M_3$ as well).

\begin{lemma} \label{inin}(a) For two input configurations $w$ and $w'$, 
of $M_3$, there exists  an 
input-input computation 
$w\to\dots\to w'$ 
iff  $u(w)=_S u(w').$ If $w\ne w',$ then
$space_{M_3}(w,w')= S_1(\max(|u(w)|,|u(w')|).$

(b) Let $C: w\to\dots\to w'$ be a reduced input-output computation of $M_3$ with
$u(w)=u.$ Then 
$space(C)\ge S_1(|u|).$
\end{lemma}
\proof  {\bf (a)} Assume that $u(w)=_S u(w'),$ and let $v$ be the least word equal to $u(w)$ (and to
$u(w')$) in
$S.$ By the definition of $M_1,$ we have a computation $C_1$ of $M_1$ connecting
the input configuration of $M_1$ with input words $u=u(w')$ and the output configuration of $M_1$
with output $v.$ By Lemma \ref{M1} (c), $space(C_1)=S_1(|u|).$ Similarly we have
$C_2$ with input word $u'=u(w')$ and the same output word $v.$ Let $C_3$ and $C_4,$
resp., be the corresponding computations of $M_2$ (see Subsection \ref{M2}). Define $C'$ to be the reduced form of the computation $C_3C_4^{-1}$ of $M_3.$ Then $C'$ connects $w$ and $w'$ and $space(C')\le S_1(\max(|u|,|u'|).$ 

Assume now that $u=u(w)\ne u'=u(w')$ and $C$ is an input-input computation $w\to\dots\to w'$ of $M_3.$
The history of $C$ is $h\equiv h_0\tau_1h_1\tau_2\dots\tau_s h_s,$ where  
 $\tau_i\equiv\theta_{12}^{\pm 1}$ or $\tau_i\equiv\theta_{23}^{\pm 1}$ for $i\le s,$ and the subwords $h_i$-s contain no connecting commands. The connecting commands $\tau_i$-s and
the subcomputations with histories $h_i$-s whose commands correspond to the commands of $M_{21}$ or to the commands of $M_{23}$ (or to inverses)  do not change the content of the input-output tape. By Lemma \ref{dvatau},
the subcomputations of the form $\tau_{i-1}h_i\tau_i,$ where $h_i$ corresponds to $M_{22},$ do not change the content of
the input-output tape modulo the relations of $S.$ Since $h_0$ and $h_s$ must consist of the commands of $M_{21}$ (or inverses) for an input-input computation, we obtain $u=_S u',$ as required.

Furthermore, 
if $u\ne u',$ then $s>1,$ and the subcomputations with histories $\tau_1h_1\tau_2$
and $\tau_{s-1}h_{s-1}\tau_s$ satisfy the
assumption of Lemma \ref{dvatau}, whence $space(C)\ge \max S_1(|u|,|u'|).$ 

The obtained inequalities for $space(C)$ and $space(C')$ complete the proof of Statement (a).

{\bf (b)} Consider the history $h\equiv h_0\tau_1h_1\tau_2\dots\tau_s h_s$ of $C.$ Since $C$ is a reduced input-output computation, we have that $h_0$ 
must consist of (positive) commands of $M_{21},$ 
$\tau_1=\theta_{12},$ and $\tau_s=\theta_{23}.$ 
So $u=u(w)=u(w\cdot h_0),$ $s\ge 2,$ and Statement (b) follows from Lemma \ref{dvatau} applied to the subcomputation
with history $\tau_1h_1\tau_2.$

\endproof

The notation $space(C)\preceq f(n),$ where $n=n(C)$ depends on the computation $C,$ 
will mean further  that  for some constants $c_1,c_2, c_3$ independent of $C,$ we have $space(C) \le c_1f(c_2n)+c_3n.$

\begin{lemma} \label{CM3}
Let $C: w_0\to\dots\to w_t$ be a computation of $M_3$ with the smallest
 space for the fixed $w_0$ and $w_t.$ Then $space(C)\preceq S_1(\max(|w_0|_a,|w_t|_a)).$

\end{lemma} 

 \proof 
 Let us say that a configuration $w$ of $M_3$ has type 1 (resp., 2 or 3) if a
 command of $M_{21}$ (resp., of $M_{22}$ or $M_{23}$) or its inverse is applicable to $w.$ 
 
 {\bf Case 1.} Assume that both $w_0$ and $w_t$ are of type 1. 
 Then using the command inverse to the command of $M_{21},$ we can start with $w_0$ and clean up the tape number $k+1$
 preserving the content $u(w_0)$ of the input-output tape: $C_1: w_0\to\dots\to \bar w.$
 Similarly we have $C_2: w_t\to\dots\to \bar{\bar w}.$ The spaces of these computations
 are $|w_0|_a$ and $|w_t|_a,$ resp. The input configurations $\bar w$ and $\bar{\bar w}$ 
 can be connected by a computation $C_1^{-1}CC_2,$ and by Lemma \ref{inin}, there
 exists an input-input computation $C_3: \bar w\to\dots\to \bar{\bar w}$ of
 space $S_1(\max(|u(\bar w)|,|u(\bar{\bar w})|)=S_1(\max(|u(w_0)|,|u(w_t)|).$ The same upper
 bound holds for the computation $C_1C_3C_2^{-1}: w_0\to\dots\to w_t,$
which proves the lemma in this case.

{\bf Case 2.} Assume that $w_0$ and $w_t$ have types 1 or 3. Taking into account the
previous case, we may assume that $w_t$ is of type 3. If all $w_i$-s have type 3 in $C,$
then their lengths monotonically increase or decrease since $M_{23}$ has only one command.
Hence $space(C)\le \max (|w_0|_a, |w_t|_a).$ Otherwise, by Lemma \ref{dvatau}, the
history of $C$ must have a suffix $\theta_{21} h'\theta_{23} h'',$ where $h''$ 
consists of the commands of $M_{23}$ and $h'$ is a history of an $M_{22}$-computation.
 It follows from the definition of $M_{22}$ that $v=u(w_t)$ is an output word. 
 
 By the definition of $M_2,$
 there exists an $M_2$-computation $C'$ which starts with an input configuration
 $w'$ with the input word $v,$ ends with the output configuration with output also
  $v,$ and has space $S_2(|v|).$  Also there is a computation $C''$ of $M_{23}$
 which deletes several letter in $w_t$ and ends with the same configuration
 as $C'.$ Hence the computation $C'' C'^{-1}$ converts the configuration $w_t$
 into an input configuration $w'$ of length $\le |w_t|,$ and has space $\max(S_2(|v|),|w_t|_a)\preceq  S_2(|w_t|_a)\sim S_1(|w_t|_a)$ by Lemma \ref{M1M2} (b). 
 Hence it suffices to obtain a desired upper estimate for 
 $space_{M_3}(w_0, w').$ But now $w'$ is of type 1. Similarly, if $w_0$ is of type 3, it can be
 replaced by a word $w''$ of type 1. Thus Case 2 reduces to Case 1.
 
 {\bf Case 3.} One of the words $w_0$, $w_t$ (or both) is of type 2. If all $w_i$-s in $C$ are
 of type 2, then the commands from $C$ do not change the lengths, and it is nothing to prove.
 Otherwise one can find $j>i$ such that $w_0,\dots,  w_i$ have equal lengths,
 $w_{j},\dots, w_t$ are of the same length too, and the subcomputation $w_i\to\dots\to w_j$
 satisfies the assumptions of Case 1 or of Case 2.  This completes the proof.
 \endproof 

\begin{lemma}\label{M2M3} The space function $S_3(n)$ and the generalized space function $S'_3(n)$ 
of the machine $M_3$
are both equivalent to the space function $S_1(n)$ of $M_1.$ 
\end{lemma}
\proof 
By Lemma \ref{CM3}, we have $S'_3(n)\preceq S_1(n).$ 
On the other hand, 
by the definition of $S_1(n),$ there is an input-output computation $C: w_0\to\dots\to w_t$ of $M_1$ with the same
input and output words $u=u(w_0)=u(w_t)$  of length $ n,$
and with space $S_1(|u|)=S_1(n).$  
Then one can construct a computation $C': w'_0\to\dots \to w'_{t'}$ of $M_2$ (and of $M_3$) of the same space $S_1(n),$ where the $M_{22}$-portion of $C'$ corresponds to $C.$   
By Lemma \ref{inin} (b), any computation of $M_3$ connecting $w'_0$ and $w'_{t'}$ has
space at least $S_1(|u(w'_0)|)=S_1(|u|)=S_1(n),$ whence $S'_3(n) \ge S_1(n),$  
and the statement of Lemma \ref{M2M3} is proved.

\endproof  

A configuration $w$ of a machine $M$ is called {\it reachable} if there is a computation
$w_0\to\dots\to w_t=w,$ where $w_0$ is an input configuration of $M$.

\begin{lemma}\label{reach} If $w$ is a reachable configuration of the machine $M_3,$
then there is a computation $w_0\to\dots\to w_t=w,$ where $w_0$ is an input configuration
and $|w_0|\le |w|.$
\end{lemma}

\proof By definition, we have a computation $C: w_0\to\dots\to w_t=w$ starting with an input
configuration $w_0.$ We will induct on $|w_t|_a,$ and for fixed $|w_t|_a,$ we will induct on $t.$ 
The base $|w_t|_a=t=0$ is obvious, and moreover, we may always assume that $t>0.$

If $w_{t-1}\to w_t$ is a transition of the machine $M_{21}$ or its inverse, then one can
obtain an input configuration $w'_0$ from $w_t=w$ using a repeated erasing of the auxiliary letter $*$
by the command $\theta_*^{-1}.$
Clearly we have $|w'_0|\le |w|,$ and the statement is true.

If $w_{t-1}\to w_t$ is a transition of the machine $M_{22}$ or its inverse, then
$|w_{t-1}|_a= |w_t|_a,$ and  it remains to apply the inductive hypothesis to the reachable word $w_{t-1}.$ 

Now we assume that $w_{t-1}\to w_t$ is a transition of the machine $M_{23}$ or its inverse. 
By Lemma \ref{dvatau}, the history of $\cal C$ must have a suffix of the form $\theta_{12} h'\theta_{23}h'',$ where
$h'$ (resp., $h''$) is a product of the commands of $M_{22}$ (of the commands of $M_{23}$ or inverses).
Since the subcomputation $C': w_r\to\dots\to w_s$  with history  $ h'$ is an input-output computation of $M_{22},$
there is a computation of $M_1$ with the input  $u=u(w_r)$ and the output $v=u(w_s).$
By Lemma \ref{M1} (a), the word $v$ is not equal in $S$ to a shorter
word. Therefore, by the definition of $M_1,$ there is an input-output computation
%$C''$ 
of $M_1$ with both input and output words equal to $v.$
Then by Lemma \ref{M1M2} (a),
there is a computation of $M_2$ (and $M_3$) with the input and the output equal to $v;$
it starts with a configuration $w',$ where $u(w')=v$ and ends with some $w'',$ with $u(w'')=v$. Note that $v=u(w_s)=u(w_t).$ Therefore
$|w'|_a=|v|\le |w_t|_a,$ and so $|w'|\le |w_t|.$ The command $\theta_{23}$ is applicable to both configurations $w''$ and $w_s,$ and so they can be connected by an $M_3$-computation, where
every command (or its inverse) is a command of $M_{23}.$ 
It follows that  there is an $M_3$-computation
$w'\to\dots\to w''\to\dots\to w_s\to\dots\to w_t,$ and the lemma is proved. \endproof

 \medskip

Below we will treat the NTM $M_3$ as a nondeterministic `input-input' machine, i.e., the 'purpose'
of $M_3$ is to transform an input configuration $w(u)$ to an input configuration $w(v).$ Consider
the following relation
$u\sim v$ on the set of words in the input alphabet:
{\it there exists a (reduced) input-input $M_3$-computation} $w(u)\to\dots\to w(v)$. This is an equivalence relation.
Indeed, it becomes reflexive if one adds computations of length $0.$ The transitivity is obvious, and its symmetry follows from the symmetry of $M_3.$ So we use the term
{\it equivalence} machine (or just $E$-{\it machine}) for a symmetric input-input NTM. 

\subsection{One-tape machine}

It is well known that any NTM is equivalent to a one-tape NTM with the same
space complexity (see Corollary 1.16 in \cite{DK}). But here we take some precautions to preserve
the {\it generalized} space complexity. 

Let $M$ be a  $k$-tape E-machine with an input alphabet $A.$ 
We will construct an equivalent (i.e., defining the same equivalence relation on the words over
$A$)
one-tape E-machine $M'.$ $M'$ has the same input alphabet $A,$ and at the preliminary
stage it inserts the endmarkers $\alpha_1,\dots \omega_k$ and the components $q_{11},\dots, q_{1k}$ of the
vector of start states $\overrightarrow s_1$ of $M,$ that is, these letters become
tape letters of $M',$ and at the first stage, $M'$ converts an input configuration
$\alpha u q_1 \omega$ of $M'$ into the configuration $\alpha\alpha_1 u q_{11}\omega_1\dots
\alpha_k q_{1k}\omega_k q_1 \omega$ (i.e., the  head of $M'$ runs to $\alpha,$ check the left endmarker, inserts
the letter $\alpha_1,$ and then returns to $\omega$ inserting the remaining extra-letters
$q_{11},\omega_1,\dots, \alpha_k, q_{1k}, \omega_k$).

Every configuration $w$ of $M$ is represented by the configuration $W=\alpha w q \omega$
of $M',$ where the state letter $q$ of $M'$ is the vector $(q^{(1)},\dots, q^{(k)})$ of states of $w$ (but all the letters of $w,$ including the extra-letters, are tape letters for $M'$). For every transition
 $w\to w\cdot\theta$ of $M$ with positive command $\theta,$ we construct a computation $C(w,\theta): W\to \dots\to W'$ of $M'$ as follows.
 The first command just changes the state $q$ by the state $q_{\theta},$ i.e., it memorizes $\theta,$ and $M'$ will remember
 $\theta$ until the computation  $C(w,\theta)$ ends. This command involves the endmarker $\omega.$ Then the head of $M'$ goes to the left
 and simulates the application of the command $\theta$ when it meets $q^{(i)}.$ For example,
 if the $i$-th part 
 of $\theta$ is
 $cq^{(i)}a\to dq'^{(i)}b$, then the corresponding
 computation of $M'$ is of the form 
 $$\dots c q^{(i)} q_{\theta}(1) a\dots \to \dots c q^{(i)}q_{\theta}(2) b\dots \to \dots c q_{\theta}(3)q'^{(i)} b\dots \to 
 \dots d q_{\theta}(4)q'^{(i)}b \dots,$$ where $q_{\theta}(1),\dots q_{\theta}(4)$ are
  auxiliary state letters of $M'.$ So the head of $M'$ must reach $\alpha$ (there is a command
  involving $\alpha$) and then returns to $\omega.$
 The last command of this computation $q'_{\theta}\omega\to q'\omega$ forgets $\theta,$ and the state letter $q'$ of $W'$ is just the
 the vector of states of the configuration $w',$ so that $W'$ corresponds to $w'.$ 
 
 \begin{rk} \label{oneone} Two different positive (or two different negative) commands of $M'$
 cannot be applicable to a configuration containing a state letter indexed by some $\theta$. 
 \end{rk}

\begin{lemma}\label{1tape} (a)The E-machines $M$ and $M'$ recognize the same equivalence relation
on the set of input words. 

(b) They have equivalent generalized space functions $S'_M(n)$ and $S'_{M'}(n).$ 

(c) If $M=M_3$ and $W$ is a reachable configuration of $M',$ then there is
a computation $W_0\to\dots\to W,$ where $W_0$ is an input configuration of $M'$ and
$|W_0|_a\le |W|_a+c$ for a constant $c$ independent of $W.$ 

(d) If a reduced computation $W_0\to\dots\to W_t$ of $M'$ has no commands involving $\alpha$ or
has no commands involving $\omega,$ then $t$ is bounded from above by $c_1|W_0|+c_2$
for some constants $c_1,c_2.$

(e) If a computation of the form $Uq_1\omega=W_0\to\dots\to W_t$, where $q_1$ is the start state of $M',$ has no commands involving $\alpha$
then $|W_0|=|W_t|$ and the computation commands  involve tape letters only from the input
alphabet. A non-empty  computation of the form $Uq_1\omega\to\dots\to U'q_1\omega$ has a command involving $\alpha.$

(f) Let an $M'$-computation starts with an input configuration $w_0=\alpha u q_1\omega$
and ends with $w_t=\alpha u' q_1\omega.$ Then $w_t$ is also an input configuration. 

\end{lemma}
\proof Observe that
the computation $C(w,\theta): W\to\dots $ exists iff one may apply $\theta$ to $w$
and $W$ is the configuration of $M'$ corresponding to $w.$  Moreover, if
two configurations $W$ and $W'$ of $M'$ represent some configurations $w$ and $w'$ of $M,$
then they can be connected by a computation $C(w,\theta)$ iff $w'=w\cdot\theta.$

{\bf (a)} By the definition of $M',$ every input-input computation of $M$ can be simulated by
$M'.$ Now let us consider a non-empty reduced input-input computation $C': W_0\to\dots\to W_t$
of $M',$ and denote by $W_{i_1},\dots, W_{i_s}$ ($0<i_1<\dots<i_s<t$) the intermediate
configurations 
representing the configurations of $M$ (i.e., $M'$ do not remember 
the commands of $M$ in these states). Since there
are no other configurations with this property between $W_{i_{m-1}}$ and $W_{i_m}$, all the commands of the subcomputation
$W_{i_{m-1}}\to\dots\to W_{i_m}$ must correspond to the same command $\theta$ of $M.$ By Remark \ref{oneone},
the history of this subcomputation has no subwords of the form $\tau^{-1}\tau'$
(of the form $\tau'\tau^{-1}$), where  both $\tau$ and $\tau'$ are positive commands
of $M'.$
Therefore the subcomputation must be of the form $C(w_{m-1},\theta)$ or $C(w_{m},\theta)^{-1}$ for a positive command $\theta$ of the machine $M,$ and $W_{i_{m-1}}, W_{i_m}$
correspond to $w_{m-1}$ and to $w_m=w_{m-1}\cdot\theta,$ resp., or to  $w_{m-1}=w_{m}\cdot \theta $ and to $w_m,$ resp.

Since the preliminary stage $W_0\to\dots\to W_{i_1}$ (resp., $W_t\to\dots\to W_{i_s}$)
is also deterministic, the pair of input words for $C'$ coincides with the pair of input
words for the computation $C: w_1\to\dots\to w_s$ of $M,$ and so $M$ and $M'$ recognize
the same binary relation.

{\bf (b)} Let now $C': W_0\to\dots\to W_t$ be an arbitrary reduced computation of $M'$ with
$\max(|W_0|_a, |W_t|_a)\le n.$ We define $W_{i_1},\dots, W_{i_s}$ as in Part (a) of the
proof. If $s=0,$ then for every $j,$ $|W_j|_a \le n+c$ for a constant $c$ independent
of the computation since the computations of the form $C(w,\theta)$ and the preliminary
computations (and their subcomputations), up to a constant, do not change the space. So we will assume that $s\ge 1.$
Then as in Part (a), the computation $C'': W_{i_1}\to \dots\to W_{i_s}$ corresponds to a
computation $w_1\to\dots\to w_s$ of $M,$ where $|w_j|$ and $|W_{i_j}|$ are almost
(up to an additive constant) equal. Therefore $w_1$ and $w_s$ can be connected by a computation
$C$ of $M$ of space at most $S'_M(n+c).$ There is a computation $C'''$ of $M'$ corresponding
to $C$ and having almost the same space. If we replace the subcomputation $C''$ of $C'$ by
$C'''$ we get a computation of $M'$ which connects $W_0$ and $W_t$ and has space $\le S'_M(n+c)+c.$
Hence $S'_{M'}(n)\preceq S'_M(n).$ 

Similarly, if we start with a computation $C: w_1\to\dots\to w_s$ of $M$ with $|w_1|_a,|w_s|_a\le n,$
then we can replace it by a computation of $M$ of space at most $S'_{M'}(n+c),$ whence $S'_{M}(n)\preceq S'_{M'}(n),$ as required. 

{\bf (c)} Assume now that $W_t=W$ and $W_0$ is an input configuration in the computation
$C'$ from (b). 
Now we consider
the computation $C: w_0\to w_1\to\dots\to w_s$ of $M,$ where $w_j$ ($j\le s$) corresponds
to the configuration $W_{i_j}$ of $M'.$ By Lemma \ref{reach}, one can find a computation
$ w_0'\to\dots\to w'_{s'}=w_s$ of $M,$ such that $w'_0$ is an input configuration and
$|w'_0|\le|w'_{s'}|.$ Then one can construct a computation $W'_0\to\dots\to W'_{i'_1}\to\dots\to W'_{i'_{s'}},$ where
each  $W'_{i'_j}$ represents $w'_j,$ and therefore 
$|W'_0|\le |W'_{i'_{s'}}|.$ 
This give a computation $$W''\to\dots\to W'_0\to\dots\to W'_{i'_{s'}}=W_{i_s}\to\dots\to W_t,$$
where $W''$ is an input configuration of $M',$ and $|W''|\le |W_t|+c,$ as required, because the
preliminary subcomputation $W''\to\dots\to W'_0$ does not decrease the space and
the subcomputation $W_{i_s}\to\dots\to W_t$ is either empty or a part of a computation $C(w_s,\theta),$
and therefore it can remove a bounded number of tape letters.

{\bf (d)} Follows from the fact that the computation of the form $C(w,\theta)$ involves
both $\alpha$ and $\omega$.

{\bf (e)} The (reduced) work of $M'$ is  deterministic in the beginning:
the head goes to the left until it reaches and checks the endmarker $\alpha.$ This implies Property (e). 

{\bf (f)} We must show that $u'$ is a word in the input alphabet $A.$ For this goal
we can (1) assume that the computation is reduced, (2) consider the inverse computation $w_t\to\dots\to w_0,$ and (3) take into account that
at the preliminary stage, the machine $M'$ verifies (when the head goes to $\alpha$)
if all the letters of the tape word belong to $A.$ 

\endproof

We will use the doubling of the tape alphabet. This is a well-known trick helpful
for simulating of machine commands by (semi)group relations (see \cite{R}). 
Let $Y$ be a tape alphabet of a one-tape machine $M.$ We denote by $Y_l$ and $Y_r$
two disjoint copies of $Y$ (`left' and `right') and replace every configuration $\alpha u q v \omega$ of $M$
by $\alpha u_l q v_r \omega,$ where $u_l$ (resp., $v_r$) is a copy of $u$ in $Y_l$ (in $Y_r$).
Respectively, one modifies every command, e.g., a command $aqb\to cq'd$ is replaced by
$a_lqb_r\to c_lq'd_r.$ The input alphabet is replaced by its copy $A_l\subset Y_l.$ Clearly, one obtain
one-to one correspondence between the computations of $M$ and the computations of the
modified TM. The constructed machine inherits the basic properties of $M.$ In particular,
it has the same generalized space function.

\begin{lemma}\label{M} Assume that a multi-tape DTM $M_0$ solves the word problem in a 
finitely generated semigroup or monoid $S$ with space function $S_0(n).$
Then there is a one-tape
%a one-tape, symmetric, input-input, NTM 
E-machine $M_5$ such that

(a) the equivalence relation recognized by $M_5$ is the set of all pairs of words $(u,v)$ in
the generators of $S$ satisfying the equality $u=_S v;$

(b) the generalized space function $S'_5(n)$ of $M_5$ is equivalent to $S_0(n);$

(c) the left and right parts of the tape alphabet of $M_5$ are disjoint;

(d) if $w$ is a reachable configuration of $M_5,$ then there is a computation
$w_0\to\dots\to w,$ where $w$ is an input configuration of $M_5$ and
$|w_0|_a\le |w|_a+c$ for an integer $c\ge 1$ independent of $w;$

(e) if a reduced computation $w_0\to\dots\to w_t$ of $M_5$ has no commands involving $\alpha$ or
has no commands involving $\omega,$ then $t$ is bounded from above by $c_1|w_0|+c_2$
for some constants $c_1,c_2;$

(f) if a computation $W_0=Uq_1\omega\to\dots\to W_t$ of $M_5$ has no commands involving $\alpha,$
%and starts with an input configuration, 
then $|W_0|=|W_t|$ and the commands of this computation do not involve letters from $Y_l\backslash A_l.$  
A non-empty reduced computation $Uq_1\omega\to\dots\to U'q_1\omega$ of $M_5$ has a command involving $\alpha.$

(g) let an $M_5$-computation starts with an input configuration $w_0=\alpha u q_1\omega$
and ends with $w_t=\alpha u' q_1\omega.$ Then $w_t$ is also an input configuration of $M_5.$

\end{lemma}

\proof  Recall that starting with the  DTM $M_0$ we have constructed
the input-output DTM-s $M_1,$ $M_2,$ and an
E-machine $M_3.$ Let us use the construction of this subsection assuming that $M_3=M$ and
and  $M_4=M'.$ Doubling the tape alphabet we get a machine $M_5$ providing Property (c). Then the statement (a) follows from Lemmas \ref{inin} (a),
\ref{1tape} (a), and from the definition of $M_5.$ The statement (b) follows from Lemmas \ref{M1} (b), 
%\ref{M1M2} (b), 
\ref{M2M3}, \ref{1tape} (b), and from the definition of $M_5.$ Lemma \ref{1tape} (c,d,e,f) implies Properties (d), (e), (f) and (g).
 
\endproof 

\section{Defining relations and derivation trapezia}

\subsection{Embedding homomorphism}

Now we define an embedding of $S$ in a finitely presented monoid $H.$
%announced in Theorem \ref{main}.
Let $A$ be a finite generator 
set of $S,$ and let the machine $M_5$ be given by Lemma \ref{M}.
We have $M_5= \langle A_l, Y_l \sqcup Y_r , Q,
\Theta, q_1 \rangle,$ where $A_l\subset Y_l$ is the input alphabet which is
the copy of $A$,  $Y_l \sqcup Y_r$ is the tape alphabet (with left and right parts),
$Q$ is the set of states of $M_5$, $\Theta$ is a set of commands, and $q_1\in Q$ is
the start state.

The set of generators of the monoid $H$ is $A_H=A\sqcup Y_l \sqcup Y_r\sqcup Q\sqcup\{\alpha, \omega, p\}$
where $\alpha$ and $\omega$ are the endmarker symbols of $M_5,$ and $p$ is one more generator. The set of defining relations
of $H$ is 
\begin{equation}\label{machine}
R_H=\{V'=V \;\;for\; every\; command\;\; V\to V'\; of\; M_5\}
\cup
\end{equation} 
\begin{equation} \label{aux}
\{pa=a_lp \;\; for\; every\;\; a\in A\;\; and\; for\; its\; copy \;\; a_l\in A_l\}\cup\{\alpha p = 1, p=q_1\omega\}
\end{equation}
\medskip

\begin{lemma} \label{hom} The identity map on the generator set $A$ of $S$ extends
to a homomorphism $\phi: S\to H.$ If $S$ has $1$ in the signature, then $\phi$ is a monoid
homomorphism (i.e. $\phi(1)=1$).
\end{lemma} 
\proof Assume that $u=_S v.$ 
We must prove that $u=_H v.$ 

By Lemma \ref{M} (a), there is an input-input
computation $C$ of $M_5$ starting with $\alpha u_lq_1\omega$ and ending with $\alpha v_l q_1\omega,$
where $u_l$ and $v_l$ are the copies of $u$ and $v$ in the input alphabet $A_l$ of $M_5.$
Since the relations $V=V'$ are included in $R_H$ for all the  commands $V\to V'$ of $M_5,$  all
configurations of $C$ are equal in $H,$ in particular, $\alpha u_lq_1\omega=_H \alpha v_l q_1\omega.$
Using the relation $q_1\omega = p,$ we obtain $\alpha u_l p=_H \alpha v_l p.$ Now applying
relations of the form $a_lp=pa,$ we have $\alpha p u =_H \alpha p v.$ Finally,
$u= _H v$ since $\alpha p=1$ by the definition of $H.$ \endproof

We will prove in Lemma \ref{inj} that $\phi$ is an injective homomorphism.

\begin{rk}\label{sp} It follows from the proof of Lemma \ref{hom} and from Lemma \ref{M}(b)
that for two equal in $S$ words $u$ and $v$ of length at most $n,$ 
%can be connected by a derivation over $H$ with space 
we have $space_H(u,v)\le S'_5(n)+3.$
\end{rk}

\medskip

\subsection{Derivation trapezia}

Assume that $\cal S =\langle{\cal A}\mid\cal R\rangle$ is a semigroup or monoid
presentation. Then every derivation over this presentation has a visual geometric
interpretation in terms of finite connected planar graphs. For group presentations,
these graphs are called van Kampen diagrams (see \cite{LS}), and semigroup
diagrams were introduced by Kashintsev (see \cite{K} and \cite{Re}). Below we
use a modified approach. Our diagrams uniquely restore derivations,
which is preferable when one compares derivations with the computations of a TM. We call
such diagrams {\it derivation trapezia} since they look similar to trapezia constructed from
bands and associated with group computations (see \cite{R}, \cite{SBR}, \cite{O1},
\cite{BORS}, etc.)

Every {\it cell} $\pi$ is a trapezium in  Euclidean plane with horizontal top
and bottom. The top and the bottom of a {\it trivial} cell are labeled by the same
letter from $\cal A.$ In the {\it relation} cell $\pi$ corresponding to a nontrivial relation
$u=v$ from $\cal R,$ the bottom is labeled by the word $u$ and the top is labeled
by $v.$ This means that the bottom (the top) is divided into $|u|$ (resp., $|v|$)
subsegments of nonzero length, each of the subsegments has a label from $\cal A,$
and one read the word $u$ (the word $v$) on the bottom (on the top) from left to right.
The sides of the trapezium $\pi$ have no labels. 
Note that $\pi$ can be a triangle if $|u|=0$ or $|v|=0$; but
we will not include the trivial relations of the form $1=1$ in $\cal R.$ Also
we assume that $\cal R$ is symmetric, i.e., a relation $u=v$ belongs to $\cal R$ iff $v=u$
is in $\cal R.$

For every 
transition $w'uw''\to w'vw'',$ where $u=v$ is a defining relation from $\cal R,$
we construct a {\it derivation band} as follows.
We draw a horizontal parallel paths in the  plane, the top and the bottom path
directed from left to right. The bottom path (the top path) has  $|w'uw''|$ (resp., $|w'vw''|$)
edges of nonzero length, each of them is labeled by a letter from $\cal A$  so that
the label of the bottom (the top) path is  $w'uw''$ (resp. $w'vw''$). Then we connect
the initial (the terminal) vertex of the subsegment labeled by $u$ in the bottom with,
respectively, the initial (the terminal) vertex of the subsegment labeled by $v$ in the top.
This gives us the relation cell corresponding to the relation $u=v.$ Finally, we connect
the corresponding vertices of the top and the bottom to obtain $|w'|+|w''|$ trivial
cells of the constructed derivation band. The left-most and the right-most connecting
segments are, respectively, the left and the right {\it sides} of the derivation band.

%TeXCAD (http://texcad.sf.net/) Picture. File: [semi1.pic]. Options on following lines.
%\grade{\on}
%\emlines{\off}
%\epic{\off}
%\beziermacro{\on}
%\reduce{\on}
%\snapping{\off}
%\pvinsert{% Your \input, \def, etc. here}
%\quality{8.000}
%\graddiff{0.005}
%\snapasp{1}
%\zoom{4.0000}
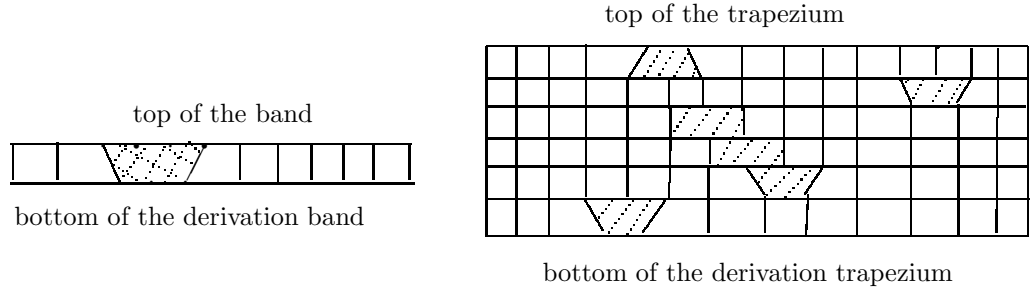
\unitlength 1mm % = 2.845pt
\linethickness{0.4pt}
\ifx\plotpoint\undefined\newsavebox{\plotpoint}\fi % GNUPLOT compatibility
\begin{picture}(100, 50)(10, 45)
\put(10,65.5){\line(1,0){53.75}}
\put(10.25,70.75){\line(0,-1){4.75}}
\put(10,70.75){\line(1,0){53.25}}
\put(63,70.75){\line(0,-1){4.75}}
\put(16.25,71){\line(0,-1){5}}
%\emline(22.25,70.75)(24.5,65.75)
\multiput(22.25,70.75)(.03358209,-.07462687){67}{\line(0,-1){.07462687}}
%\end
\put(35.75,70.5){\line(-1,-2){2.25}}
\put(22.25,70.5){\circle*{.5}}
\put(35.75,70.5){\circle*{.707}}
\put(33.25,65.75){\circle*{.5}}
\put(24.5,65.75){\circle*{.5}}
\put(28.75,65.75){\circle*{.5}}
\put(26.75,70.5){\circle*{.707}}
\put(31.25,70.5){\circle*{.5}}
\put(40.5,70.75){\line(0,-1){4.75}}
\put(45.5,70.5){\line(0,-1){5}}
\put(50,70.5){\line(0,-1){4.5}}
\put(54.25,70.75){\line(0,-1){4.75}}
\put(58.25,70.5){\line(0,-1){4.5}}
\put(123.5,205.75){\line(0,1){0}}
\put(85,204){\line(0,-1){.5}}
\put(73.25,83.75){\line(0,-1){25.25}}
\put(73.25,83.75){\line(1,0){72}}
\put(145.25,83.75){\line(0,-1){25.25}}
\put(145.25,58.5){\line(-1,0){71.75}}
\put(73.25,79.5){\line(1,0){72}}
\put(73.25,75.75){\line(1,0){71.75}}
\put(73.25,71.5){\line(1,0){72}}
\put(73.25,67.75){\line(1,0){72}}
\put(73.25,63.5){\line(1,0){72.25}}
%\emline(86.25,63.5)(88.75,59.25)
\multiput(86.25,63.5)(.03333333,-.05666667){75}{\line(0,-1){.05666667}}
%\end
%\emline(97,63.25)(94,59)
\multiput(97,63.25)(-.03370787,-.04775281){89}{\line(0,-1){.04775281}}
%\end
%\emline(107.75,67.75)(110.25,64)
\multiput(107.75,67.75)(.03333333,-.05){75}{\line(0,-1){.05}}
%\end
%\emline(118,67.75)(115.5,64)
\multiput(118,67.75)(-.03333333,-.05){75}{\line(0,-1){.05}}
%\end
\put(112.75,71.5){\line(0,-1){3.5}}
\put(103,71.5){\line(0,-1){3.25}}
\put(97.75,75.75){\line(0,-1){3.75}}
\put(107.5,75.75){\line(0,-1){3.75}}
%\emline(92.25,79.75)(94.75,83.75)
\multiput(92.25,79.75)(.03333333,.05333333){75}{\line(0,1){.05333333}}
%\end
%\emline(100,83.5)(101.75,79.75)
\multiput(100,83.5)(.03365385,-.07211538){52}{\line(0,-1){.07211538}}
%\end
\put(92,79.25){\line(0,-1){15.5}}
\put(97.5,79.5){\line(0,-1){4}}
%\emline(97.75,72)(97.5,63.5)
\multiput(97.75,72)(-.03125,-1.0625){8}{\line(0,-1){1.0625}}
%\end
\put(86.5,84){\line(0,-1){20.25}}
\put(81.5,83.75){\line(0,-1){25.25}}
\put(77.25,83.75){\line(0,-1){25}}
\put(102.75,67.75){\line(0,-1){8.75}}
\put(107.25,83.75){\line(0,-1){8.5}}
\put(110.25,63.75){\line(0,-1){4.75}}
%\emline(116,64)(115.75,58.75)
\multiput(116,64)(-.03125,-.65625){8}{\line(0,-1){.65625}}
%\end
\put(112.75,83.75){\line(0,-1){12.5}}
\put(118,83.5){\line(0,-1){15.5}}
\put(122.5,83.75){\line(0,-1){25.25}}
%\emline(133.25,83.5)(133.5,83.75)
\multiput(133.25,83.5)(.03125,.03125){8}{\line(0,1){.03125}}
%\end
\put(102,79.5){\line(0,-1){3.5}}
%\emline(128.25,79.5)(129.75,76.25)
\multiput(128.25,79.5)(.03333333,-.07222222){45}{\line(0,-1){.07222222}}
%\end
%\emline(137.75,79.5)(135.75,76.25)
\multiput(137.75,79.5)(-.03333333,-.05416667){60}{\line(0,-1){.05416667}}
%\end
\put(133,83.75){\line(0,-1){3.75}}
\put(128.25,83.75){\line(0,-1){3.75}}
\put(137.75,83.5){\line(0,-1){3.75}}
%\emline(141.25,83.5)(141,59)
\multiput(141.25,83.5)(-.03125,-3.0625){8}{\line(0,-1){3.0625}}
%\end
\put(136,76){\line(0,-1){17.5}}
\put(129.75,76.25){\line(0,-1){17.5}}
%\dottedline(90,63.25)(88.25,61)
\multiput(89.93,63.18)(-.4375,-.5625){5}{{\rule{.4pt}{.4pt}}}
%\end
%\dottedline(92.75,63.25)(89.5,59)
\multiput(92.68,63.18)(-.46429,-.60714){8}{{\rule{.4pt}{.4pt}}}
%\end
%\dottedline(95,63.25)(92,58.75)
\multiput(94.93,63.18)(-.42857,-.64286){8}{{\rule{.4pt}{.4pt}}}
%\end
%\dottedline(100.25,75.5)(98.25,72.5)
\multiput(100.18,75.43)(-.4,-.6){6}{{\rule{.4pt}{.4pt}}}
%\end
%\dottedline(98.25,72.5)(98.25,72.5)
\multiput(98.18,72.43)(0,0){3}{{\rule{.4pt}{.4pt}}}
%\end
%\dottedline(103,75.25)(100.5,72)
\multiput(102.93,75.18)(-.5,-.65){6}{{\rule{.4pt}{.4pt}}}
%\end
%\dottedline(105.75,75.5)(103,71.5)
\multiput(105.68,75.43)(-.45833,-.66667){7}{{\rule{.4pt}{.4pt}}}
%\end
%\dottedline(107.25,74)(104,68)
\multiput(107.18,73.93)(-.40625,-.75){9}{{\rule{.4pt}{.4pt}}}
%\end
%\dottedline(108.5,71.25)(106.25,68.5)
\multiput(108.43,71.18)(-.45,-.55){6}{{\rule{.4pt}{.4pt}}}
%\end
%\dottedline(111.25,71.5)(108,67.25)
\multiput(111.18,71.43)(-.46429,-.60714){8}{{\rule{.4pt}{.4pt}}}
%\end
%\dottedline(112.5,69.75)(109.5,65.75)
\multiput(112.43,69.68)(-.42857,-.57143){8}{{\rule{.4pt}{.4pt}}}
%\end
%\dottedline(114,67.5)(110.5,64)
\multiput(113.93,67.43)(-.58333,-.58333){7}{{\rule{.4pt}{.4pt}}}
%\end
%\dottedline(116.75,67.75)(113,63.75)
\multiput(116.68,67.68)(-.53571,-.57143){8}{{\rule{.4pt}{.4pt}}}
%\end
%\dottedline(96.75,83.75)(94.25,79.75)
\multiput(96.68,83.68)(-.41667,-.66667){7}{{\rule{.4pt}{.4pt}}}
%\end
%\dottedline(98.5,83.75)(96.25,80)
\multiput(98.43,83.68)(-.45,-.75){6}{{\rule{.4pt}{.4pt}}}
%\end
%\dottedline(100,82.75)(98.25,79.75)
\multiput(99.93,82.68)(-.35,-.6){6}{{\rule{.4pt}{.4pt}}}
%\end
%\dottedline(101,80.75)(100.75,79.75)
\multiput(100.93,80.68)(-.125,-.5){3}{{\rule{.4pt}{.4pt}}}
%\end
%\dottedline(131.25,79.25)(129.75,76.25)
\multiput(131.18,79.18)(-.3,-.6){6}{{\rule{.4pt}{.4pt}}}
%\end
%\dottedline(134.25,79.25)(132.25,75.75)
\multiput(134.18,79.18)(-.4,-.7){6}{{\rule{.4pt}{.4pt}}}
%\end
%\dottedline(136.5,79.25)(134.75,76.25)
\multiput(136.43,79.18)(-.35,-.6){6}{{\rule{.4pt}{.4pt}}}
%\end
%\dottedline(129.5,79.5)(129,79)
\multiput(129.43,79.43)(-.25,-.25){3}{{\rule{.4pt}{.4pt}}}
%\end
%\dottedline(25.25,70.5)(23.5,68)
\multiput(25.18,70.43)(-.4375,-.625){5}{{\rule{.4pt}{.4pt}}}
%\end
%\dottedline(28.25,70.75)(25,66.5)
\multiput(28.18,70.68)(-.46429,-.60714){8}{{\rule{.4pt}{.4pt}}}
%\end
%\dottedline(30.75,70.75)(28,66.25)
\multiput(30.68,70.68)(-.45833,-.75){7}{{\rule{.4pt}{.4pt}}}
%\end
%\dottedline(33.25,70.5)(30.5,65.75)
\multiput(33.18,70.43)(-.45833,-.79167){7}{{\rule{.4pt}{.4pt}}}
%\end
%\dottedline(23.75,68)(26.75,66)
\multiput(23.68,67.93)(.6,-.4){6}{{\rule{.4pt}{.4pt}}}
%\end
%\dottedline(24.5,70.5)(32,66.5)
\multiput(24.43,70.43)(.75,-.4){11}{{\rule{.4pt}{.4pt}}}
%\end
%\dottedline(29.75,70.75)(34,68.25)
\multiput(29.68,70.68)(.70833,-.41667){7}{{\rule{.4pt}{.4pt}}}
%\end
%\dottedline(33.75,71)(35,69.75)
\multiput(33.68,70.93)(.4167,-.4167){4}{{\rule{.4pt}{.4pt}}}
%\end
\put(26.25,73.75){top of the band}
\put(89,87){top of the trapezium}
\put(10.5,60){bottom of the derivation band}
\put(80.75,52.5){bottom of the derivation trapezium}
\end{picture}

It is obvious that every band with at most one transition cell corresponds to an 
elementary transition $w\to w'.$ (We allow trivial transitions $w\to w$. The corresponding
bands have only trivial cells. Note that the band corresponding to the transition
$1\to 1$ has no cells, but it has unlabeled side edges.)

Let $w_0\to w_1\to\dots\to w_t$ be a derivation over $\cal S.$  Then the {\it derivation
trapezium} $\Delta$ of height $t$ corresponding to this derivation is composed of $t$ derivation bands,
where the bottom of the derivation band ${\cal T}_{i+1}$ corresponding to the transition $w_i\to w_{i+1}$
coincides with the top of the derivation band ${\cal T}_i$ corresponding to $w_{i-1}\to w_i$
($i=1,\dots, t-1$). Thus the label of the bottom (of the top) of $\Delta$ is $w_0$ (resp., $w_t$).
The left (the right) sides of the derivation bands ${\cal T}_i$-s form the {\it left side}
({\it the right side}) of $\Delta.$

We see that every derivation produces a derivation trapezium, and vice versa, every
trapezium composed of derivation bands as above, is a derivation trapezium for some derivation
(which may admit trivial transitions).
Every horizontal edge of a derivation trapezium is labeled, and every vertical one (i.e., connecting
the top and the bottom of a derivation band) is unlabeled. 

A path is {\it vertical} if every its edge is vertical and different edges cross different
derivation bands.

We call a derivation trapezium $\Delta$
{\it  indivisible} 
%\footnote{Mozhno pozzhe on etom?} 
if the only vertical paths connecting the top and
the bottom of $\Delta$ are the left and the right sides of $\Delta.$

\begin{rk}\label{indiv}

If $\Delta$ corresponds to a derivation $w_0\to\dots\to w_t,$ and it is divisible, 
then $\Delta$ is a union of two derivation trapezia of the same height: $\Delta_1$
and $\Delta_2,$ where $w_0(1)$ and $w_t(1)$ are bottom and top labels of $\Delta_1,$
$w_0(2)$ and  $w_t(2)$ are bottom and top labels of $\Delta_2,$ $w_0=w_0(1)w_0(2),$ $w_t=w_t(1)w_t(2).$
The derivation trapezium $\Delta_1$ (resp., $\Delta_2$) corresponds to a derivation $w_0(1)\to\dots\to w_t(1)$ 
(to $w_0(2)\to\dots\to w_t(2)$), where some transitions may be trivial. This observation reduces the study of the properties 
of derivation trapezia to indivisible ones. 

Also we can apply the following {\it time separation trick} to the divisible derivation trapezium $\Delta.$
Since  for every $i\le t,$ either transition $w_{i-1}(1)\to w_i(1)$ or the transition
$w_{i-1}(2)\to w_i(2)$ is trivial (does not change the word), one may switch the
order of the corresponding transitions in the derivation $w_0\to\dots\to w_t$ as follows:
\begin{equation}\label{div}
w_0(1)w_0(2)\to\dots\to w_t(1)w_0(2)\to\dots\to w_t(1)w_t(2),
\end{equation}
where the length of this derivation  is $t$ (not $2t$) since some
trivial transitions are now omitted in  first and in the second parts of (\ref{div}). 
Thus, the corresponding derivation trapezium $\Delta'$ has the same height as $\Delta$ and the same
bottom and top labels.
The derivation subtrapezia
$\Delta'_1$ and $\Delta'_2$ of $\Delta'$ are {\it time separated}: the derivation bands corresponding
to the nontrivial transitions in $\Delta_2$ follow after the transitions bands corresponding
to the nontrivial transitions in $\Delta_1$ (or vice versa). 
Note that the space of Derivation (\ref{div}) can be greater that the space of 
the original derivation.

\end{rk} 

\subsection{Vertical bands in trapezia over the monoid $H$}

derivation bands are horizontal. Now we consider derivation trapezia over the
presentation \\ $H=\langle A_H\mid R_H\rangle$ and define vertical bands, namely
$q$-bands, $\alpha$-bands, $\omega$-bands and $a$-bands. 

By definition, a $q$-{\it letter} is a letter from $Q\cup\{p\}.$ A $q$-{\it edge} is an
edge labeled by a $q$-letter, a $q$-{\it cell} is a cell, having a $q$-edge in its top or
bottom. (So every $p$-{\it cell}, i.e., having a boundary edge labeled by $p,$ is also a $q$-cell.) A $q$-{\it band} of length $n$ in a derivation trapezia $\Delta$ is a sequence of $q$-cells $\pi_1,\dots,\pi_n$ such that
if a cell $\pi_i$ belongs to a derivation band ${\cal T}_j,$ then $\pi_{i+1}$ belongs
to ${\cal T}_{j+1}$ and these two cells share a $q$-edge ($i=1,\dots, n-1$).

A $q$-band $\cal C$ is called {\it maximal} if it is not contained in a longer $q$-band.
It follows from the list of defining relations of $H$ that the first cell $\pi_1$ (the last
cell $\pi_n$) of $\cal C$ either shares a $q$-edge with the bottom of $\Delta$ (resp. with
the top of $\Delta$) or it is an $\alpha p$-{\it cell}, i.e., a cell corresponding to the
relation $1=\alpha p$ (to the relation $\alpha p = 1,$ resp.).

Similarly one defines $\alpha$- and $\omega$-{\it edges}, $\alpha$-{\it bands} and $\omega$-{\it bands}. The properties of the first
and the last cells of maximal $\alpha$-bands are similar to the properties of the
maximal $q$-bands mentioned above. The first cell $\pi_1$ (the last
cell $\pi_n$) of a maximal $\omega$-band  either shares an $\omega$-edge with the bottom of $\Delta$ (resp. with
the top of $\Delta$) or it is a $q_1 \omega$-{\it cell}, i.e., a cell corresponding to the
relation $p= q_1\omega$ (to the relation $q_1\omega = p,$ resp.).

An $a$-{\it edge} is an edge labeled by a letter of the alphabet $A\cup Y_l\cup Y_l,$ 
which, by definition, consists of $a$-{\it letters}. By definition, an $a$-band consists of
trivial $a$-cells $\pi_i$-s with one $a$-letter written on the bottom and with the same letter labeling
the top of $\pi_i$-s. A maximal $a$-band must start (end) either on the bottom (resp., top) of $\Delta$ 
or on the boundary of a $q$-cell having an $a$-letter in its top (resp., bottom) label.

The above definitions imply that a cell cannot belong to two different maximal $q$-bands (resp., $\alpha$-bands, $\omega$-bands,
$a$-bands). If an $\alpha$-band $\cal B$ and a $q$-band $\cal C$ start with the same
$\alpha p$-cell $\pi,$ 
then any derivation band going from left to right and crossing both $\cal B$ and $\cal C,$
must first cross $\cal B$ and then it crosses $\cal C$ (or it crosses a cell
shared by $\cal B$ and $\cal C$).
So we may say that the band $\cal C$ is disposed from the right of $\cal B.$
This simple observation leads to

\begin{lemma} \label{alfaq} Assume that an $\alpha$-band $\cal B$ and a $q$-band $\cal C$ of a derivation trapezium $\Delta$ start or end with the same $\alpha p$-cell $\pi.$ Then either they end (resp., start) with the same
$\alpha p$-cell $\pi'$ or they both reaches the top (resp., bottom) of $\Delta.$
\end{lemma}
\proof Proving by contradiction, we assume that these two bands start with $\pi,$ and band $\cal B$
is not longer than $\cal C$. (The other cases are similar.) Then $\cal B$ is disposed
from the left of $\cal C,$ and the last cell of $\cal B$ is an $\alpha p$-cell $\pi_1$.
 Then some $q$-band ${\cal C}_1$ must also terminate at $\pi_1,$ and
${\cal C}_1$ has to be placed from the right of $\cal B$ and from the left of $\cal C.$
Since it cannot start with $\pi,$ ${\cal C}_1$ has to start with an $\alpha p $-cell
$\pi_2$ belonging to a derivation band situated above the derivation band containing
the cell $\pi.$ Similarly, an $\alpha$-band ${\cal B}_1$ must start with $\pi_2,$ it
is disposed from the left of ${\cal C}_1$ and from the right of $\cal B,$ and its
last cell is an $\alpha p$-cell $\pi_3\ne \pi_1.$ Reasoning this way, we can get arbitrarily
many cells in $\Delta,$ a contradiction. \endproof 

%TeXCAD (http://texcad.sf.net/) Picture. File: [semi2.pic]. Options on following lines.
%\grade{\on}
%\emlines{\off}
%\epic{\off}
%\beziermacro{\on}
%\reduce{\on}
%\snapping{\off}
%\pvinsert{% Your \input, \def, etc. here}
%\quality{8.000}
%\graddiff{0.005}
%\snapasp{1}
%\zoom{4.0000}
\unitlength 1mm % = 2.845pt
\linethickness{0.4pt}
\ifx\plotpoint\undefined\newsavebox{\plotpoint}\fi % GNUPLOT compatibility
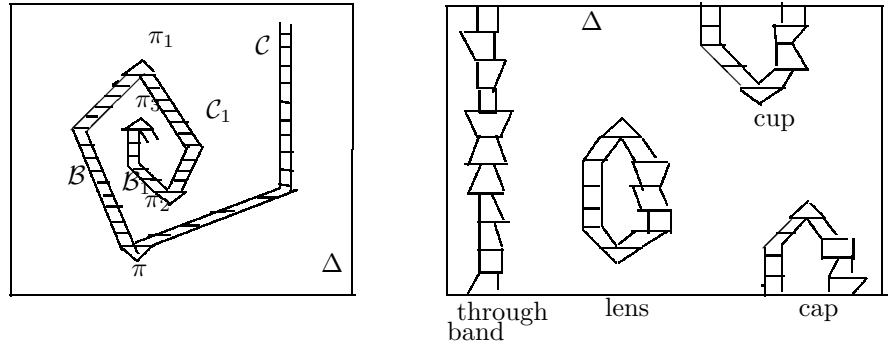
\begin{picture}(135.25,58.25)(5,12)
\put(33.5,47.25){\line(1,0){4.25}}
%\emline(33.75,47.25)(36,48.75)
\multiput(33.75,47.25)(.05,.03333333){45}{\line(1,0){.05}}
%\end
%\emline(36,48.75)(37.5,47.5)
\multiput(36,48.75)(.03947368,-.03289474){38}{\line(1,0){.03947368}}
%\end
\put(34.25,47){\line(0,-1){4.5}}
%\emline(34.25,42.5)(37.75,39)
\multiput(34.25,42.5)(.033653846,-.033653846){104}{\line(0,-1){.033653846}}
%\end
\put(37.75,39){\line(1,0){4.25}}
%\emline(38,38.75)(39.75,37.5)
\multiput(38,38.75)(.04605263,-.03289474){38}{\line(1,0){.04605263}}
%\end
\put(39.75,37.5){\line(0,-1){.25}}
%\emline(39.75,37.25)(42,39)
\multiput(39.75,37.25)(.04326923,.03365385){52}{\line(1,0){.04326923}}
%\end
\put(35.75,47.25){\line(0,-1){4.5}}
%\emline(35.75,42.75)(39.5,39.25)
\multiput(35.75,42.75)(.036057692,-.033653846){104}{\line(1,0){.036057692}}
%\end
%\emline(39.75,39.5)(42.25,45)
\multiput(39.75,39.5)(.03333333,.07333333){75}{\line(0,1){.07333333}}
%\end
\put(42.25,45){\line(1,0){.25}}
%\emline(42,44.75)(35.75,54.5)
\multiput(42,44.75)(-.033602151,.052419355){186}{\line(0,1){.052419355}}
%\end
%\emline(41.5,39.25)(44.25,45)
\multiput(41.5,39.25)(.03353659,.07012195){82}{\line(0,1){.07012195}}
%\end
%\emline(44,45)(37.75,54.5)
\multiput(44,45)(-.033602151,.051075269){186}{\line(0,1){.051075269}}
%\end
\put(38,54.5){\line(-1,0){4.5}}
%\emline(33.5,54.5)(36.25,56.5)
\multiput(33.5,54.5)(.04583333,.03333333){60}{\line(1,0){.04583333}}
%\end
%\emline(36.25,56.5)(37.75,55)
\multiput(36.25,56.5)(.03333333,-.03333333){45}{\line(0,-1){.03333333}}
%\end
%\emline(36,47)(37,45.5)
\multiput(36,47)(.0333333,-.05){30}{\line(0,-1){.05}}
%\end
%\emline(37.5,47.25)(38.25,46)
\multiput(37.5,47.25)(.0326087,-.0543478){23}{\line(0,-1){.0543478}}
%\end
\put(36.5,54.75){\line(-1,-1){7.75}}
%\emline(28.75,47)(35.5,31)
\multiput(28.75,47)(.03358209,-.07960199){201}{\line(0,-1){.07960199}}
%\end
%\emline(35.5,31)(35.75,30.75)
\multiput(35.5,31)(.03125,-.03125){8}{\line(0,-1){.03125}}
%\end
\put(34.25,54.75){\line(-1,0){.25}}
\put(34,54.5){\line(-1,-1){7.25}}
%\emline(27,47.25)(33.5,31.25)
\multiput(27,47.25)(.033678756,-.082901554){193}{\line(0,-1){.082901554}}
%\end
\put(33.25,31.75){\line(1,0){4.5}}
%\emline(33.5,31.5)(35.5,29.75)
\multiput(33.5,31.5)(.03846154,-.03365385){52}{\line(1,0){.03846154}}
%\end
%\emline(35.5,29.75)(37.5,31.75)
\multiput(35.5,29.75)(.03333333,.03333333){60}{\line(0,1){.03333333}}
%\end
%\emline(35.25,32.25)(54.5,39.5)
\multiput(35.25,32.25)(.089534884,.03372093){215}{\line(1,0){.089534884}}
%\end
\put(54.5,39.5){\line(0,1){22}}
%\emline(37.25,31.75)(56.75,39.25)
\multiput(37.25,31.75)(.087443946,.033632287){223}{\line(1,0){.087443946}}
%\end
%\emline(56.75,39.25)(55.75,39)
\multiput(56.75,39.25)(-.125,-.03125){8}{\line(-1,0){.125}}
%\end
\put(56,39){\line(0,1){22.25}}
%\emline(32.25,53.25)(33,53.75)
\multiput(32.25,53.25)(.05,.0333333){15}{\line(1,0){.05}}
%\end
\put(32.25,53){\line(1,0){2.25}}
\put(36.75,52.75){\line(1,0){1.75}}
\put(31.25,51.75){\line(1,0){2}}
\put(38,51.25){\line(1,0){1.75}}
\put(30,50.25){\line(1,0){2.25}}
\put(38.75,49.75){\line(1,0){2.25}}
\put(28.75,49){\line(1,0){1.5}}
\put(39.75,48.25){\line(1,0){2.5}}
\put(27.25,47.5){\line(1,0){1.75}}
\put(41,46.75){\line(1,0){1.75}}
\put(42,45.25){\line(1,0){2}}
\put(34.5,45.5){\line(1,0){1.5}}
\put(28,45.25){\line(1,0){1.75}}
\put(28.75,43.5){\line(1,0){1.25}}
\put(34.25,43.75){\line(1,0){1.25}}
\put(41.25,43.25){\line(1,0){2.25}}
\put(29.25,41.75){\line(1,0){1.5}}
\put(35,42){\line(1,0){2}}
\put(41,41.75){\line(1,0){1.5}}
\put(30.25,40.25){\line(1,0){1.5}}
%\emline(35.75,40.75)(38.75,40.5)
\multiput(35.75,40.75)(.375,-.03125){8}{\line(1,0){.375}}
%\end
\put(40.5,40.5){\line(1,0){1.5}}
\put(31,38.5){\line(1,0){1.25}}
\put(31.5,37){\line(1,0){2}}
\put(32.25,35){\line(1,0){1.75}}
\put(32.75,33.5){\line(1,0){1.75}}
\put(37,33.25){\line(1,0){3}}
\put(41.25,34.5){\line(1,0){2.5}}
%\emline(45.25,36.25)(48.5,36.5)
\multiput(45.25,36.25)(.40625,.03125){8}{\line(1,0){.40625}}
%\end
\put(51.5,38.25){\line(1,0){2.5}}
%\emline(54.25,40)(55.75,39.75)
\multiput(54.25,40)(.1875,-.03125){8}{\line(1,0){.1875}}
%\end
%\emline(54.5,42)(56.25,41.75)
\multiput(54.5,42)(.21875,-.03125){8}{\line(1,0){.21875}}
%\end
\put(54.75,44.25){\line(1,0){1.5}}
\put(54.75,46.5){\line(1,0){1.25}}
\put(54.5,48.5){\line(1,0){1.25}}
%\emline(54.5,51)(56,50.75)
\multiput(54.5,51)(.1875,-.03125){8}{\line(1,0){.1875}}
%\end
\put(54.5,53.5){\line(1,0){1.5}}
\put(54.5,55.75){\line(1,0){1.5}}
\put(54.5,58.25){\line(1,0){1.5}}
\put(54.5,60){\line(1,0){1.5}}
\put(34.75,27.75){$\pi$}
\put(37,58.75){$\pi_{1}$}
\put(36.5,37){$\pi_{2}$}
\put(35.25,50.25){$\pi_{3}$}
\put(51.25,57.25){$\cal C$}
\put(44.75,48.75){${\cal C}_{1}$}
\put(26.25,40){$\cal B$}
\put(33.5,39.25){${\cal B}_1$}
\put(18.75,62.25){\line(0,-1){37.25}}
\put(18.5,25.25){\line(1,0){46}}
\put(18.75,62){\line(0,1){1.75}}
\put(18.5,64){\line(1,0){45.75}}
%\emline(64.25,64)(64,25.25)
\multiput(64.25,64)(-.03125,-4.84375){8}{\line(0,-1){4.84375}}
%\end
\put(60,28.25){$\Delta$}
\put(76.5,63.75){\line(1,0){58.25}}
\put(76.75,63.75){\line(0,-1){38.5}}
\put(76.75,25.25){\line(1,0){58.5}}
\put(134.5,63.5){\line(0,-1){38}}
\put(80.75,63.75){\line(0,-1){3.5}}
\put(79,60.25){\line(1,0){4.5}}
\put(83.5,63.5){\line(0,-1){3.25}}
%\emline(78.5,60.25)(80.5,56.5)
\multiput(78.5,60.25)(.03333333,-.0625){60}{\line(0,-1){.0625}}
%\end
\put(80.5,56.5){\line(1,0){.5}}
\put(81,56.5){\line(1,0){3.5}}
\put(83,60.25){\line(0,-1){3.75}}
%\emline(84.5,56.5)(82.75,52.25)
\multiput(84.5,56.5)(-.03365385,-.08173077){52}{\line(0,-1){.08173077}}
%\end
\put(80.75,56.25){\line(0,-1){3.5}}
\put(80.75,52.75){\line(1,0){2}}
\put(83.25,52.75){\line(0,-1){3.25}}
\put(81,52.5){\line(0,-1){3}}
\put(81,49.5){\line(1,0){2.25}}
\put(79,49.75){\line(1,0){6.5}}
%\emline(79.25,49.5)(80.75,46.25)
\multiput(79.25,49.5)(.03333333,-.07222222){45}{\line(0,-1){.07222222}}
%\end
%\emline(85.25,49.5)(83.75,46.25)
\multiput(85.25,49.5)(-.03333333,-.07222222){45}{\line(0,-1){.07222222}}
%\end
\put(80.75,46.25){\line(1,0){3.25}}
%\emline(81,46.5)(79.5,42.5)
\multiput(81,46.5)(-.03333333,-.08888889){45}{\line(0,-1){.08888889}}
%\end
%\emline(83.25,46.5)(84.5,43)
\multiput(83.25,46.5)(.03289474,-.09210526){38}{\line(0,-1){.09210526}}
%\end
\put(79,42.75){\line(1,0){5.75}}
%\emline(81.25,42.75)(79.5,38.75)
\multiput(81.25,42.75)(-.03365385,-.07692308){52}{\line(0,-1){.07692308}}
%\end
%\emline(83,42.5)(84.25,39.25)
\multiput(83,42.5)(.03289474,-.08552632){38}{\line(0,-1){.08552632}}
%\end
\put(79.5,38.75){\line(1,0){4.75}}
\put(81,38.75){\line(0,-1){3.5}}
\put(81,35.25){\line(0,1){0}}
%\emline(82.75,38.5)(84.25,35.25)
\multiput(82.75,38.5)(.03333333,-.07222222){45}{\line(0,-1){.07222222}}
%\end
\put(80.75,35){\line(1,0){4.25}}
\put(81,35.25){\line(0,-1){3.75}}
\put(81,31.5){\line(0,1){0}}
%\emline(83,35)(84,31.75)
\multiput(83,35)(.0333333,-.1083333){30}{\line(0,-1){.1083333}}
%\end
\put(80.75,31.25){\line(1,0){3.5}}
\put(84,31.5){\line(0,-1){2.75}}
\put(84,28.75){\line(0,-1){.5}}
\put(81,31){\line(0,-1){2.75}}
\put(81,28.25){\line(1,0){3}}
%\emline(81,28.25)(79.5,25.5)
\multiput(81,28.25)(-.03333333,-.06111111){45}{\line(0,-1){.06111111}}
%\end
\put(83.5,28.25){\line(0,-1){2.25}}
%\emline(100,48.75)(97.25,46.25)
\multiput(100,48.75)(-.03666667,-.03333333){75}{\line(-1,0){.03666667}}
%\end
\put(97.25,46.25){\line(1,0){5.25}}
%\emline(102.5,46.25)(100,48.75)
\multiput(102.5,46.25)(-.03333333,.03333333){75}{\line(0,1){.03333333}}
%\end
%\emline(97,46.25)(94.75,43)
\multiput(97,46.25)(-.03358209,-.04850746){67}{\line(0,-1){.04850746}}
%\end
\put(94.75,43){\line(0,-1){9}}
%\emline(94.75,34)(97.5,31.25)
\multiput(94.75,34)(.03353659,-.03353659){82}{\line(0,-1){.03353659}}
%\end
%\emline(99.75,46)(97.25,42.75)
\multiput(99.75,46)(-.03333333,-.04333333){75}{\line(0,-1){.04333333}}
%\end
\put(97.25,42.75){\line(0,-1){8.25}}
%\emline(97.25,34.5)(99.25,31.75)
\multiput(97.25,34.5)(.03333333,-.04583333){60}{\line(0,-1){.04583333}}
%\end
%\emline(100,46)(101.75,43)
\multiput(100,46)(.03365385,-.05769231){52}{\line(0,-1){.05769231}}
%\end
%\emline(102.75,46)(104.5,43.5)
\multiput(102.75,46)(.03365385,-.04807692){52}{\line(0,-1){.04807692}}
%\end
\put(95,43.25){\line(1,0){3}}
\put(101.25,43){\line(1,0){4.75}}
%\emline(101.5,42.75)(102.75,39.25)
\multiput(101.5,42.75)(.03289474,-.09210526){38}{\line(0,-1){.09210526}}
%\end
\put(102.75,39.25){\line(0,1){0}}
%\emline(106.25,42.75)(105,39.5)
\multiput(106.25,42.75)(-.03289474,-.08552632){38}{\line(0,-1){.08552632}}
%\end
\put(105,39.5){\line(0,-1){.25}}
\put(94.75,39.75){\line(1,0){2.25}}
\put(102.5,39.75){\line(1,0){2.75}}
\put(101,39.75){\line(1,0){1.75}}
%\emline(101.5,39.5)(102.5,36.75)
\multiput(101.5,39.5)(.0333333,-.0916667){30}{\line(0,-1){.0916667}}
%\end
%\emline(105,39.5)(106,37.25)
\multiput(105,39.5)(.0333333,-.075){30}{\line(0,-1){.075}}
%\end
\put(95,36.75){\line(1,0){2}}
\put(102,36.5){\line(1,0){4.5}}
\put(103.5,36.25){\line(0,-1){2.75}}
\put(106.5,36.5){\line(0,-1){3}}
\put(95,34.25){\line(1,0){2}}
\put(103.25,33.75){\line(1,0){3.25}}
\put(101.25,33.75){\line(1,0){2.5}}
%\emline(101.75,33.75)(99.5,32)
\multiput(101.75,33.75)(-.04326923,-.03365385){52}{\line(-1,0){.04326923}}
%\end
%\emline(102.75,31.75)(106,33.75)
\multiput(102.75,31.75)(.05416667,.03333333){60}{\line(1,0){.05416667}}
%\end
\put(97,31.5){\line(1,0){5.75}}
%\emline(97.25,31.5)(99.25,29.5)
\multiput(97.25,31.5)(.03333333,-.03333333){60}{\line(0,-1){.03333333}}
%\end
%\emline(102.5,31.5)(99.5,29.5)
\multiput(102.5,31.5)(-.05,-.03333333){60}{\line(-1,0){.05}}
%\end
\put(115.75,52.75){\line(1,0){6}}
%\emline(115.75,52.5)(118.5,50.5)
\multiput(115.75,52.5)(.04583333,-.03333333){60}{\line(1,0){.04583333}}
%\end
%\emline(118.25,50.75)(121.25,52.75)
\multiput(118.25,50.75)(.05,.03333333){60}{\line(1,0){.05}}
%\end
\put(110.5,64.25){\line(0,-1){5.75}}
\put(110.5,58.5){\line(1,-1){5.25}}
\put(113,63.75){\line(0,-1){5}}
%\emline(113,58.75)(118.25,53.25)
\multiput(113,58.75)(.033653846,-.03525641){156}{\line(0,-1){.03525641}}
%\end
%\emline(118.25,53.5)(120.25,55.25)
\multiput(118.25,53.5)(.03846154,.03365385){52}{\line(1,0){.03846154}}
%\end
\put(120.25,55.25){\line(1,0){4.5}}
%\emline(122.75,55.25)(121.25,52.75)
\multiput(122.75,55.25)(-.03333333,-.05555556){45}{\line(0,-1){.05555556}}
%\end
\put(120.5,58.5){\line(0,-1){2.75}}
%\emline(122.75,58.5)(124.5,55.5)
\multiput(122.75,58.5)(.03365385,-.05769231){52}{\line(0,-1){.05769231}}
%\end
\put(120.5,58.75){\line(1,0){2.25}}
%\emline(122.75,59)(124.5,61)
\multiput(122.75,59)(.03365385,.03846154){52}{\line(0,1){.03846154}}
%\end
\put(124.5,61){\line(1,0){.25}}
\put(124.75,61){\line(-1,0){4.75}}
%\emline(120.25,61.25)(120.5,58.25)
\multiput(120.25,61.25)(.03125,-.375){8}{\line(0,-1){.375}}
%\end
\put(121.25,64.25){\line(0,-1){2.75}}
\put(124.25,63.75){\line(0,-1){2.75}}
\put(110.75,61){\line(1,0){2.25}}
\put(110.5,58.5){\line(1,0){3.25}}
\put(113.25,55.75){\line(1,0){3.25}}
%\emline(124.75,37.75)(122.25,35.5)
\multiput(124.75,37.75)(-.03731343,-.03358209){67}{\line(-1,0){.03731343}}
%\end
\put(122.25,35.5){\line(0,1){0}}
\put(122.25,35.5){\line(1,0){4.5}}
%\emline(125,37.5)(127,35.5)
\multiput(125,37.5)(.03333333,-.03333333){60}{\line(0,-1){.03333333}}
%\end
\put(122.5,35.25){\line(-1,-1){3.75}}
\put(119,26.5){\line(-1,0){.5}}
\put(119,31.5){\line(0,-1){6.5}}
%\emline(124.5,35.5)(121.25,32)
\multiput(124.5,35.5)(-.033505155,-.036082474){97}{\line(0,-1){.036082474}}
%\end
\put(121.25,32){\line(0,-1){6.25}}
%\emline(124.5,35.5)(126.5,32.75)
\multiput(124.5,35.5)(.03333333,-.04583333){60}{\line(0,-1){.04583333}}
%\end
\put(126.5,32.75){\line(1,0){3.5}}
\put(130,32.75){\line(0,1){0}}
%\emline(126.75,35.5)(128.5,33.25)
\multiput(126.75,35.5)(.03365385,-.04326923){52}{\line(0,-1){.04326923}}
%\end
\put(130.75,32.75){\line(0,-1){2.75}}
\put(126.5,32.5){\line(0,-1){2.5}}
\put(126.5,30.25){\line(1,0){4.25}}
%\emline(128.75,30.25)(127.5,27.75)
\multiput(128.75,30.25)(-.03289474,-.06578947){38}{\line(0,-1){.06578947}}
%\end
\put(130.5,30.5){\line(0,-1){2.25}}
\put(127.5,27.5){\line(1,0){5}}
%\emline(132.5,27.5)(130.5,25.25)
\multiput(132.5,27.5)(-.03333333,-.0375){60}{\line(0,-1){.0375}}
%\end
\put(127.5,27.75){\line(0,-1){2}}
\put(120.5,33.5){\line(1,0){2.5}}
\put(119,30.5){\line(1,0){2.25}}
\put(119,27.75){\line(1,0){2.5}}
\put(78,22){through}
\put(76.75,19.25){band}
\put(97.75,22.5){lens}
\put(123.5,22.5){cap}
\put(117.5,47.75){cup}
\put(94.5,60.5){$\Delta$}
\end{picture}

Lemma \ref{alfaq} implies that there can exist maximal $\alpha$- and $q$-bands
of  three types in a trapezium $\Delta:$

(1) The bands connecting an $\alpha$-edge (or a $q$-edge) of the bottom of $\Delta$ with an 
$\alpha$-edge (or a $q$-edge) of the top. We call such bands {\it through bands}. 

(2) Pairs formed by an $\alpha$-band and a $q$-band, sharing their  the first and the
last cells. We call such a pair an $\alpha q$-{\it lens}.

(3) Pairs formed by an $\alpha$-band and a $q$-band, sharing the  first  (the 
last) cell and terminating (resp., starting) on the $\alpha$- and $q$-edges of the
top (resp, bottom) of $\Delta.$ We say that such a pair form an $\alpha q$-{\it cup} (resp., $\alpha q$-{\it cap}).

\begin{lemma} \label{omega} Let $\pi$ and $\pi'$ be, resp., the first cell and the last cell  of a
maximal  $\omega$-band $\cal D.$ Then
%\footnote{Proverit', nuzhna li lemma} 

(a) $\pi$ and $\pi'$ 
cannot belong to different maximal $q$-bands.

(b) If $\pi$ belongs to the maximal $q$-band $\cal C$ of a $\alpha q$-lens, then $\pi'$ also belongs to $\cal C.$

(c) If an $\alpha q$-cap (or cup) $\Gamma$ surrounds no smaller caps (resp. cups), then $\Gamma$
surrounds no $\omega$-bands starting on the bottom (resp, on the top) of $\Delta.$ 
\end{lemma}
\proof {\bf (a)} Arguing by contradiction we assume that $\cal D$ is a shortest counter-example.
It starts on some maximal $q$-band $\cal C$ and ends on a maximal $q$-band $\cal C'.$
Since $\cal D$ is situated from the right of both $\cal C$ and $\cal C',$ and these two $q$-bands
do not cross, either $\cal C$ does not reach the top of the trapezium $\Delta$ or $\cal C'$ does not
start on the bottom of $\Delta.$ Choosing the former case, we deduce that
$\cal C$ ends with a $\alpha p$-cell $\pi_0,$ and the subband ${\cal C}_1$ of $\cal C$
with the first cell $\pi$ and the lats one $\pi_0$ 
%(where ${\cal C}_1$ does not include $\pi$ and $\pi_0$ themselves)
is shorter than $\cal D.$ 

%TeXCAD (http://texcad.sf.net/) Picture. File: [semi3.pic]. Options on following lines.
%\grade{\on}
%\emlines{\off}
%\epic{\off}
%\beziermacro{\on}
%\reduce{\on}
%\snapping{\off}
%\pvinsert{% Your \input, \def, etc. here}
%\quality{8.000}
%\graddiff{0.005}
%\snapasp{1}
%\zoom{4.0000}
\unitlength 1mm % = 2.845pt
\linethickness{0.4pt}
\ifx\plotpoint\undefined\newsavebox{\plotpoint}\fi % GNUPLOT compatibility
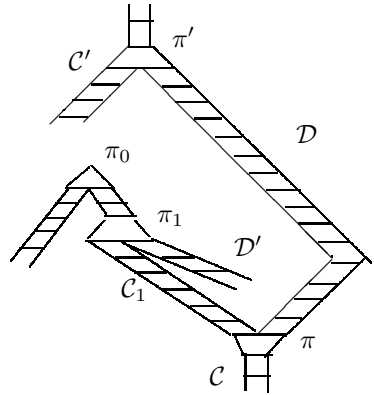
\begin{picture}(107.25,60.5)(15,10)
\put(75,68.5){\line(0,-1){5.75}}
\put(75,62.75){\line(1,0){3.25}}
\put(77.75,68.5){\line(0,-1){5.5}}
\put(75,63){\line(-1,-1){10.5}}
%\emline(78.25,62.75)(107.25,34)
\multiput(78.25,62.75)(.03399765533,-.0337045721){853}{\line(1,0){.03399765533}}
%\end
\put(72,60){\line(1,0){9}}
\put(76.5,60){\line(-1,0){.25}}
\put(76,59.75){\line(-1,-1){7}}
%\emline(106.5,34.75)(95.75,24.25)
\multiput(106.5,34.75)(-.0344551282,-.0336538462){312}{\line(-1,0){.0344551282}}
%\end
\put(88.75,24.5){\line(1,0){7.25}}
%\emline(95.75,24.5)(93,21.5)
\multiput(95.75,24.5)(-.03353659,-.03658537){82}{\line(0,-1){.03658537}}
%\end
%\emline(89,24.5)(90.25,21.75)
\multiput(89,24.5)(.03289474,-.07236842){38}{\line(0,-1){.07236842}}
%\end
\put(90.25,21.75){\line(1,0){3.75}}
\put(93.5,21.75){\line(0,-1){5}}
\put(90,21.75){\line(1,0){.5}}
\put(90.5,22){\line(0,-1){5}}
\put(76.75,60){\line(1,-1){25.25}}
\put(102,34.75){\line(0,1){.25}}
\put(92.25,24.75){\line(1,1){9.75}}
%\emline(88.5,24.5)(69.5,37)
\multiput(88.5,24.5)(-.051212938,.0336927224){371}{\line(-1,0){.051212938}}
%\end
\put(69.25,37){\line(1,0){9}}
%\emline(74,36.75)(89.25,30.5)
\multiput(74,36.75)(.081989247,-.033602151){186}{\line(1,0){.081989247}}
%\end
%\emline(78,37.25)(91.25,31.75)
\multiput(78,37.25)(.080792683,-.033536585){164}{\line(1,0){.080792683}}
%\end
%\emline(74.25,36.75)(91.75,25)
\multiput(74.25,36.75)(.0501432665,-.0336676218){349}{\line(1,0){.0501432665}}
%\end
%\emline(71.75,39.75)(69.5,37.25)
\multiput(71.75,39.75)(-.03358209,-.03731343){67}{\line(0,-1){.03731343}}
%\end
%\emline(75.75,39.75)(77.75,37.5)
\multiput(75.75,39.75)(.03333333,-.0375){60}{\line(0,-1){.0375}}
%\end
\put(71.75,40.25){\line(1,0){4.25}}
%\emline(72.75,44)(75.5,40.5)
\multiput(72.75,44)(.03353659,-.04268293){82}{\line(0,-1){.04268293}}
%\end
%\emline(69.75,43.75)(72,40.5)
\multiput(69.75,43.75)(.03358209,-.04850746){67}{\line(0,-1){.04850746}}
%\end
\put(66.5,44){\line(1,0){6.5}}
%\emline(69.75,46.75)(66.75,44.25)
\multiput(69.75,46.75)(-.04,-.03333333){75}{\line(-1,0){.04}}
%\end
%\emline(70,47)(73,44)
\multiput(70,47)(.03370787,-.03370787){89}{\line(0,-1){.03370787}}
%\end
%\emline(69.75,44)(61.75,33.5)
\multiput(69.75,44)(-.033613445,-.044117647){238}{\line(0,-1){.044117647}}
%\end
%\emline(67,44)(59.25,34.25)
\multiput(67,44)(-.033695652,-.042391304){230}{\line(0,-1){.042391304}}
%\end
\put(97.25,50){$\cal D$}
\put(85.75,17.75){$\cal C$}
\put(97.75,22.75){$\pi$}
\put(78.75,39.5){$\pi_1$}
\put(89,35){$\cal D'$}
\put(74,29.5){${\cal C}_1$}
\put(67,59.75){$\cal C'$}
\put(80.75,62.75){$\pi'$}
\put(71.75,48){$\pi_0$}
\put(75,66.25){\line(1,0){3.25}}
\put(70,57.75){\line(1,0){3.75}}
\put(79.5,57.5){\line(1,0){3.75}}
\put(67.5,55.75){\line(1,0){4}}
\put(81.25,55.5){\line(1,0){4.25}}
%\emline(66.25,54.25)(66.5,54.5)
\multiput(66.25,54.25)(.03125,.03125){8}{\line(0,1){.03125}}
%\end
\put(65.75,53.5){\line(1,0){3.75}}
\put(83.75,53.25){\line(1,0){3.5}}
\put(86.5,50.75){\line(1,0){3.5}}
\put(89,48.25){\line(1,0){3.75}}
\put(91,46){\line(1,0){3.75}}
\put(93.25,43.75){\line(1,0){4.5}}
\put(95.75,41.5){\line(1,0){4}}
\put(70.75,42){\line(1,0){4}}
\put(65.25,42){\line(1,0){2.75}}
\put(63.75,40){\line(1,0){3}}
\put(62.5,38.25){\line(1,0){3}}
\put(60.5,36){\line(1,0){3.25}}
\put(97.5,39.25){\line(1,0){4.5}}
\put(100.25,36.75){\line(1,0){4}}
\put(73.25,34.75){\line(1,0){3.75}}
\put(79.5,34.75){\line(1,0){4.25}}
\put(75.75,32.75){\line(1,0){4.25}}
\put(83.75,33){\line(1,0){4.75}}
\put(102,34){\line(1,0){3.75}}
\put(79.5,30.5){\line(1,0){4.25}}
\put(99,31.25){\line(1,0){3.75}}
\put(83.25,28.25){\line(1,0){3.5}}
\put(97,29.25){\line(1,0){3.75}}
\put(85.75,26.25){\line(1,0){4.25}}
\put(94.5,26.75){\line(1,0){3.5}}
\put(90.25,19){\line(1,0){3.5}}
\end{picture}

Since $\pi$ has a top edge  labeled by $q_1$
%${\cal C}_1$ starts with an edge labeled by $q_1$ 
and $\pi_0$ has a bottom edge labeled by $p,$ 
%ends with an edge labeled by $p,$
there must be a cell in ${\cal C}_1$ which corresponds to the relation $q_1\omega= p.$
Moreover, the number of such cells in ${\cal C}_1\backslash\pi$ must be greater than the number of cells
corresponding to the transition $p\to q_1 \omega.$  Therefore there is a cell in ${\cal C}_1,$ say $\pi_1,$ such
that a maximal $\omega$-band $\cal D'$ ends with $\pi_1$ but it does not start on ${\cal C}_1.$
This $\omega$-band $\cal D'$ is situated from the right of $\cal C$ and from the left of $\cal D.$
Since the bands $\cal C$ and $\cal D$ have the common cell $\pi$, the band $\cal D'$ is
shorter than ${\cal C}_1$, and consequently, it is shorter than $\cal D.$

We come to a contradiction with the choice of $\cal D,$ and Claim (a) is proved.

{\bf (b)} The assumption that $\cal D$ starts on $\cal C$ and ends on the top of $\Delta$
provides us, as in the proof of (a), with an $\omega$-band $\cal D'$ connecting two different maximal
$q$-bands. Thus Property (b) is proved by contradiction.

{\bf (c)} Follows from (b) since an $\omega$-band cannot start and end on the bottom (resp., on the top) of $\Delta.$
\endproof

\subsection{Minimal trapezia}\label{min}
 
If a $q$-band $\cal C$ of a derivation trapezium $\Delta$ over $H$ has $k$ $q_1\omega$-cells, then we say that
$\cal C$ has type $k.$ 
Suppose $\Delta$ has $\tau_i$ through  $q$-bands of type $i,$ $\sigma_i$ maximal $q$-bands of type $i$ in the $\alpha q$-caps and $\alpha q$-cups, and $\rho_i$ maximal $q$-bands of type $i$ in the $\alpha q$-lenses, $i=0,..,k$,  and $\Delta$ has no $q$-bands of  types $>k.$ Then we say that $\Delta$ is a trapezium
of type $\tau(\Delta)=(\tau_0,\sigma_0,\rho_0,\dots,\tau_k,\sigma_k,\rho_k,0,0,0\dots).$

Assume that $\tau(\Delta')=(\tau'_0,\sigma'_0,\rho'_0,\dots,\tau'_{k'},\sigma'_{k'}, \rho'_{k'}, 0, 0,0,\dots).$ Then by definition $\tau(\Delta)>\tau(\Delta')$
if there is $l$ such that $\tau_l>\tau'_l$ or $\tau_l=\tau'_l$ and $\sigma_l>\sigma'_l,$ or
$\tau_l=\tau'_l$ and $\sigma_l=\sigma'_l,$ but $\rho_l>\rho'_l,$
and $\tau_m=\tau'_m,$ $\sigma_m=\sigma'_m,$ $\rho_m=\rho'_m$ 
for every $m\ge l.$

Clearly, the defined order on  derivation trapezia over $H$ satisfies
the descending chain condition, and so there is a trapezium having
the smallest type among all trapezia with the same bottom an top labels.
Such a derivation trapezium is called a {\it minimal trapezium}. 

\begin{rk}\label{mini} It is easy to see that the time separation trick from Remark \ref{indiv}
preserves the numbers of  $\alpha p$-cells and   $q_1\omega$-cells in every maximal $q$-band, and so it does not change the types of maximal $q$-bands. Therefore 
it preserves the minimality of a trapezium. The same is true if one rebuilds two derivation band of a derivation trapezium replacing a subderivation $w\to w\to w'$ by $w\to w'\to w'$ or vice versa.  
\end{rk}

\begin{lemma}\label{unlab} Assume that ${\bf p}$ is a simple closed path in a minimal
trapezium $\Delta,$ and every edge of ${\bf p}$ is unlabeled. Then 
the closed region $O$ of $\Delta$ bounded by ${\bf p}$ contains no $q$-edges.
\end{lemma}
\proof 
Let us shrink to a point every labeled (horizontal) edge which is inside
$O.$ If after this surgery some unlabeled (vertical) edges connect the same vertices, 
we identify such edges. 

%TeXCAD (http://texcad.sf.net/) Picture. File: [semi4.pic]. Options on following lines.
%\grade{\on}
%\emlines{\off}
%\epic{\off}
%\beziermacro{\on}
%\reduce{\on}
%\snapping{\off}
%\pvinsert{% Your \input, \def, etc. here}
%\quality{8.000}
%\graddiff{0.005}
%\snapasp{1}
%\zoom{4.0000}
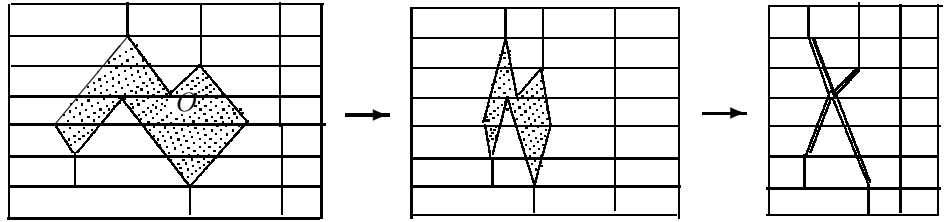
\unitlength 1mm % = 2.845pt
\linethickness{0.4pt}
\ifx\plotpoint\undefined\newsavebox{\plotpoint}\fi % GNUPLOT compatibility
\begin{picture}(145.5,40.75)(10,21)
\put(21.5,56.5){\line(0,-1){28.5}}
\put(21.75,56.5){\line(1,0){41.5}}
\put(63,56.5){\line(0,-1){28.5}}
\put(21.25,28){\line(1,0){41.5}}
\put(21.5,44.25){\line(1,0){41.5}}
\put(21.5,36.25){\line(1,0){41.5}}
\put(21.5,52.25){\line(1,0){41.25}}
\put(21.75,48.25){\line(1,0){41.25}}
\put(21.5,40.5){\line(1,0){42}}
\put(21.5,32.25){\line(1,0){41.5}}
\put(27.5,40.5){\line(5,6){10}}
\put(46,41){\line(0,1){0}}
%\emline(37.25,52.25)(43,44.5)
\multiput(37.25,52.25)(.033625731,-.045321637){171}{\line(0,-1){.045321637}}
%\end
%\emline(43,44.75)(46.75,48.25)
\multiput(43,44.75)(.036057692,.033653846){104}{\line(1,0){.036057692}}
%\end
\put(46.75,48.25){\line(0,1){.25}}
%\emline(46.75,48.5)(53.25,40.75)
\multiput(46.75,48.5)(.033678756,-.04015544){193}{\line(0,-1){.04015544}}
%\end
%\emline(53.25,41)(45.5,32.25)
\multiput(53.25,41)(-.033695652,-.038043478){230}{\line(0,-1){.038043478}}
%\end
%\emline(45.5,32.25)(36.25,44.25)
\multiput(45.5,32.25)(-.0336363636,.0436363636){275}{\line(0,1){.0436363636}}
%\end
%\emline(36.5,44)(30.5,36.75)
\multiput(36.5,44)(-.033707865,-.040730337){178}{\line(0,-1){.040730337}}
%\end
%\emline(27.75,40.25)(30,36.75)
\multiput(27.75,40.25)(.03358209,-.05223881){67}{\line(0,-1){.05223881}}
%\end
\put(57.75,40.5){\line(0,-1){4}}
\put(57.75,37.25){\line(0,-1){8.75}}
\put(57.5,52.5){\line(0,-1){12.25}}
\put(57.5,56.75){\line(0,-1){4.25}}
\put(47,48.5){\line(0,1){8}}
\put(45.5,32.25){\line(0,-1){3.75}}
\put(37.25,56.75){\line(0,-1){4.25}}
\put(30.25,36.75){\line(0,-1){4.25}}
\put(75,56){\line(1,0){35.5}}
\put(110.5,56){\line(0,-1){28}}
\put(75,56){\line(0,-1){28}}
\put(75,28.5){\line(1,0){35.25}}
\put(75.25,52){\line(1,0){35.25}}
\put(75.25,48){\line(1,0){35.25}}
\put(75,44){\line(1,0){35.5}}
\put(75,40.25){\line(1,0){35.5}}
\put(75,36){\line(1,0){35.5}}
\put(75,32.25){\line(1,0){35.75}}
%\emline(87.5,52)(84.5,40.75)
\multiput(87.5,52)(-.03370787,-.12640449){89}{\line(0,-1){.12640449}}
%\end
%\emline(87.5,52)(89,44.25)
\multiput(87.5,52)(.03333333,-.17222222){45}{\line(0,-1){.17222222}}
%\end
%\emline(89.25,44.5)(92.25,48)
\multiput(89.25,44.5)(.03370787,.03932584){89}{\line(0,1){.03932584}}
%\end
%\emline(92.25,47.75)(93.5,40.5)
\multiput(92.25,47.75)(.03289474,-.19078947){38}{\line(0,-1){.19078947}}
%\end
%\emline(93.5,40.5)(91.25,32.25)
\multiput(93.5,40.5)(-.03358209,-.12313433){67}{\line(0,-1){.12313433}}
%\end
%\emline(87.75,44)(91.25,32.75)
\multiput(87.75,44)(.033653846,-.108173077){104}{\line(0,-1){.108173077}}
%\end
%\emline(87.75,44.25)(85.75,36.5)
\multiput(87.75,44.25)(-.03333333,-.12916667){60}{\line(0,-1){.12916667}}
%\end
%\emline(84.75,40.5)(85.5,36.25)
\multiput(84.75,40.5)(.0326087,-.1847826){23}{\line(0,-1){.1847826}}
%\end
\put(87.5,56){\line(0,-1){3.75}}
\put(92.5,55.75){\line(0,-1){8}}
\put(102,56){\line(0,-1){27}}
\put(91.25,32.25){\line(0,-1){3.5}}
\put(85.75,36){\line(0,-1){3.5}}
\put(122.5,56.25){\line(0,-1){27.75}}
\put(122.5,56.25){\line(1,0){23}}
\put(145,56.5){\line(0,-1){28}}
\put(122.25,28.5){\line(1,0){22.75}}
\put(122.5,52){\line(1,0){22.5}}
\put(122.5,48){\line(1,0){22.5}}
\put(122.75,44){\line(1,0){22.25}}
\put(122.75,40.25){\line(1,0){22.5}}
\put(122.75,36.25){\line(1,0){22.25}}
\put(122.25,32){\line(1,0){23}}
%\emline(127.75,52.25)(130.5,44.5)
\multiput(127.75,52.25)(.03353659,-.0945122){82}{\line(0,-1){.0945122}}
%\end
%\emline(130.5,44.5)(127.25,36)
\multiput(130.5,44.5)(-.033505155,-.087628866){97}{\line(0,-1){.087628866}}
%\end
%\emline(130.75,44)(128,36.75)
\multiput(130.75,44)(-.03353659,-.08841463){82}{\line(0,-1){.08841463}}
%\end
\put(130.75,45.25){\line(0,1){.25}}
%\emline(130.5,44.75)(130.75,45)
\multiput(130.5,44.75)(.03125,.03125){8}{\line(0,1){.03125}}
%\end
\put(131,45){\line(0,1){.25}}
\put(131,44.75){\line(0,-1){.25}}
%\emline(130.75,44.75)(131,45)
\multiput(130.75,44.75)(.03125,.03125){8}{\line(0,1){.03125}}
%\end
%\emline(131.25,45)(136,32.75)
\multiput(131.25,45)(.033687943,-.086879433){141}{\line(0,-1){.086879433}}
%\end
%\emline(131,44)(135.75,32.5)
\multiput(131,44)(.033687943,-.081560284){141}{\line(0,-1){.081560284}}
%\end
\put(127.75,56.25){\line(0,-1){4.25}}
\put(134.5,56.75){\line(0,-1){9.25}}
\put(127.25,36.5){\line(0,-1){4.5}}
\put(135.75,32.5){\line(0,-1){4}}
\put(140,56.5){\line(0,-1){27.75}}
%\dottedline(34.75,48.75)(38.5,51)
\multiput(34.68,48.68)(.75,.45){6}{{\rule{.4pt}{.4pt}}}
%\end
%\dottedline(32,45.5)(38.75,50)
\multiput(31.93,45.43)(.75,.5){10}{{\rule{.4pt}{.4pt}}}
%\end
%\dottedline(29.75,43.25)(39.25,49)
\multiput(29.68,43.18)(.79167,.47917){13}{{\rule{.4pt}{.4pt}}}
%\end
%\dottedline(28.75,41.75)(40.25,48.5)
\multiput(28.68,41.68)(.82143,.48214){15}{{\rule{.4pt}{.4pt}}}
%\end
%\dottedline(28.5,40.5)(40.5,47.75)
\multiput(28.43,40.43)(.8,.48333){16}{{\rule{.4pt}{.4pt}}}
%\end
%\dottedline(28.75,39.5)(41,47)
\multiput(28.68,39.43)(.81667,.5){16}{{\rule{.4pt}{.4pt}}}
%\end
%\dottedline(29.5,38.5)(34.75,42.5)
\multiput(29.43,38.43)(.65625,.5){9}{{\rule{.4pt}{.4pt}}}
%\end
%\dottedline(34.75,42.5)(41.25,46.5)
\multiput(34.68,42.43)(.72222,.44444){10}{{\rule{.4pt}{.4pt}}}
%\end
%\dottedline(29.75,37.75)(32.5,39.75)
\multiput(29.68,37.68)(.55,.4){6}{{\rule{.4pt}{.4pt}}}
%\end
%\dottedline(37.75,43.25)(41.5,45.75)
\multiput(37.68,43.18)(.75,.5){6}{{\rule{.4pt}{.4pt}}}
%\end
%\dottedline(38,42.75)(42,45.25)
\multiput(37.93,42.68)(.66667,.41667){7}{{\rule{.4pt}{.4pt}}}
%\end
%\dottedline(38.75,41.75)(43,45)
\multiput(38.68,41.68)(.60714,.46429){8}{{\rule{.4pt}{.4pt}}}
%\end
%\dottedline(39,41.25)(47.25,47.25)
\multiput(38.93,41.18)(.6875,.5){13}{{\rule{.4pt}{.4pt}}}
%\end
%\dottedline(39.5,40.25)(42.25,43)
\multiput(39.43,40.18)(.55,.55){6}{{\rule{.4pt}{.4pt}}}
%\end
%\dottedline(45.75,45.25)(47.75,46.75)
\multiput(45.68,45.18)(.5,.375){5}{{\rule{.4pt}{.4pt}}}
%\end
%\dottedline(40.25,39.75)(43,42)
\multiput(40.18,39.68)(.55,.45){6}{{\rule{.4pt}{.4pt}}}
%\end
%\dottedline(46.25,44.5)(48.25,45.75)
\multiput(46.18,44.43)(.5,.3125){5}{{\rule{.4pt}{.4pt}}}
%\end
%\dottedline(40.75,39)(44.5,41.75)
\multiput(40.68,38.93)(.75,.55){6}{{\rule{.4pt}{.4pt}}}
%\end
%\dottedline(46.25,43)(48.75,45.25)
\multiput(46.18,42.93)(.5,.45){6}{{\rule{.4pt}{.4pt}}}
%\end
%\dottedline(41.25,38.5)(49.75,44.5)
\multiput(41.18,38.43)(.70833,.5){13}{{\rule{.4pt}{.4pt}}}
%\end
%\dottedline(41.5,37.75)(50.5,44)
\multiput(41.43,37.68)(.75,.52083){13}{{\rule{.4pt}{.4pt}}}
%\end
%\dottedline(41.75,37)(50.75,43)
\multiput(41.68,36.93)(.75,.5){13}{{\rule{.4pt}{.4pt}}}
%\end
%\dottedline(42.75,36.75)(51.25,42.5)
\multiput(42.68,36.68)(.77273,.52273){12}{{\rule{.4pt}{.4pt}}}
%\end
%\dottedline(43,35.75)(52.25,42.5)
\multiput(42.93,35.68)(.77083,.5625){13}{{\rule{.4pt}{.4pt}}}
%\end
%\dottedline(44,35.25)(52.5,42)
\multiput(43.93,35.18)(.70833,.5625){13}{{\rule{.4pt}{.4pt}}}
%\end
%\dottedline(44.25,34.75)(52,40.25)
\multiput(44.18,34.68)(.775,.55){11}{{\rule{.4pt}{.4pt}}}
%\end
%\dottedline(44.5,34)(50.25,38)
\multiput(44.43,33.93)(.71875,.5){9}{{\rule{.4pt}{.4pt}}}
%\end
%\dottedline(45.25,33.5)(48,35.5)
\multiput(45.18,33.43)(.55,.4){6}{{\rule{.4pt}{.4pt}}}
%\end
%\dottedline(86.75,49)(88,49.75)
\multiput(86.68,48.93)(.625,.375){3}{{\rule{.4pt}{.4pt}}}
%\end
%\dottedline(86.75,49.75)(88,50.25)
\multiput(86.68,49.68)(.625,.25){3}{{\rule{.4pt}{.4pt}}}
%\end
%\dottedline(86.25,47.25)(88,48.5)
\multiput(86.18,47.18)(.5833,.4167){4}{{\rule{.4pt}{.4pt}}}
%\end
%\dottedline(85.75,45.75)(88.25,48)
\multiput(85.68,45.68)(.5,.45){6}{{\rule{.4pt}{.4pt}}}
%\end
%\dottedline(85.5,44.25)(88.5,47)
\multiput(85.43,44.18)(.6,.55){6}{{\rule{.4pt}{.4pt}}}
%\end
%\dottedline(85.25,43)(88.75,45.75)
\multiput(85.18,42.93)(.7,.55){6}{{\rule{.4pt}{.4pt}}}
%\end
%\dottedline(85,41.5)(87.5,43.5)
\multiput(84.93,41.43)(.5,.4){6}{{\rule{.4pt}{.4pt}}}
%\end
%\dottedline(84.75,41)(85.25,41)
\multiput(84.68,40.93)(.25,0){3}{{\rule{.4pt}{.4pt}}}
%\end
%\dottedline(85.25,41)(87.75,42.5)
\multiput(85.18,40.93)(.625,.375){5}{{\rule{.4pt}{.4pt}}}
%\end
%\dottedline(85,39.5)(86.75,40.5)
\multiput(84.93,39.43)(.5833,.3333){4}{{\rule{.4pt}{.4pt}}}
%\end
%\dottedline(85.5,39)(86,39.75)
\multiput(85.43,38.93)(.25,.375){3}{{\rule{.4pt}{.4pt}}}
%\end
%\dottedline(85.25,38.25)(86.25,39)
\multiput(85.18,38.18)(.5,.375){3}{{\rule{.4pt}{.4pt}}}
%\end
%\dottedline(88.5,43.25)(92.25,46.25)
\multiput(88.43,43.18)(.625,.5){7}{{\rule{.4pt}{.4pt}}}
%\end
%\dottedline(88.75,42.25)(92.25,45.25)
\multiput(88.68,42.18)(.58333,.5){7}{{\rule{.4pt}{.4pt}}}
%\end
%\dottedline(89,41.25)(92.5,44.25)
\multiput(88.93,41.18)(.58333,.5){7}{{\rule{.4pt}{.4pt}}}
%\end
%\dottedline(89.5,40.25)(92.75,43.5)
\multiput(89.43,40.18)(.54167,.54167){7}{{\rule{.4pt}{.4pt}}}
%\end
%\dottedline(89.75,39.5)(93,42.5)
\multiput(89.68,39.43)(.54167,.5){7}{{\rule{.4pt}{.4pt}}}
%\end
%\dottedline(90,38.5)(93,41.25)
\multiput(89.93,38.43)(.6,.55){6}{{\rule{.4pt}{.4pt}}}
%\end
%\dottedline(90.25,37.75)(93.5,40.25)
\multiput(90.18,37.68)(.65,.5){6}{{\rule{.4pt}{.4pt}}}
%\end
%\dottedline(90.75,37.25)(92.75,39)
\multiput(90.68,37.18)(.5,.4375){5}{{\rule{.4pt}{.4pt}}}
%\end
%\dottedline(90.75,36.25)(92.5,37.75)
\multiput(90.68,36.18)(.5833,.5){4}{{\rule{.4pt}{.4pt}}}
%\end
%\dottedline(91,35.25)(92.25,36.25)
\multiput(90.93,35.18)(.4167,.3333){4}{{\rule{.4pt}{.4pt}}}
%\end
%\dottedline(91.5,34.5)(92.25,35)
\multiput(91.43,34.43)(.375,.25){3}{{\rule{.4pt}{.4pt}}}
%\end
\thicklines
\put(66.25,41.75){\vector(1,0){5.75}}
\put(113.75,42){\vector(1,0){5.75}}
\put(43.75,42.25){$O$}
%\emline(131.25,45)(134.25,48)
\multiput(131.25,45)(.03370787,.03370787){89}{\line(0,1){.03370787}}
%\end
\put(131.25,44.5){\line(0,1){.25}}
%\emline(131.25,44.25)(134.5,47.75)
\multiput(131.25,44.25)(.033505155,.036082474){97}{\line(0,1){.036082474}}
%\end
%\emline(128.5,51.75)(131.25,44.5)
\multiput(128.5,51.75)(.03353659,-.08841463){82}{\line(0,-1){.08841463}}
%\end
\end{picture}

It is clear that we replace every derivation band of $\Delta$
by a derivation band of the obtained trapezium $\Delta'$ (but the cells belonging
to $O$ are removed), $\Delta'$ has the same top and bottom labels as $\Delta,$ and  $\tau(\Delta')<\tau(\Delta)$ if $O$ has at least
one $q$-edge (and therefore contains a maximal $q$-band). Since $\Delta$ is a minimal trapezium, the lemma is proved. \endproof

%We say that a closed region $O$ in a derivation trapezium is {\it generated} by a thick lens $E$
%if (1) $O$ contains $E$; (2) if $O$ contains  an edge labeled by a letter from $A$ of a cell $\pi,$ %then $O$ contains $\pi;$ (3) if $O$ contains  an edge from the outer boundary of some thick lens $E'$,
%then $O$ contains $E'$ (4) $O$ is minimal with respect to (1)--(3). 

%\begin{lemma}\label{O} Let a region $O$ of a derivation trapezium is {\it generated} by a lens $E.$
%Then every edge in the outer boundary component of $O$ is either unlabeled or has a label from $A.$ 
%\end{lemma} 
%\proof 
If a maximal $\omega$-band $\cal D$ starts on the right side of the maximal $q$-band $\cal C$ of a $\alpha q$-lens $E$ then $D$ also terminates on $\cal C$ by Lemma \ref{omega} (b). Let us attach all such maximal $\omega$-bands to $\cal C$
and call the obtained figure $\Gamma$ a {\it thick lens}.

%TeXCAD (http://texcad.sf.net/) Picture. File: [semi5.pic]. Options on following lines.
%\grade{\on}
%\emlines{\off}
%\epic{\off}
%\beziermacro{\on}
%\reduce{\on}
%\snapping{\off}
%\pvinsert{% Your \input, \def, etc. here}
%\quality{8.000}
%\graddiff{0.005}
%\snapasp{1}
%\zoom{4.0000}
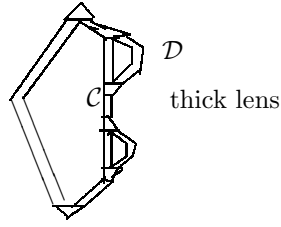
\unitlength 1mm % = 2.845pt
\linethickness{0.4pt}
\ifx\plotpoint\undefined\newsavebox{\plotpoint}\fi % GNUPLOT compatibility
\begin{picture}(63.25,45.25)(-10,30)
\put(49.25,66){\line(1,0){3.75}}
\put(54.75,63.25){\line(0,1){.25}}
\put(54.75,63.5){\line(1,0){1.75}}
\put(54.5,63.5){\line(0,-1){6.25}}
\put(54.75,57.5){\line(1,0){2.25}}
\put(55.75,63.25){\line(0,-1){5.5}}
%\emline(51.25,65.75)(54.75,63.5)
\multiput(51.25,65.75)(.05223881,-.03358209){67}{\line(1,0){.05223881}}
%\end
%\emline(52.75,65.75)(54.75,64.75)
\multiput(52.75,65.75)(.0666667,-.0333333){30}{\line(1,0){.0666667}}
%\end
%\emline(53.75,64.75)(52.5,65)
\multiput(53.75,64.75)(-.15625,.03125){8}{\line(-1,0){.15625}}
%\end
\put(54.5,65){\line(1,0){.25}}
%\emline(54,65.25)(58,64)
\multiput(54,65.25)(.10526316,-.03289474){38}{\line(1,0){.10526316}}
%\end
%\emline(55.75,63.75)(57.75,64)
\multiput(55.75,63.75)(.25,.03125){8}{\line(1,0){.25}}
%\end
%\emline(55.5,63.5)(58.25,62.25)
\multiput(55.5,63.5)(.07236842,-.03289474){38}{\line(1,0){.07236842}}
%\end
\put(58.25,62.25){\line(0,-1){2.25}}
%\emline(57.5,63.75)(59.75,62.25)
\multiput(57.5,63.75)(.05,-.03333333){45}{\line(1,0){.05}}
%\end
%\emline(59.75,62.25)(59.25,58.75)
\multiput(59.75,62.25)(-.0333333,-.2333333){15}{\line(0,-1){.2333333}}
%\end
%\emline(58,60)(55.5,57.5)
\multiput(58,60)(-.03333333,-.03333333){75}{\line(0,-1){.03333333}}
%\end
%\emline(59.25,59.25)(56.75,57.25)
\multiput(59.25,59.25)(-.04166667,-.03333333){60}{\line(-1,0){.04166667}}
%\end
\put(56.75,57.25){\line(0,1){0}}
%\emline(56.75,57.25)(55.75,56)
\multiput(56.75,57.25)(-.0333333,-.0416667){30}{\line(0,-1){.0416667}}
%\end
\put(55.75,56){\line(-1,0){1}}
\put(54.5,57.75){\line(0,-1){1.5}}
\put(54.5,56.25){\line(0,-1){3.75}}
\put(55.5,56){\line(0,-1){3}}
\put(55.5,53){\line(-1,0){1.25}}
\put(54.5,53){\line(0,-1){2}}
%\emline(55,53.25)(56.25,51.5)
\multiput(55,53.25)(.03289474,-.04605263){38}{\line(0,-1){.04605263}}
%\end
\put(54,51.25){\line(1,0){2.5}}
\put(54.5,51.25){\line(0,-1){5.5}}
%\emline(55,50.75)(55.5,51)
\multiput(55,50.75)(.0625,.03125){8}{\line(1,0){.0625}}
%\end
\put(55.5,51){\line(-1,0){.25}}
\put(55.5,51.25){\line(0,-1){5}}
\put(54.25,46.25){\line(1,0){2.75}}
%\emline(55.5,51)(57.5,49.5)
\multiput(55.5,51)(.04444444,-.03333333){45}{\line(1,0){.04444444}}
%\end
\put(57.5,49.5){\line(0,-1){1.75}}
%\emline(57.5,47.75)(55.25,46.25)
\multiput(57.5,47.75)(-.05,-.03333333){45}{\line(-1,0){.05}}
%\end
%\emline(56.25,51.25)(58.75,49.5)
\multiput(56.25,51.25)(.04807692,-.03365385){52}{\line(1,0){.04807692}}
%\end
%\emline(58.75,49.5)(58.25,47.75)
\multiput(58.75,49.5)(-.0333333,-.1166667){15}{\line(0,-1){.1166667}}
%\end
%\emline(58.25,47.75)(56.75,46.5)
\multiput(58.25,47.75)(-.03947368,-.03289474){38}{\line(-1,0){.03947368}}
%\end
%\emline(57,46)(55.5,45)
\multiput(57,46)(-.05,-.0333333){30}{\line(-1,0){.05}}
%\end
\put(54.5,46.5){\line(0,-1){2.25}}
\put(54.5,44.25){\line(0,1){0}}
%\emline(56.5,46.5)(55.5,44.5)
\multiput(56.5,46.5)(-.0333333,-.0666667){30}{\line(0,-1){.0666667}}
%\end
%\emline(55.5,44.5)(54.25,44.75)
\multiput(55.5,44.5)(-.15625,.03125){8}{\line(-1,0){.15625}}
%\end
%\emline(54.75,45)(50.25,41)
\multiput(54.75,45)(-.037815126,-.033613445){119}{\line(-1,0){.037815126}}
%\end
%\emline(50.25,41)(50.5,41.25)
\multiput(50.25,41)(.03125,.03125){8}{\line(0,1){.03125}}
%\end
%\emline(55.75,44.5)(51.75,41)
\multiput(55.75,44.5)(-.038461538,-.033653846){104}{\line(-1,0){.038461538}}
%\end
\put(47.75,41.25){\line(1,0){4.25}}
%\emline(52,41.25)(49.5,39.5)
\multiput(52,41.25)(-.04807692,-.03365385){52}{\line(-1,0){.04807692}}
%\end
%\emline(48,41.25)(49.75,39.75)
\multiput(48,41.25)(.03888889,-.03333333){45}{\line(1,0){.03888889}}
%\end
%\emline(49.25,66)(51.5,68.25)
\multiput(49.25,66)(.03358209,.03358209){67}{\line(0,1){.03358209}}
%\end
%\emline(51.5,68.25)(53.25,66.25)
\multiput(51.5,68.25)(.03365385,-.03846154){52}{\line(0,-1){.03846154}}
%\end
%\emline(51.5,66)(43.75,55.5)
\multiput(51.5,66)(-.033695652,-.045652174){230}{\line(0,-1){.045652174}}
%\end
%\emline(49.75,66.25)(42,55.25)
\multiput(49.75,66.25)(-.033695652,-.047826087){230}{\line(0,-1){.047826087}}
%\end
\put(42.25,55.25){\line(2,-5){5.5}}
\put(43.75,55.75){\line(2,-5){5.5}}
\put(52.25,54.5){$\cal C$}
\put(62.25,60.75){$\cal D$}
\put(63.25,54.25){thick lens}
\end{picture}

\begin{lemma}\label{thick} Every edge of the outer boundary component $\bf x$ of a thick lens $\Gamma$ is either unlabeled or labeled by a letter from the alphabet $A.$
\end{lemma}

\proof Every edge of an $\alpha$-cell (of an $\omega$-cell) of $\Gamma$ lying on $\bf x$ is unlabeled since
$\alpha$ (resp., $\omega$) can occur only as the left-most (resp., the right-most) letter
in the relator words of $H.$ It follows from the definitions of bands and $\Gamma$
that the edges of $\alpha p$-cells and of $q_1\omega$-cells belonging to $\bf x$ are unlabeled too.

Since the cells of $\cal C$ corresponding to the relation $p=q_1\omega$ and to $q_1\omega=p$ alternate
in its maximal $q$-band $\cal C$ of $\Gamma,$ every maximal $\omega$-band of $\Gamma$ must start with a cell of $\cal C$ corresponding to $p=q_1\omega$ and end on the next $q_1\omega$-cell of $\cal C$ corresponding to $q_1\omega=p.$ 
Hence the only cells of $\cal C$ having edges in $\bf x$ are $p$-cells, and labeled edges of $\bf x$ are their $a$-edges on
the right side of $\cal C.$ So they are labeled by letters from $A$ (see Relations (\ref{aux})).
\endproof

We say that a closed region $O$ in a derivation trapezium is {\it generated} by a thick lens $\Gamma$
(or by the $\alpha q$-lens $E$ defining $\Gamma$)
if (1) $O$ contains $\Gamma$; (2) if $O$ contains  an edge $e$ of a cell $\pi$ and $e$ is labeled by a letter from $A,$  then $O$ contains $\pi;$ (3) if $O$ contains  an edge from the outer boundary of some thick lens $\Gamma'$,
then $O$ contains $\Gamma'$ (4) $O$ is minimal with respect to (1)--(3). 

\begin{lemma}\label{O} Let a region $O$ of a derivation trapezium is {\it generated} by a $\alpha q$-lens $E$ or by a thick lens $\Gamma.$
Then every edge in the outer boundary component of $O$ is either unlabeled or has a label from $A.$ 
\end{lemma} 

\proof By the definition of $O,$ it constructed from several thick lenses and several maximal $a$-bands
which start/terminate on the thick lenses and correspond to $a$-letters from $A.$ So Lemma \ref{thick}
completes the proof. 
\endproof

\begin{lemma}\label{lens} Let $\Delta$ be a minimal trapezium over $H.$ Then 

(a) An $\alpha q$-lens $\Gamma$ of $\Delta$
encloses no other $\alpha q$-lenses and no $\omega$-cells.

(b) Let $\Delta$ have a through $q$-band $\cal C,$
and assume that the top and bottom edges of $\Delta$ from the left of $\cal C$ are labeled by
letters from $\{\{\alpha\}\cup Y_l\}.$
Then there are no $\alpha q$-lenses and no $\omega$-cells 
from the left of $\cal C.$

(c) Assume that an $\alpha q$-cap (or cup) $\Gamma$ of $\Delta$ encloses 
a $\alpha q$-lens $E.$ Then the closed region $O$ generated by $E$ does not share any labeled edge with $\Gamma.$ 
 \end{lemma}

\proof {\bf (a)} Assume that an $\alpha q$-lens $E$ is enclosed in $\Gamma$, and
there are no bigger $\alpha q$-lenses enclosed in $\Gamma$ and surrounding $E.$ 
Note that the region $O$ generated by the $\alpha q$-lens $E$ is also enclosed in $\Gamma.$
By Lemma \ref{O}, every labeled edge $e$ of the outer
boundary component $\bf p$ of $O$ must be connected by a maximal $a$-band
$\cal A$ (of length $\ge 0$) with an $a$-edge $f$ on the left side 
of  the maximal $q$-band
$\cal C$ of $\Gamma.$ However $f$ is labeled by a letter from $Y_l$ while $e$ is
labeled by a letter from $A.$  This contradicts to the condition
$Y_l\cap A=\emptyset,$ and so $\bf p$ has no labeled edges.

Then Lemma \ref{unlab} gives another contradiction since $O$ contains $q$-edges
of the $\alpha q$-lens $E.$  Hence our assumption false, and $\Gamma$ surrounds no $\alpha q$-lenses. 

The second assertion of (a)  is also true since an $\omega$-band enclosed in $\Gamma$ 
cannot start/end on $\cal C.$ 

{\bf (b)}  The same proof as for (a), but now $\Gamma$ is the part of $\Delta$
from the left of $\cal C.$

{\bf (c)} Follows from Lemma \ref{O}. 
\endproof

\begin{lemma}\label{qa}  (a) Assume that an $a$-band $\cal A$ starts and ends on
a $q$-band $\cal C$ of a minimal trapezium, 
and $\cal C$ has no edges labeled by $p.$
Then $\cal A$ and $\cal C$
surrounds no  $\alpha q$-lenses and no $\omega$-cells. 

(b) Let an $\omega$-band $\cal D$ start and end on
a $q$-band $\cal C$ of a minimal trapezium $\Delta.$ Then $\cal C$ and $\cal D$
surround no $\omega$-, $\alpha$-, or $q$-cells.
\end{lemma}
\proof {\bf (a)}
It follows from the assumption of the lemma that the $a$-band $\cal A$ corresponds to
a letter from $Y_l$ (from $Y_r$) if it is disposed from the left (resp., from the right) of $\cal C.$
%The `left' version can be considered as in Lemma \ref{lens}: 
Now the assumption that
$\cal A$ and $\cal C$ surround an $\alpha q$-lens $E$ gives a contradiction as in Lemma \ref{lens}. 
%To obtain a similar
%contradiction for the `right' version we should explain that no $a$-band ${\cal A}_1$ starting
%(or terminating) on $\cal C$ from the right can terminate (or start) on the maximal
%$q$-band ${\cal C}_1$ of $E$ (or on the $q$-band of other surrounded $\alpha q$-lens).
%Indeed, since the $a$-band ${\cal A}_1$ corresponds to a letter $a\in Y_r,$ it should
%terminate on a cell $\pi$ of ${\cal C}_1$ having no edges labeled by $p.$ But by the definition
%of $\alpha q$-lens, ${\cal C}_1$ itself starts and terminates with cells having $p$ in the boundary
%labels, and only the $q_1\omega$-relation change $p$ by a different $q$-letter (or vice versa). It %follows from Lemma \ref{omega}  that $\pi$ must be disposed on ${\cal C}_1$ between
%two $q_1\omega$-cells connected from the right by an $\omega$-band. We come to a contradiction
%since ${\cal A}_1$ cannot cross this $\omega$-band. \footnote{Risunok?}
The second claim is  obvious since $\cal C$ has
no $p$-edges, and so no $\omega$-band can start on the band $\cal C.$

{\bf (b)}
Let $\cal D$ start with a $q_1\omega$-cell $\pi$ of $\cal C.$ Then $\pi$ corresponds
to the relation $p=q_1\omega.$ Therefore the next $q_1\omega$-cell $\pi'$ of $\cal C$
must correspond to the relation $q_1\omega= p,$ and some maximal $\omega$-band $\cal D'$
terminates at $\pi'.$ Since by Lemma \ref{omega}, $\cal D'$ must also start on 
$\cal C,$ and different maximal $\omega$-bands cannot cross each other, we conclude
that $\cal D'=\cal D.$ In other words, if  
the closed region $\Gamma$
bounded by $\cal C$ and $\cal D$ (where the cells from $\cal C$ and $\cal D$ do
not belong to $\Gamma$) encloses a maximal $\omega$-band, then such a band must connect
two cells of an $\alpha q$-lens enclosed in $\Gamma.$
If $\Gamma$ contains
$\alpha$- or $q$-cells, then $\Gamma$ contains $\alpha q$-lenses as well. But
such an assumption leads to a contradiction exactly as in the `right' version of Part (a). 
\endproof 

\subsection{Types of $q$-bands in minimal trapezia}

A derivation trapezium $\Delta$ will be called a {\it machine trapezium} if the top
or the bottom label $w$ of $\Delta$ is a configuration of the machine $M_5$ and every
nontrivial cell corresponds to one of the machine relation (\ref{machine}) (i.e.,
there are no cells corresponding to the auxiliary relations (\ref{aux})).

Let $w$ and $w'$ be the bottom and the top labels of a derivation band of a machine
trapezium $\Delta.$ Then it follows by the induction on the height of $\Delta$ that either $w'=w$ or
this band correspond to a transition $w\to w'$ of $M_5.$  Therefore both the bottom
and the top labels of $\Delta$ are configurations of $M_5,$ and they can be connected by a computation of $M_5.$

The definition of {\it peeled machine trapezium} is similar but now $w$ can be a configuration
of $M_5$ without one of the endmarkers $\alpha$ or $\omega,$ or without both. Hence the bottom
and top labels of a peeled machine trapezium plus the additional letter $\alpha$ in the beginning or/and the letter $\omega$ at the end of them are connected by a  computation 
of $M_5$ without commands  involving $\alpha$ or $\omega,$ or both, resp.

\begin{lemma} \label{type2} Assume that $\Delta$ be a minimal trapezium over $H.$

(a) Let $\Gamma$ be an  $\alpha q$-lens in $\Delta$ formed by an $\alpha$-band $\cal B$
and $q$-band $\cal C.$  Then the type of $\cal C$ is $2.$

(b) Assume that $\Delta$  has a through $q$-band $\cal C$ and  a through $\alpha$-band
$\cal B$ from the left of $\cal C.$ Then there is no horizontal path $\bf x$ starting with
an $\alpha$-edge of $\cal B,$ ending with an $p$-edge of $\cal C,$ and having label
$\alpha Up,$ where $U$ is a word in the alphabet $A_l.$ 

(c) Assume that $\Delta$ has an $\alpha q$-cap or an $\alpha q$-cup $\Gamma$ with
maximal $\alpha$-band $\cal B$ and maximal $q$-band $\cal C.$ Also assume that there
are neither $\alpha q$-lenses nor $\alpha q$-caps/cups enclosed in $\Gamma.$ Then the type of $\cal C$ is at most $1.$
\end{lemma}

\proof {\bf (a)} Denote by $\pi_1$ and $\pi_m$ the first and the last cells of $\cal C.$
They correspond to the relations $1=\alpha p$ and $\alpha p =1,$ resp. Therefore
$\cal C$ has equal number of cells corresponding to the relation $p=q_1\omega$ and
to $q_1\omega =p,$ in particular, the type $t$ of $\cal C$ is even.

{\bf Case 1.} Assume that $t=0.$ Then each of the cells $\pi_2,\dots\pi_{m-1}$ is
either trivial or corresponds to a relation $pa = a_lp$ (or to $a_lp=pa$).
Note that by Lemma \ref{lens}, there are neither $\alpha q$-lenses nor $\omega$-cells enclosed in $\Gamma.$
Hence every cell between $\cal B$ and $\cal C$ is a trivial $a$-cell and it belongs to an
$a$-band starting and ending on $\cal C$ and corresponding to a letter from $A_l.$
A right-most $a$-band $\cal A$ enclosed in $\Gamma$ and a part of $\cal C$ form 
a derivation (sub)trapezium with top and bottom labels  equal to $pa$ ($a\in A$).
Therefore one can replace this subtrapezium by a trapezium having only trivial cells.
%TeXCAD (http://texcad.sf.net/) Picture. File: [semi7.pic]. Options on following lines.
%\grade{\on}
%\emlines{\off}
%\epic{\off}
%\beziermacro{\on}
%\reduce{\on}
%\snapping{\off}
%\pvinsert{% Your \input, \def, etc. here}
%\quality{8.000}
%\graddiff{0.005}
%\snapasp{1}
%\zoom{4.0000}
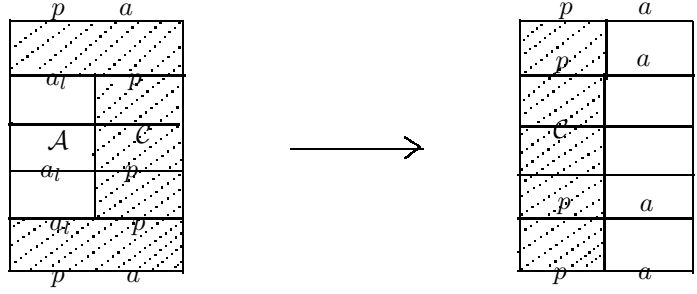
\unitlength 1mm % = 2.845pt
\linethickness{0.4pt}
\ifx\plotpoint\undefined\newsavebox{\plotpoint}\fi % GNUPLOT compatibility
\begin{picture}(133.75,40.5)(20,50)
\put(42.75,90.75){\line(0,-1){33.5}}
\put(110.5,90.75){\line(0,-1){33.5}}
\put(43,90.5){\line(1,0){22.75}}
\put(110.75,90.5){\line(1,0){22.75}}
\put(65.75,90.5){\line(0,-1){33.5}}
\put(133.5,90.5){\line(0,-1){33.5}}
\put(42.5,57.25){\line(1,0){23.25}}
\put(110.25,57.25){\line(1,0){23.25}}
\put(42.75,83.25){\line(1,0){23.25}}
\put(110.5,83.25){\line(1,0){23.25}}
\put(42.5,64.25){\line(1,0){23.25}}
\put(110.25,64.25){\line(1,0){23.25}}
\put(54,83.5){\line(0,-1){19}}
\put(121.75,83.5){\line(0,-1){19}}
\put(121.75,65){\line(0,-1){7.5}}
\put(122,90.75){\line(0,-1){8}}
\put(47.75,73.5){$\cal A$}
\put(59.5,74.25){$\cal C$}
\put(47.5,82){$a_l$}
\put(58.5,82){$p$}
\put(48,62.75){$a_l$}
\put(59,62.75){$p$}
\put(48.25,55.75){$p$}
\put(58.25,55.75){$a$}
\put(48.25,91.25){$p$}
\put(57.25,91.25){$a$}
\put(115.75,91.5){$p$}
\put(115.25,84.25){$p$}
\put(115.5,65.5){$p$}
\put(115,56.25){$p$}
\put(126.25,91.5){$a$}
\put(126,84.5){$a$}
\put(126.5,65.25){$a$}
\put(126.25,56.25){$a$}
\put(115,75){$\cal C$}
%\dottedline(42.75,86.75)(47.5,90.25)
\multiput(42.68,86.68)(.67857,.5){8}{{\rule{.4pt}{.4pt}}}
%\end
%\dottedline(43,84.25)(51.25,90.5)
\multiput(42.93,84.18)(.6875,.52083){13}{{\rule{.4pt}{.4pt}}}
%\end
%\dottedline(46.75,84.25)(54.25,90.5)
\multiput(46.68,84.18)(.68182,.56818){12}{{\rule{.4pt}{.4pt}}}
%\end
%\dottedline(49.5,83.75)(56.75,90.25)
\multiput(49.43,83.68)(.65909,.59091){12}{{\rule{.4pt}{.4pt}}}
%\end
%\dottedline(52,83.5)(59.75,90.25)
\multiput(51.93,83.43)(.70455,.61364){12}{{\rule{.4pt}{.4pt}}}
%\end
%\dottedline(55.25,83.5)(62.75,90.25)
\multiput(55.18,83.43)(.68182,.61364){12}{{\rule{.4pt}{.4pt}}}
%\end
%\dottedline(58.75,84)(65.25,90.25)
\multiput(58.68,83.93)(.65,.625){11}{{\rule{.4pt}{.4pt}}}
%\end
%\dottedline(54.5,80.25)(56.25,82)
\multiput(54.43,80.18)(.5833,.5833){4}{{\rule{.4pt}{.4pt}}}
%\end
%\dottedline(54.5,77.75)(65.25,87)
\multiput(54.43,77.68)(.71667,.61667){16}{{\rule{.4pt}{.4pt}}}
%\end
%\dottedline(54,75.5)(65.25,84.5)
\multiput(53.93,75.43)(.70313,.5625){17}{{\rule{.4pt}{.4pt}}}
%\end
%\dottedline(53.75,72.75)(65.75,82.25)
\multiput(53.68,72.68)(.70588,.55882){18}{{\rule{.4pt}{.4pt}}}
%\end
%\dottedline(54.25,70.5)(65,79)
\multiput(54.18,70.43)(.71667,.56667){16}{{\rule{.4pt}{.4pt}}}
%\end
%\dottedline(54.25,68.25)(65,76.5)
\multiput(54.18,68.18)(.71667,.55){16}{{\rule{.4pt}{.4pt}}}
%\end
%\dottedline(54,66)(65.25,74.25)
\multiput(53.93,65.93)(.75,.55){16}{{\rule{.4pt}{.4pt}}}
%\end
%\dottedline(54.5,64.5)(65.5,72)
\multiput(54.43,64.43)(.73333,.5){16}{{\rule{.4pt}{.4pt}}}
%\end
%\dottedline(53.75,64)(45.5,57.5)
\multiput(53.68,63.93)(-.6875,-.54167){13}{{\rule{.4pt}{.4pt}}}
%\end
%\dottedline(42.75,58)(51.25,64)
\multiput(42.68,57.93)(.70833,.5){13}{{\rule{.4pt}{.4pt}}}
%\end
%\dottedline(43,60.75)(48.25,64.25)
\multiput(42.93,60.68)(.75,.5){8}{{\rule{.4pt}{.4pt}}}
%\end
%\dottedline(43.25,63.25)(44.5,64)
\multiput(43.18,63.18)(.625,.375){3}{{\rule{.4pt}{.4pt}}}
%\end
%\dottedline(48.5,58)(65.25,70)
\multiput(48.43,57.93)(.76136,.54545){23}{{\rule{.4pt}{.4pt}}}
%\end
%\dottedline(50.75,57.5)(65.25,68)
\multiput(50.68,57.43)(.76316,.55263){20}{{\rule{.4pt}{.4pt}}}
%\end
%\dottedline(54.25,57.5)(65.25,65.75)
\multiput(54.18,57.43)(.73333,.55){16}{{\rule{.4pt}{.4pt}}}
%\end
%\dottedline(57.25,57.75)(65.5,63.75)
\multiput(57.18,57.68)(.6875,.5){13}{{\rule{.4pt}{.4pt}}}
%\end
%\dottedline(60,57.75)(65.5,62.25)
\multiput(59.93,57.68)(.6875,.5625){9}{{\rule{.4pt}{.4pt}}}
%\end
%\dottedline(61.75,57.5)(65.5,60.5)
\multiput(61.68,57.43)(.625,.5){7}{{\rule{.4pt}{.4pt}}}
%\end
%\dottedline(64.75,57.75)(65.75,58.75)
\multiput(64.68,57.68)(.3333,.3333){4}{{\rule{.4pt}{.4pt}}}
%\end
%\dottedline(111,88)(114,90.5)
\multiput(110.93,87.93)(.6,.5){6}{{\rule{.4pt}{.4pt}}}
%\end
%\dottedline(111.25,86.25)(116.5,90.25)
\multiput(111.18,86.18)(.65625,.5){9}{{\rule{.4pt}{.4pt}}}
%\end
%\dottedline(111.25,83.75)(119.75,90.25)
\multiput(111.18,83.68)(.70833,.54167){13}{{\rule{.4pt}{.4pt}}}
%\end
%\dottedline(111,80)(122,89.25)
\multiput(110.93,79.93)(.6875,.57813){17}{{\rule{.4pt}{.4pt}}}
%\end
%\dottedline(110.5,77.75)(121.75,86.25)
\multiput(110.43,77.68)(.75,.56667){16}{{\rule{.4pt}{.4pt}}}
%\end
%\dottedline(110.75,75.5)(121.25,82.75)
\multiput(110.68,75.43)(.75,.51786){15}{{\rule{.4pt}{.4pt}}}
%\end
%\dottedline(110.75,73)(121.75,81)
\multiput(110.68,72.93)(.73333,.53333){16}{{\rule{.4pt}{.4pt}}}
%\end
%\dottedline(110.75,70.75)(120.75,78.25)
\multiput(110.68,70.68)(.71429,.53571){15}{{\rule{.4pt}{.4pt}}}
%\end
%\dottedline(110.5,68.75)(121,76)
\multiput(110.43,68.68)(.75,.51786){15}{{\rule{.4pt}{.4pt}}}
%\end
%\dottedline(110.75,67.25)(121.25,73.75)
\multiput(110.68,67.18)(.80769,.5){14}{{\rule{.4pt}{.4pt}}}
%\end
%\dottedline(110.75,64.75)(121.75,72)
\multiput(110.68,64.68)(.73333,.48333){16}{{\rule{.4pt}{.4pt}}}
%\end
%\dottedline(110.5,62.25)(121.25,69.25)
\multiput(110.43,62.18)(.76786,.5){15}{{\rule{.4pt}{.4pt}}}
%\end
%\dottedline(110.25,59.75)(121.75,67.25)
\multiput(110.18,59.68)(.76667,.5){16}{{\rule{.4pt}{.4pt}}}
%\end
%\dottedline(110.75,57.5)(121.5,65)
\multiput(110.68,57.43)(.76786,.53571){15}{{\rule{.4pt}{.4pt}}}
%\end
%\dottedline(114.5,58)(121.25,62.75)
\multiput(114.43,57.93)(.75,.52778){10}{{\rule{.4pt}{.4pt}}}
%\end
%\dottedline(117.5,58)(121.5,61)
\multiput(117.43,57.93)(.57143,.42857){8}{{\rule{.4pt}{.4pt}}}
%\end
%\dottedline(119.5,58)(121.75,59.75)
\multiput(119.43,57.93)(.5625,.4375){5}{{\rule{.4pt}{.4pt}}}
%\end
%\dottedline(43,77)(42.75,77)
\multiput(42.93,76.93)(-.125,0){3}{{\rule{.4pt}{.4pt}}}
%\end
\put(42.5,76.75){\line(1,0){22.75}}
\put(42.75,70.5){\line(1,0){22.75}}
\put(110.5,76.5){\line(1,0){22.75}}
\put(110.25,70){\line(1,0){23.5}}
\put(46.75,69.75){$a_l$}
\put(58,69.75){$p$}
\put(80,73.5){\line(1,0){17.75}}
%\emline(95.25,75)(97.25,73.75)
\multiput(95.25,75)(.05263158,-.03289474){38}{\line(1,0){.05263158}}
%\end
%\emline(97.25,73.75)(95.25,72.25)
\multiput(97.25,73.75)(-.04444444,-.03333333){45}{\line(-1,0){.04444444}}
%\end
\end{picture}

This surgery reduces the number of nontrivial cells in the $q$-band $\cal C$. 
Finally, we will have an $\alpha q$-lens with unlabeled (outer) boundary. Hence the $\alpha q$-lens can
be removed from $\Delta$ by Lemma \ref{unlab}. Since $\Delta$ is a minimal trapezium,
the case $t=0$ is not possible.

{\bf Case 2.} Assume now that $t\ge 4,$ and let $\pi_i$ and $\pi_j$ be the first cells of 
$\cal C$ corresponding to the relations  $p=q_1\omega$ and $q_1\omega =p,$ respectively.
Let $\Delta_0$ be the trapezium formed by the derivation bands ${\cal T}_k,\dots,{\cal T}_{k+j-i}$ of $\Delta$ containing $\pi_i,\dots, \pi_j,$ resp. There are two vertical paths dividing
$\Delta_0$: The left side of $\cal B$ and 
%the path formed by the right sides of $\pi_i$, $\pi_j$, and by 
the right side of the $\omega$-band starting with $\pi_i$ and ending with $\pi_j$ (see Lemma \ref{omega}).
These paths divide $\Delta_0$ into $3$ subtrapezia $\Delta_1,
\Delta_2, \Delta_3$ (from left to right).

%TeXCAD (http://texcad.sf.net/) Picture. File: [semi6.pic]. Options on following lines.
%\grade{\on}
%\emlines{\off}
%\epic{\off}
%\beziermacro{\on}
%\reduce{\on}
%\snapping{\off}
%\pvinsert{% Your \input, \def, etc. here}
%\quality{8.000}
%\graddiff{0.005}
%\snapasp{1}
%\zoom{4.0000}
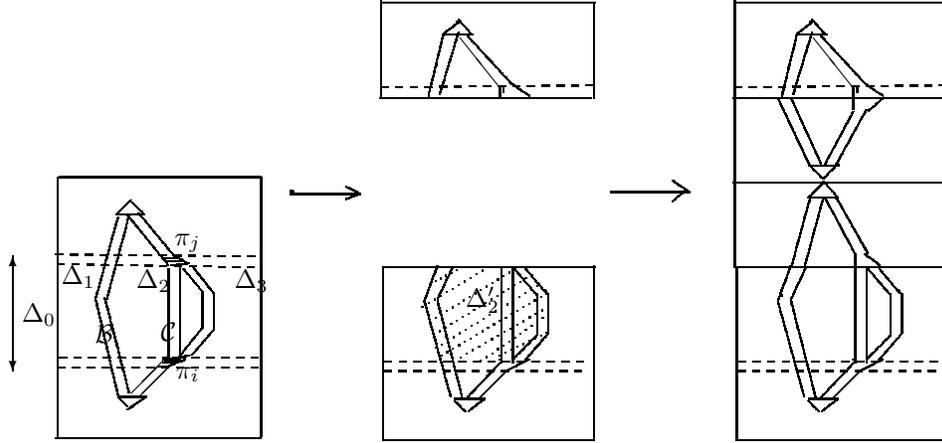
\unitlength 1mm % = 2.845pt
\linethickness{0.4pt}
\ifx\plotpoint\undefined\newsavebox{\plotpoint}\fi % GNUPLOT compatibility
\begin{picture}(148.75,65.75)(20,30)
\put(30.25,72.25){\line(0,-1){34.75}}
%\emline(40,69.25)(38,67.25)
\multiput(40,69.25)(-.03333333,-.03333333){60}{\line(0,-1){.03333333}}
%\end
\put(38,67.25){\line(1,0){3.75}}
%\emline(41.75,67.25)(40,69)
\multiput(41.75,67.25)(-.03365385,.03365385){52}{\line(0,1){.03365385}}
%\end
%\emline(38.25,67)(35.25,56)
\multiput(38.25,67)(-.03370787,-.12359551){89}{\line(0,-1){.12359551}}
%\end
%\emline(39.75,67.25)(36.75,56)
\multiput(39.75,67.25)(-.03370787,-.12640449){89}{\line(0,-1){.12640449}}
%\end
\put(38.25,43){\line(1,0){3.75}}
%\emline(42,43)(39.5,41.25)
\multiput(42,43)(-.04807692,-.03365385){52}{\line(-1,0){.04807692}}
%\end
%\emline(38.25,43)(39.5,41.75)
\multiput(38.25,43)(.03289474,-.03289474){38}{\line(0,-1){.03289474}}
%\end
\put(40,43.5){\line(1,1){5}}
\put(41.25,43.25){\line(1,1){3.75}}
%\emline(44.75,47)(47.25,48)
\multiput(44.75,47)(.0833333,.0333333){30}{\line(1,0){.0833333}}
%\end
\put(47.25,48){\line(0,1){0}}
\put(43.75,47){\line(1,0){1}}
%\emline(39.5,67.25)(45,60.5)
\multiput(39.5,67.25)(.033536585,-.041158537){164}{\line(0,-1){.041158537}}
%\end
\put(45,60.5){\line(0,1){0}}
%\emline(41,67.5)(45.75,61.75)
\multiput(41,67.5)(.033687943,-.040780142){141}{\line(0,-1){.040780142}}
%\end
%\emline(45.75,61.75)(47.5,60.75)
\multiput(45.75,61.75)(.0583333,-.0333333){30}{\line(1,0){.0583333}}
%\end
\put(30.5,72.25){\line(1,0){26.75}}
\put(57.25,72.25){\line(0,-1){34.5}}
\put(30,37.5){\line(1,0){27.25}}
%\dashline{1}(30.25,60.75)(57.5,60.25)
\put(30.18,60.68){\line(1,0){.9732}}
\put(32.126,60.644){\line(1,0){.9732}}
\put(34.073,60.608){\line(1,0){.9732}}
\put(36.019,60.573){\line(1,0){.9732}}
\put(37.965,60.537){\line(1,0){.9732}}
\put(39.912,60.501){\line(1,0){.9732}}
\put(41.858,60.465){\line(1,0){.9732}}
\put(43.805,60.43){\line(1,0){.9732}}
\put(45.751,60.394){\line(1,0){.9732}}
\put(47.698,60.358){\line(1,0){.9732}}
\put(49.644,60.323){\line(1,0){.9732}}
\put(51.59,60.287){\line(1,0){.9732}}
\put(53.537,60.251){\line(1,0){.9732}}
\put(55.483,60.215){\line(1,0){.9732}}
%\end
%\dashline{1}(30,62)(57.5,61.75)
\put(29.93,61.93){\line(1,0){.9821}}
\put(31.894,61.912){\line(1,0){.9821}}
\put(33.858,61.894){\line(1,0){.9821}}
\put(35.823,61.876){\line(1,0){.9821}}
\put(37.787,61.858){\line(1,0){.9821}}
\put(39.751,61.84){\line(1,0){.9821}}
\put(41.715,61.823){\line(1,0){.9821}}
\put(43.68,61.805){\line(1,0){.9821}}
\put(45.644,61.787){\line(1,0){.9821}}
\put(47.608,61.769){\line(1,0){.9821}}
\put(49.573,61.751){\line(1,0){.9821}}
\put(51.537,61.733){\line(1,0){.9821}}
\put(53.501,61.715){\line(1,0){.9821}}
\put(55.465,61.698){\line(1,0){.9821}}
%\end
%\dashline{1}(30.25,47)(57,47)
\put(30.18,46.93){\line(1,0){.9907}}
\put(32.161,46.93){\line(1,0){.9907}}
\put(34.143,46.93){\line(1,0){.9907}}
\put(36.124,46.93){\line(1,0){.9907}}
\put(38.106,46.93){\line(1,0){.9907}}
\put(40.087,46.93){\line(1,0){.9907}}
\put(42.069,46.93){\line(1,0){.9907}}
\put(44.05,46.93){\line(1,0){.9907}}
\put(46.032,46.93){\line(1,0){.9907}}
\put(48.013,46.93){\line(1,0){.9907}}
\put(49.995,46.93){\line(1,0){.9907}}
\put(51.976,46.93){\line(1,0){.9907}}
\put(53.957,46.93){\line(1,0){.9907}}
\put(55.939,46.93){\line(1,0){.9907}}
%\end
%\dashline{1}(30.5,48.25)(57,48.25)
\put(30.43,48.18){\line(1,0){.9815}}
\put(32.393,48.18){\line(1,0){.9815}}
\put(34.356,48.18){\line(1,0){.9815}}
\put(36.319,48.18){\line(1,0){.9815}}
\put(38.282,48.18){\line(1,0){.9815}}
\put(40.245,48.18){\line(1,0){.9815}}
\put(42.207,48.18){\line(1,0){.9815}}
\put(44.17,48.18){\line(1,0){.9815}}
\put(46.133,48.18){\line(1,0){.9815}}
\put(48.096,48.18){\line(1,0){.9815}}
\put(50.059,48.18){\line(1,0){.9815}}
\put(52.022,48.18){\line(1,0){.9815}}
\put(53.985,48.18){\line(1,0){.9815}}
\put(55.948,48.18){\line(1,0){.9815}}
%\end
\put(46.5,60.25){\line(0,-1){11.75}}
\put(45,60){\line(0,-1){11.75}}
%\emline(48.25,60.25)(51,56.75)
\multiput(48.25,60.25)(.03353659,-.04268293){82}{\line(0,-1){.04268293}}
%\end
\put(51,56.75){\line(0,-1){4.75}}
%\emline(51,52)(48,48.5)
\multiput(51,52)(-.03370787,-.03932584){89}{\line(0,-1){.03932584}}
%\end
%\emline(47,60)(49.5,57)
\multiput(47,60)(.03333333,-.04){75}{\line(0,-1){.04}}
%\end
\put(49.5,57){\line(0,-1){4.75}}
%\emline(49.5,52.25)(46.5,48.5)
\multiput(49.5,52.25)(-.03370787,-.04213483){89}{\line(0,-1){.04213483}}
%\end
\put(73.5,60.25){\line(1,0){28.25}}
\put(120.5,60.25){\line(1,0){28.25}}
\put(73.5,60.25){\line(0,-1){23}}
\put(120.5,60.25){\line(0,-1){23}}
\put(73.5,37.25){\line(1,0){28}}
\put(120.5,37.25){\line(1,0){28}}
\put(101.5,37.25){\line(0,1){23.25}}
\put(148.5,37.25){\line(0,1){23.25}}
\put(89.25,60.25){\line(0,-1){12.75}}
\put(136.25,60.25){\line(0,-1){12.75}}
\put(90.75,60.25){\line(0,-1){12.75}}
\put(137.75,60.25){\line(0,-1){12.75}}
%\dashline{1}(73.75,47.75)(101.25,47.75)
\put(73.68,47.68){\line(1,0){.9821}}
\put(75.644,47.68){\line(1,0){.9821}}
\put(77.608,47.68){\line(1,0){.9821}}
\put(79.573,47.68){\line(1,0){.9821}}
\put(81.537,47.68){\line(1,0){.9821}}
\put(83.501,47.68){\line(1,0){.9821}}
\put(85.465,47.68){\line(1,0){.9821}}
\put(87.43,47.68){\line(1,0){.9821}}
\put(89.394,47.68){\line(1,0){.9821}}
\put(91.358,47.68){\line(1,0){.9821}}
\put(93.323,47.68){\line(1,0){.9821}}
\put(95.287,47.68){\line(1,0){.9821}}
\put(97.251,47.68){\line(1,0){.9821}}
\put(99.215,47.68){\line(1,0){.9821}}
%\end
%\dashline{1}(120.75,47.75)(148.25,47.75)
\put(120.68,47.68){\line(1,0){.9821}}
\put(122.644,47.68){\line(1,0){.9821}}
\put(124.608,47.68){\line(1,0){.9821}}
\put(126.573,47.68){\line(1,0){.9821}}
\put(128.537,47.68){\line(1,0){.9821}}
\put(130.501,47.68){\line(1,0){.9821}}
\put(132.465,47.68){\line(1,0){.9821}}
\put(134.43,47.68){\line(1,0){.9821}}
\put(136.394,47.68){\line(1,0){.9821}}
\put(138.358,47.68){\line(1,0){.9821}}
\put(140.323,47.68){\line(1,0){.9821}}
\put(142.287,47.68){\line(1,0){.9821}}
\put(144.251,47.68){\line(1,0){.9821}}
\put(146.215,47.68){\line(1,0){.9821}}
%\end
%\dashline{1}(73.75,46.5)(101,46.5)
\put(73.68,46.43){\line(1,0){.9732}}
\put(75.626,46.43){\line(1,0){.9732}}
\put(77.573,46.43){\line(1,0){.9732}}
\put(79.519,46.43){\line(1,0){.9732}}
\put(81.465,46.43){\line(1,0){.9732}}
\put(83.412,46.43){\line(1,0){.9732}}
\put(85.358,46.43){\line(1,0){.9732}}
\put(87.305,46.43){\line(1,0){.9732}}
\put(89.251,46.43){\line(1,0){.9732}}
\put(91.198,46.43){\line(1,0){.9732}}
\put(93.144,46.43){\line(1,0){.9732}}
\put(95.09,46.43){\line(1,0){.9732}}
\put(97.037,46.43){\line(1,0){.9732}}
\put(98.983,46.43){\line(1,0){.9732}}
%\end
%\dashline{1}(120.75,46.5)(148,46.5)
\put(120.68,46.43){\line(1,0){.9732}}
\put(122.626,46.43){\line(1,0){.9732}}
\put(124.573,46.43){\line(1,0){.9732}}
\put(126.519,46.43){\line(1,0){.9732}}
\put(128.465,46.43){\line(1,0){.9732}}
\put(130.412,46.43){\line(1,0){.9732}}
\put(132.358,46.43){\line(1,0){.9732}}
\put(134.305,46.43){\line(1,0){.9732}}
\put(136.251,46.43){\line(1,0){.9732}}
\put(138.198,46.43){\line(1,0){.9732}}
\put(140.144,46.43){\line(1,0){.9732}}
\put(142.09,46.43){\line(1,0){.9732}}
\put(144.037,46.43){\line(1,0){.9732}}
\put(145.983,46.43){\line(1,0){.9732}}
%\end
%\dashline{1}(89,48)(89.25,48)
\put(88.93,47.93){\line(1,0){.125}}
%\end
%\dashline{1}(136,48)(136.25,48)
\put(135.93,47.93){\line(1,0){.125}}
%\end
%\emline(89.25,47.75)(84.5,42.75)
\multiput(89.25,47.75)(-.033687943,-.035460993){141}{\line(0,-1){.035460993}}
%\end
%\emline(136.25,47.75)(131.5,42.75)
\multiput(136.25,47.75)(-.033687943,-.035460993){141}{\line(0,-1){.035460993}}
%\end
%\emline(93,48)(90,46.5)
\multiput(93,48)(-.06666667,-.03333333){45}{\line(-1,0){.06666667}}
%\end
%\emline(140,48)(137,46.5)
\multiput(140,48)(-.06666667,-.03333333){45}{\line(-1,0){.06666667}}
%\end
\put(90,46.5){\line(-1,-1){3.75}}
\put(137,46.5){\line(-1,-1){3.75}}
\put(83,42.75){\line(1,0){4}}
\put(130,42.75){\line(1,0){4}}
%\emline(87,42.75)(84.75,41)
\multiput(87,42.75)(-.04326923,-.03365385){52}{\line(-1,0){.04326923}}
%\end
%\emline(134,42.75)(131.75,41)
\multiput(134,42.75)(-.04326923,-.03365385){52}{\line(-1,0){.04326923}}
%\end
%\emline(83,42.5)(84.5,41)
\multiput(83,42.5)(.03333333,-.03333333){45}{\line(0,-1){.03333333}}
%\end
%\emline(130,42.5)(131.5,41)
\multiput(130,42.5)(.03333333,-.03333333){45}{\line(0,-1){.03333333}}
%\end
%\emline(91,60)(94,56.5)
\multiput(91,60)(.03370787,-.03932584){89}{\line(0,-1){.03932584}}
%\end
%\emline(138,60)(141,56.5)
\multiput(138,60)(.03370787,-.03932584){89}{\line(0,-1){.03932584}}
%\end
%\emline(92.5,60)(95.25,56.75)
\multiput(92.5,60)(.03353659,-.03963415){82}{\line(0,-1){.03963415}}
%\end
%\emline(139.5,60)(142.25,56.75)
\multiput(139.5,60)(.03353659,-.03963415){82}{\line(0,-1){.03963415}}
%\end
%\emline(95.25,56.75)(95,57.25)
\multiput(95.25,56.75)(-.03125,.0625){8}{\line(0,1){.0625}}
%\end
%\emline(142.25,56.75)(142,57.25)
\multiput(142.25,56.75)(-.03125,.0625){8}{\line(0,1){.0625}}
%\end
\put(94,56.75){\line(0,-1){5.5}}
\put(141,56.75){\line(0,-1){5.5}}
%\emline(94,51.25)(90.75,48)
\multiput(94,51.25)(-.033505155,-.033505155){97}{\line(0,-1){.033505155}}
%\end
%\emline(141,51.25)(137.75,48)
\multiput(141,51.25)(-.033505155,-.033505155){97}{\line(0,-1){.033505155}}
%\end
\put(95.5,56.75){\line(0,-1){5.25}}
\put(142.5,56.75){\line(0,-1){5.25}}
\put(95.5,51.5){\line(0,1){0}}
\put(142.5,51.5){\line(0,1){0}}
%\emline(95.5,51.5)(92.25,47.75)
\multiput(95.5,51.5)(-.033505155,-.038659794){97}{\line(0,-1){.038659794}}
%\end
%\emline(142.5,51.5)(139.25,47.75)
\multiput(142.5,51.5)(-.033505155,-.038659794){97}{\line(0,-1){.038659794}}
%\end
%\emline(79,55.75)(82.5,43)
\multiput(79,55.75)(.033653846,-.122596154){104}{\line(0,-1){.122596154}}
%\end
%\emline(126,55.75)(129.5,43)
\multiput(126,55.75)(.033653846,-.122596154){104}{\line(0,-1){.122596154}}
%\end
%\emline(80.75,55.75)(84,43)
\multiput(80.75,55.75)(.033505155,-.131443299){97}{\line(0,-1){.131443299}}
%\end
%\emline(127.75,55.75)(131,43)
\multiput(127.75,55.75)(.033505155,-.131443299){97}{\line(0,-1){.131443299}}
%\end
%\emline(79,55.75)(80.75,60.5)
\multiput(79,55.75)(.03365385,.09134615){52}{\line(0,1){.09134615}}
%\end
%\emline(126,55.75)(127.75,60.5)
\multiput(126,55.75)(.03365385,.09134615){52}{\line(0,1){.09134615}}
%\end
%\emline(80.75,56)(82.75,60.25)
\multiput(80.75,56)(.03333333,.07083333){60}{\line(0,1){.07083333}}
%\end
%\emline(127.75,56)(129.75,60.25)
\multiput(127.75,56)(.03333333,.07083333){60}{\line(0,1){.07083333}}
%\end
\put(73.25,82.75){\line(0,1){13}}
\put(120.25,82.75){\line(0,1){13}}
\put(73.25,95.5){\line(1,0){28.5}}
\put(120.25,95.5){\line(1,0){28.5}}
\put(101.5,95.5){\line(0,-1){12.75}}
\put(148.5,95.5){\line(0,-1){12.75}}
\put(73,82.75){\line(1,0){28.25}}
\put(120,82.75){\line(1,0){28.25}}
%\dashline{1}(73.25,84.25)(101.5,84.5)
\put(73.18,84.18){\line(1,0){.9741}}
\put(75.128,84.197){\line(1,0){.9741}}
\put(77.076,84.214){\line(1,0){.9741}}
\put(79.025,84.231){\line(1,0){.9741}}
\put(80.973,84.249){\line(1,0){.9741}}
\put(82.921,84.266){\line(1,0){.9741}}
\put(84.869,84.283){\line(1,0){.9741}}
\put(86.818,84.3){\line(1,0){.9741}}
\put(88.766,84.318){\line(1,0){.9741}}
\put(90.714,84.335){\line(1,0){.9741}}
\put(92.662,84.352){\line(1,0){.9741}}
\put(94.611,84.369){\line(1,0){.9741}}
\put(96.559,84.387){\line(1,0){.9741}}
\put(98.507,84.404){\line(1,0){.9741}}
\put(100.456,84.421){\line(1,0){.9741}}
%\end
%\dashline{1}(120.25,84.25)(148.5,84.5)
\put(120.18,84.18){\line(1,0){.9741}}
\put(122.128,84.197){\line(1,0){.9741}}
\put(124.076,84.214){\line(1,0){.9741}}
\put(126.025,84.231){\line(1,0){.9741}}
\put(127.973,84.249){\line(1,0){.9741}}
\put(129.921,84.266){\line(1,0){.9741}}
\put(131.869,84.283){\line(1,0){.9741}}
\put(133.818,84.3){\line(1,0){.9741}}
\put(135.766,84.318){\line(1,0){.9741}}
\put(137.714,84.335){\line(1,0){.9741}}
\put(139.662,84.352){\line(1,0){.9741}}
\put(141.611,84.369){\line(1,0){.9741}}
\put(143.559,84.387){\line(1,0){.9741}}
\put(145.507,84.404){\line(1,0){.9741}}
\put(147.456,84.421){\line(1,0){.9741}}
%\end
%\dashline{1}(83.5,93)(83,92.75)
\multiput(83.43,92.93)(-.0625,-.03125){4}{\line(-1,0){.0625}}
%\end
%\dashline{1}(130.5,93)(130,92.75)
\multiput(130.43,92.93)(-.0625,-.03125){4}{\line(-1,0){.0625}}
%\end
%\emline(83.75,93)(81.5,91.25)
\multiput(83.75,93)(-.04326923,-.03365385){52}{\line(-1,0){.04326923}}
%\end
%\emline(130.75,93)(128.5,91.25)
\multiput(130.75,93)(-.04326923,-.03365385){52}{\line(-1,0){.04326923}}
%\end
%\emline(83.5,93.25)(85.25,91.25)
\multiput(83.5,93.25)(.03365385,-.03846154){52}{\line(0,-1){.03846154}}
%\end
%\emline(130.5,93.25)(132.25,91.25)
\multiput(130.5,93.25)(.03365385,-.03846154){52}{\line(0,-1){.03846154}}
%\end
\put(81.5,91.25){\line(1,0){4}}
\put(128.5,91.25){\line(1,0){4}}
%\emline(81.5,91.25)(79.5,83)
\multiput(81.5,91.25)(-.03333333,-.1375){60}{\line(0,-1){.1375}}
%\end
%\emline(128.5,91.25)(126.5,83)
\multiput(128.5,91.25)(-.03333333,-.1375){60}{\line(0,-1){.1375}}
%\end
%\emline(83.5,91.5)(81,83.25)
\multiput(83.5,91.5)(-.03333333,-.11){75}{\line(0,-1){.11}}
%\end
%\emline(130.5,91.5)(128,83.25)
\multiput(130.5,91.5)(-.03333333,-.11){75}{\line(0,-1){.11}}
%\end
\put(83.75,90.75){\line(4,-5){5}}
\put(130.75,90.75){\line(4,-5){5}}
\put(89.5,84.25){\line(0,-1){.75}}
\put(136.5,84.25){\line(0,-1){.75}}
\put(89,84.25){\line(0,-1){.5}}
\put(136,84.25){\line(0,-1){.5}}
\put(89,83.75){\line(0,-1){.75}}
\put(136,83.75){\line(0,-1){.75}}
%\emline(84.75,91)(90.75,84.5)
\multiput(84.75,91)(.033707865,-.036516854){178}{\line(0,-1){.036516854}}
%\end
%\emline(131.75,91)(137.75,84.5)
\multiput(131.75,91)(.033707865,-.036516854){178}{\line(0,-1){.036516854}}
%\end
%\emline(90.75,84.5)(93,83)
\multiput(90.75,84.5)(.05,-.03333333){45}{\line(1,0){.05}}
%\end
%\emline(137.75,84.5)(140,83)
\multiput(137.75,84.5)(.05,-.03333333){45}{\line(1,0){.05}}
%\end
\put(120.25,83.25){\line(0,-1){23}}
\put(148.5,83){\line(0,-1){23.5}}
\put(120,71.5){\line(1,0){28.25}}
%\emline(130.25,73.75)(132,72)
\multiput(130.25,73.75)(.03365385,-.03365385){52}{\line(0,-1){.03365385}}
%\end
%\emline(132,72)(134,74)
\multiput(132,72)(.03333333,.03333333){60}{\line(0,1){.03333333}}
%\end
\put(134,74){\line(0,-1){.25}}
\put(134,73.75){\line(0,1){0}}
\put(134,73.75){\line(-1,0){4}}
\put(136,83.25){\line(0,-1){2}}
%\emline(140,82.75)(138,81)
\multiput(140,82.75)(-.03846154,-.03365385){52}{\line(-1,0){.03846154}}
%\end
%\emline(137.75,81.25)(133.75,74)
\multiput(137.75,81.25)(-.033613445,-.06092437){119}{\line(0,-1){.06092437}}
%\end
%\emline(136,81)(132.25,74)
\multiput(136,81)(-.033482143,-.0625){112}{\line(0,-1){.0625}}
%\end
%\emline(126,82.75)(130.25,73.75)
\multiput(126,82.75)(.033730159,-.071428571){126}{\line(0,-1){.071428571}}
%\end
%\emline(127.75,82.75)(131.75,74)
\multiput(127.75,82.75)(.033613445,-.073529412){119}{\line(0,-1){.073529412}}
%\end
%\emline(132,71.5)(130,69.5)
\multiput(132,71.5)(-.03333333,-.03333333){60}{\line(0,-1){.03333333}}
%\end
\put(130,69.5){\line(1,0){.25}}
%\emline(132,71.75)(134,69.5)
\multiput(132,71.75)(.03333333,-.0375){60}{\line(0,-1){.0375}}
%\end
\put(130.25,69.5){\line(1,0){3.75}}
%\emline(130.25,69.5)(127.75,60.25)
\multiput(130.25,69.5)(-.03333333,-.12333333){75}{\line(0,-1){.12333333}}
%\end
%\emline(132,69.5)(129.5,60.25)
\multiput(132,69.5)(-.03333333,-.12333333){75}{\line(0,-1){.12333333}}
%\end
\put(136.25,62){\line(0,-1){2}}
%\emline(137.5,62)(139.5,60.5)
\multiput(137.5,62)(.04444444,-.03333333){45}{\line(1,0){.04444444}}
%\end
%\emline(132.25,69.5)(136.25,62)
\multiput(132.25,69.5)(.033613445,-.06302521){119}{\line(0,-1){.06302521}}
%\end
%\emline(133.75,69.5)(138,62)
\multiput(133.75,69.5)(.033730159,-.05952381){126}{\line(0,-1){.05952381}}
%\end
\put(44.5,61.5){\line(1,0){1.5}}
\put(44.75,61){\line(1,0){2}}
%\emline(45.25,60.5)(47.5,60.75)
\multiput(45.25,60.5)(.28125,.03125){8}{\line(1,0){.28125}}
%\end
\put(47.5,60.75){\line(1,0){.25}}
%\emline(44.75,48)(47.75,48.5)
\multiput(44.75,48)(.2,.0333333){15}{\line(1,0){.2}}
%\end
%\emline(44.5,48)(47.25,48.25)
\multiput(44.5,48)(.34375,.03125){8}{\line(1,0){.34375}}
%\end
%\emline(44.25,48)(47,48.25)
\multiput(44.25,48)(.34375,.03125){8}{\line(1,0){.34375}}
%\end
\put(44.5,47.75){\line(1,0){2.25}}
\put(44.25,47.5){\line(1,0){1.5}}
\put(46,45.5){$\pi_i$}
\put(45.75,63.25){$\pi_j$}
\put(35.25,50.25){$\cal B$}
\put(44,50.5){$\cal C$}
\put(53.75,57.5){$\Delta_3$}
\put(30.75,58){$\Delta_1$}
%\vector[both](24.25,61.75)(24.25,46.5)
\put(24.25,46.5){\vector(0,-1){.07}}\put(24.25,61.75){\vector(0,1){.07}}\put(24.25,61.75){\line(0,-1){15.25}}
%\end
\put(25.5,53){$\Delta_0$}
%\emline(35.25,56.25)(38.25,43.25)
\multiput(35.25,56.25)(.03370787,-.14606742){89}{\line(0,-1){.14606742}}
%\end
%\emline(36.5,56.25)(39.75,43.5)
\multiput(36.5,56.25)(.033505155,-.131443299){97}{\line(0,-1){.131443299}}
%\end
\put(40.75,57.5){$\Delta_2$}
\put(84.5,55){$\Delta'_2$}
%\dottedline(126.5,56.5)(126.25,56.5)
\multiput(126.43,56.43)(-.125,0){3}{{\rule{.4pt}{.4pt}}}
%\end
%\dottedline(79.75,58)(84,60.25)
\multiput(79.68,57.93)(.70833,.375){7}{{\rule{.4pt}{.4pt}}}
%\end
%\dottedline(78.75,55)(87.5,60.25)
\multiput(78.68,54.93)(.79545,.47727){12}{{\rule{.4pt}{.4pt}}}
%\end
%\dottedline(79.75,53)(91,60)
\multiput(79.68,52.93)(.75,.46667){16}{{\rule{.4pt}{.4pt}}}
%\end
%\dottedline(80.25,51.25)(93.25,59.25)
\multiput(80.18,51.18)(.76471,.47059){18}{{\rule{.4pt}{.4pt}}}
%\end
%\dottedline(80.5,49.5)(94,57.75)
\multiput(80.43,49.43)(.79412,.48529){18}{{\rule{.4pt}{.4pt}}}
%\end
%\dottedline(81.25,48.5)(95.25,56.75)
\multiput(81.18,48.43)(.77778,.45833){19}{{\rule{.4pt}{.4pt}}}
%\end
%\dottedline(83.75,48.25)(95.5,55)
\multiput(83.68,48.18)(.83929,.48214){15}{{\rule{.4pt}{.4pt}}}
%\end
%\dottedline(86.5,47.75)(94.75,53.5)
\multiput(86.43,47.68)(.75,.52273){12}{{\rule{.4pt}{.4pt}}}
%\end
%\dottedline(89,48.5)(95.25,51.5)
\multiput(88.93,48.43)(.78125,.375){9}{{\rule{.4pt}{.4pt}}}
%\end
\thicklines
%\dottedline(61.25,70)(61.5,70)
\multiput(61.18,69.93)(.125,0){3}{{\rule{.8pt}{.8pt}}}
%\end
%\dottedline(61.25,70.25)(61.5,70)
\multiput(61.18,70.18)(.125,-.125){3}{{\rule{.8pt}{.8pt}}}
%\end
\thinlines
\put(61.25,70){\line(1,0){9}}
%\emline(68.5,71)(70.5,70)
\multiput(68.5,71)(.0666667,-.0333333){30}{\line(1,0){.0666667}}
%\end
%\emline(70.5,70)(68.5,69.25)
\multiput(70.5,70)(-.0869565,-.0326087){23}{\line(-1,0){.0869565}}
%\end
\put(103.75,70.25){\line(1,0){10.5}}
\put(114.25,70.25){\line(0,1){0}}
%\emline(114.25,70.25)(112,71.75)
\multiput(114.25,70.25)(-.05,.03333333){45}{\line(-1,0){.05}}
%\end
%\emline(113.75,70.5)(112,69)
\multiput(113.75,70.5)(-.03888889,-.03333333){45}{\line(-1,0){.03888889}}
%\end
\end{picture}

Applying the time separation trick (see Remarks \ref{indiv} and \ref{mini}) with possible decrease of the length of $\cal C$, we may assume that the derivation corresponding
to $\Delta_2$ has no trivial transitions,
in particular, the middle part of this derivation of length
$j-i-1$ has this property too. It corresponds to the trapezium $\Delta'_2$ obtained from $\Delta_2$ by the deleting of the first and the last derivation bands. 

Observe that by Lemma \ref{lens}, $\Delta'_2$ has no cells corresponding to the auxiliary relations (\ref{aux}), and the bottom label of it is of the form $\alpha u q_1\omega$ for some
word $u$ in the tape alphabet of $M_5$. In fact, $u$ is a word in $A_l$ since by the definition
of $\pi_i,$ an
 $a$-band ending on the bottom of $\Delta'_2$  starts on a cell of $\cal C$ having
an edge labeled by $p.$  Therefore $\alpha u q_1\omega$ is a configuration of $M_5$, and
moreover, it is an input configuration, and $\Delta'_2$ is a machine trapezium.

The derivation trapezium $\Delta'_2$
corresponds to a 
computation $C$ of $M_5.$ By Lemma \ref{M} (g), the computation $C$ must be
an input-input computation, since it ends with a configuration  $\alpha u' q_1\omega.$

Thus using (\ref{aux}) we can construct the following  derivation $D$ starting with the bottom label of ${\cal T}_{k+j-i}$ (this word
contains the subword $\alpha u' q_1\omega$):

$$ (\dots \alpha u' q_1\omega\dots)\to (\dots\alpha u' p\dots)\to \dots\to (\dots \alpha p u''\dots)
\to (\dots u''\dots),$$
where $u''$ is the copy of $u'$ in the alphabet $A.$
This property makes possible the following surgery with $\Delta.$ We cut $\Delta$ along
the bottom path  of ${\cal T}_{k+j-i},$ and insert mutually mirror trapezia corresponding to
the derivation $D$ and to its inverse. Since $D$ removes the $\alpha$- and  $q$-letters
in the distinguished subword, this surgery replaces the $\alpha q$-lens $\Gamma$ by an $\alpha q$-lens
with maximal $q$-band of type $2$ and an $\alpha q$-lens with maximal $q$-band of type $t-2$. Since all maximal $q$-bands of $\Delta$,
except for $\cal C$ are untouched by this surgery (more precisely, we added several trivial cells to some of them), the obtained trapezium has
smaller type than $\Delta,$ a contradiction.

{\bf (b)} 
Proving by contradiction, we may assume that the top or bottom of some derivation band contains the subpath $\bf x$.
So its label is of the form $\dots\alpha Up\dots .$ 
Again due to Relations (\ref{aux}), $\alpha U p = U',$ where $U'$ is the copy of $U$ in $A.$ Hence one
can use the same trick as in (a) and replace the through bands $\cal B$ and $\cal C$ by 
a cap and a cup, and the sum of types of their $q$-bands is equal to the type of $\cal C,$ contrary to the minimality of $\Delta.$ 

{\bf (c)} We may assume that $\Gamma$ is a cup. It follows from the
assumptions of the lemma that there are no cells corresponding to Relations (\ref{aux})
surrounded  by $\Gamma$ and the top of $\Delta.$ If the type of $\cal C$ is at least two,
then, as in Case (a), one can consider two $q_1\omega$-cells $\pi_i$ and $\pi_j$
in $\cal C$ and then using similar surgery, replace the cup $\Gamma$
by a cup of smaller type and an $\alpha q$-lens. So our assumption leads to a contradiction with
the minimality of $\Delta.$ 

The lemma is proved. \endproof 

\begin{rk}\label{parts} (a) By Lemma \ref{type2}, the maximal $q$-band $\cal C$ of an $\alpha q$-lens
$E$ in a minimal trapezium $\Delta$ has exactly two $q_1\omega$-cells, say, $\pi_i$ and $\pi_j.$
By Lemma \ref{omega} (a), these two cells are connected (from the right of $\cal C$) by a maximal
$\omega$-band $\cal D.$ We obtain a thick lens $\Gamma$ by adding $\cal D$ to $\Gamma.$ 
The cells of $\cal C$ under $\pi_i$ and above $\pi_j$ (including $\pi_i$ and $\pi_j$)
correspond to the auxiliary relations (\ref{aux}). So the edges of the outer boundary of the thick lens
are either unlabeled or labeled by letters from $A.$ 

(b) One can argue as in Case 2 of
Lemma \ref{type2} (though $t=2$ now) and obtain the subtrapezium $\Delta_0$ and
its parts $\Delta'_2\subset \Delta_2.$ As in the proof of Lemma \ref{type2}, we may assume that $\Delta'_2$
is a machine trapezium, and every derivation band of it corresponds to a (non-trivial) command of an input-input computation of $M_5,$ and
so the top and the bottom of $\Delta_0$ have labels of the form $w'u_lw''$ and $w'v_lw'',$ 
where $u_l$ and $v_l$ are words in the alphabet $A_l,$ and the copies $ u$ and $ v$
of these words in the alphabet $A$ are equal in $S$ by Lemma \ref{M} (a). We will call $\Delta'_2=M(\Gamma)$
the {\it machine part} of 
%$\Delta$ associated with the 
the thick lens $\Gamma;$ $\Delta_2 = \bar M(\Gamma)$ is the {\it augmented machine part} of $\Gamma.$
(It worths to note that it contains  nontrivial $q_1 \omega$-cells in the first and in the last derivation bands.) 
\end{rk}

\begin{lemma} \label{inj} The homomorphism $\phi: S \to H$ defined in Lemma \ref{M} is injective.
\end{lemma}
\proof Let $w$ and $w'$ be two words in generators of $S,$ i.e., in the alphabet $A.$
Assuming that $w=_H w',$ we must prove that $w=_S w'.$ So we have a derivation
$w=w_0\to\dots\to w_t=w'$ over $H$ and denote by $\Delta$ the corresponding minimal derivation trapezium.
Since the boundary labels $w$ and $w'$ of $\Delta$ have neither $\alpha$- nor $q$-letters,
all maximal $\alpha$- and $q$-bands (if any) are paired in some
$\alpha q$-lenses $E_1,\dots, E_k,$ and neither of the corresponding thick lenses $\Gamma_1,\dots, \Gamma_k$ is enclosed in another
one by Lemmas \ref{lens} and \ref{qa} (b). 

We may assume that the first derivation band of $\Delta$ containing a $q_1\omega$-cell from $\cup_{l=1}^k E_l$ (it exists by Lemma \ref{type2} if $k>0$) does contain a cell from $E_1,$ and so
it does not contain other $q_1\omega$-cells. Let $\bar M(\Gamma_1)$ be the augmented machine part of the lens $\Gamma_1$  (see Remark
\ref{parts}).
As in Remark \ref{parts}, we may use the time separation trick, 
 and have each of the  lowest $q_1\omega$-cells of $\Gamma_2\dots \Gamma_k$
 disposed in the derivation bands of $\Delta$ with higher numbers than the derivation bands containing any of the cells of $\Gamma_1.$
 Therefore the time separation trick can now be applied to $\Gamma_2.$ This reconstruction does not touch
 $\bar M (\Gamma_1)$ and creates $\bar M(\Gamma_2)$ with cells disposed 
 above the derivation bands
of $\Delta$ crossing $\bar M(\Gamma_1).$
Finally, we replace $\Delta$ by a minimal trapezium with the same top and bottom
labels, where the augmented machine part $\bar M(\Gamma_i)$ lies above  $\bar M(\Gamma_{i-1})$ for $i=2,\dots,k$.
  
We will keep the same notation $\Delta$ for the obtained trapezium.
In every word $w_i$ of the derivation $w_0\to\dots\to w_t$ corresponding to $\Delta,$
we delete all letters which do not belong to $A\cup A_l,$ replace every letter
from $A_l$ by its copy from $A$ and denote the obtained word from $A^*$ by $\psi(w_i)= W_i.$

By Remark \ref{parts}, $W_r=_S W_s$ if $w_r$ and $w_s$ include, resp., the top and the bottom
labels of some $M(\Gamma_i)$. If $E$ is a trapezium formed by the derivation bands of
$\Delta$ situated between $M(\Gamma_{i-1})$ and $M(\Gamma_{i})$
(or between the bottom (the top) of $\Delta$ and $M(\Gamma_1)$ (and $M(\Gamma_k)$)),
and $W_r$ and $W_s$ are $\psi$-images of the top and the bottom labels of $E,$ then $W_r=W_s,$
because the derivation $w_s\to\dots\to w_r$ uses only the auxiliary relation (\ref{aux}).

Consequently, $W_0=_S W_t,$ and so $w_0=W_0=_S W_t=w_t,$ as required. \endproof 

\begin{lemma}\label{oneq} Let $\Delta$ be a minimal derivation trapezium over $H$ with
the bottom label $w_0=(\alpha) UpV$ and the top label $w_t=(\alpha) U'pV',$ where $U, U'$ are
words in $A_l,$  $V,V'$ are words in $A,$ and $\alpha$ can be absent in both labels. 
We assume that $\Delta$ has a through 
$q$-band $\cal C.$ 
Then using notation of Lemma \ref{inj}, we have $\psi(w_0)=_S \psi(w_t).$
\end{lemma} 
\proof 
Let $\Gamma_1,\dots,\Gamma_k$ be all the thick lenses of $\Delta$ ($k\ge 0$).
By Lemma \ref{lens} (b) none of them is placed from the left of $\cal C.$ 
If the type of $\cal C$ is equal to $2l\ge 0,$ then
it has $2l$ $q_1\omega$-cells, and using these cells one can define $l$ (peeled) augmented machine
trapezia $\bar M_1,\dots \bar M_l,$ where each $\bar M_j$ is bounded from the left by a portion of 
a through $\alpha$-band $\cal B$ (or by the left side of $\Delta$ if there is no
$\alpha$-bands in $\Delta$)
and bounded from the right by an $\omega$-band connecting some  $q_1\omega$-cells
of $\cal C.$ 
%Adding (a part of) one transition band with the $q_1\omega$-cell to the top and (a part of) one %transition band
%to the bottom of each of $M_i$ we have augmented (peeled) machine trapezia $\bar M_1,\dots\bar M_l.$  
Note that by Lemma \ref{qa} (b) and the above observation, these (peeled)  trapezia have no lenses.

Now one can apply the time separation trick to the system $\bar M_1,\dots, \bar M_l,
\bar M(\Gamma_1),\dots, \bar M(\Gamma_k)$ as this was done for the augmented machine parts of
$\Gamma_1,\dots,\Gamma_k$ in the proof of Lemma \ref{inj}. So as there, we
will have $\psi(w_0)=_S\psi(w_t),$ 
%If the top and bottom labels do not contain $\alpha,$ then we can argue as above
%but $M_1,\dots, M_l$ are now peeled machine trapezia. Again we get $\psi(w_0)=\psi(w_t),$
and the lemma is proved. 
\endproof  

\subsection{A-triangles in derivation trapezia}

Assume that $\bf x$ is a nontrivial subpath of the bottom (or of the top) of a trapezium
$\Delta,$ and two vertical paths $\bf y$ and $\bf z$ start  at ${\bf x}_-$ (the original vertex of the path $\bf x$) and
${\bf x}_+$ (the terminal vertex of $\bf x$), resp. If ${\bf y}_+={\bf z}_+,$ and there are no other common vertices of $\bf y$ and $\bf z,$ then we say that $\bf x, y,z$ bound a
triangle subtrapezium $\Delta_0$ of $\Delta.$  (A triangle trapezium    
corresponds to a derivation ending or starting with the empty word $1.$)
If the {\it base \bf x} is labeled by a word in $A,$ we say that $\Delta_0$ is
an $A$-{\it triangle}.

\begin{lemma} \label{capcup} (a) Assume that $\Gamma$ is an  $\alpha q$-cap or an $\alpha q$-cup 
in a minimal trapezium $\Delta,$ and
there are no other $\alpha q$-caps (resp., cups) enclosed in $\Gamma$. Assume that there is an $\alpha q$-lens $E$
enclosed in $\Gamma.$ Then there is an $A$- triangle $\Delta_0$ in $\Delta,$ containing $E$
and enclosed in $\Gamma.$

(b) Let $\Gamma$ be a  triangle in a minimal trapezium $\Delta.$ Assume that
there are no  $\alpha q$-caps or $\alpha q$-cups  but there is an $\alpha q$-lens $E$
enclosed in $\Gamma.$ Then there is an $A$-triangle $\Delta_0$ containing $E$ and enclosed in $\Gamma.$

(c) Assume that a $q$-band $\cal C$ and an $\omega$-band $\cal D$ start (or end) with the
same $q_1\omega$-cell and end (resp. start) on the top (resp., on the bottom) of $\Delta.$
If there are no  $\alpha q$-caps or $\alpha q$-cups  but there is an $\alpha q$-lens  surrounded by these two bands and by the top (by the bottom) of $\Delta,$
then $\Delta$ has an $A$-triangle containing $E.$

(d) Let $\cal C$ be a through $q$-band without $p$-edges in a minimal indivisible
trapezium $\Delta$ and $\cal D$  a through $\omega$-band from the right of $\cal C.$
Suppose there are neither $\alpha q$-caps, nor $\alpha q$-cups, nor through
bands between $\cal C$ and $\cal D,$ but there is an $\alpha q$-lens $E$ between them.
Then there is an $A$- triangle $\Delta_0$ containing $E$ between $\cal C$ and $\cal D.$

%(e)\footnote{Proverit', est' li ssylki na etot punkt} Let $\cal C$ be a through $q$-band in a minimal %indivisible
%trapezium $\Delta$ and $\cal B$  a through $\alpha$-band from the left of $\cal C.$
%Suppose there are neither $\alpha q$-caps, nor $\alpha q$-cups, nor through
%bands between $\cal B$ and $\cal C,$ but there is an $\alpha q$-lens $E$ between them.
%Then there is an $A$- triangle $\Delta_0$ containing $E$ between $\cal B$ and $\cal C.$

\end{lemma}     
\proof {\bf (a)}
We will assume that $\Gamma$ is an $\alpha q$-cap.
Let $O$ be the closed region generated by $E$ (see Subsection \ref{min}).
Note that by Lemma \ref{lens}(c), 
no labeled edge of $\Gamma$ belongs to $O,$ and by Lemma \ref{O} and the definition of $O$,
every labeled edge of the outer boundary of $O$ belongs
to the bottom of $\Delta.$

If $O$ has no edges on the bottom of $\Delta,$ then the (outer) boundary of $O$
has no labeled edges.  
This would contradicts Lemma \ref{unlab}.
So the region $O$ is enclosed in a simple loop $\bf x p,$ where $\bf x$ is the minimal subpath
of the bottom  of $\Delta$ containing all bottom edges belonging to $O,$ and 
$\bf p$ is an unlabeled path on the boundary of $O.$ 

We select a factorization  ${\bf p=zy}^{-1},$ where the last edges of both $\bf y$ and $\bf z$
go upward, and consider  two cases.

{\bf Case 1.} Assume that both $\bf z$ and $\bf y$ are vertical paths. Then $\bf x, y, z$
bound a triangle trapezium $\Delta_0$ . By Lemma \ref{omega} (b), $\bf x$ has no $\omega$-edges
since the maximal $\omega$-band starting on such and edge could not end anywhere.
So every  band starting on $\bf x$ must reach 
an edge of $O$ (or an edge of a closed region generated by another $\alpha q$-lens). 
Therefore the label of this edge belongs to $A$ by Lemma \ref{O}.
Hence $\Delta_0$ is an $A$-triangle.

{\bf Case 2.} One of the paths $\bf y,$ $\bf z,$ say $\bf y$ is not vertical. Then there
is a subpath $eg_1\dots g_lf$ ($l\ge 1$) in $\bf y,$ where the edges $e$ and $f$ go upward, but all 
the edges $g_1,\dots, g_l$ go downward. 

%TeXCAD (http://texcad.sf.net/) Picture. File: [semi8.pic]. Options on following lines.
%\grade{\on}
%\emlines{\off}
%\epic{\off}
%\beziermacro{\on}
%\reduce{\on}
%\snapping{\off}
%\pvinsert{% Your \input, \def, etc. here}
%\quality{8.000}
%\graddiff{0.005}
%\snapasp{1}
%\zoom{4.0000}
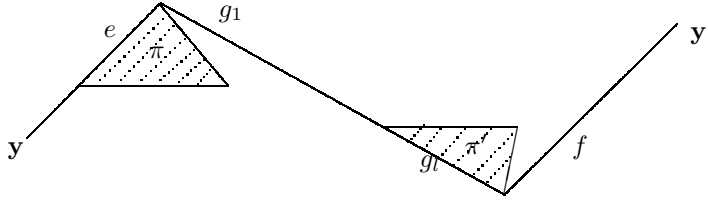
\unitlength 1mm % = 2.845pt
\linethickness{0.4pt}
\ifx\plotpoint\undefined\newsavebox{\plotpoint}\fi % GNUPLOT compatibility
\begin{picture}(117.25,37.25)(0,30)
%\emline(29,49.25)(46.75,67.25)
\multiput(29,49.25)(.0336812144,.0341555977){527}{\line(0,1){.0341555977}}
%\end
%\emline(46.75,67.25)(92.5,41.75)
\multiput(46.75,67.25)(.06051587302,-.03373015873){756}{\line(1,0){.06051587302}}
%\end
%\emline(92.5,41.75)(115.5,64.5)
\multiput(92.5,41.75)(.03407407407,.0337037037){675}{\line(1,0){.03407407407}}
%\end
\put(35.75,56.25){\line(1,0){20}}
%\emline(55.75,56.25)(46.5,67.25)
\multiput(55.75,56.25)(-.0336363636,.04){275}{\line(0,1){.04}}
%\end
\put(76.75,50.75){\line(1,0){17.5}}
\put(94.25,50.75){\line(-1,-5){1.75}}
\put(45.25,60){$\pi$}
\put(87.25,47.25){$\pi'$}
\put(39.25,62.75){$e$}
\put(26.5,47.25){$\bf y$}
\put(117.25,62.5){$\bf y$}
\put(54.5,65.5){$g_1$}
\put(101.5,47.25){$f$}
\put(81.25,45.5){$g_l$}
%\dottedline(38.75,57)(47.75,65.75)
\multiput(38.68,56.93)(.64286,.625){15}{{\rule{.4pt}{.4pt}}}
%\end
%\dottedline(41.5,57)(48.5,64.25)
\multiput(41.43,56.93)(.58333,.60417){13}{{\rule{.4pt}{.4pt}}}
%\end
%\dottedline(44.25,57)(50,62.5)
\multiput(44.18,56.93)(.63889,.61111){10}{{\rule{.4pt}{.4pt}}}
%\end
%\dottedline(46.75,57)(51.25,61)
\multiput(46.68,56.93)(.5625,.5){9}{{\rule{.4pt}{.4pt}}}
%\end
%\dottedline(49.25,56.75)(52.5,59.75)
\multiput(49.18,56.68)(.54167,.5){7}{{\rule{.4pt}{.4pt}}}
%\end
%\dottedline(51.75,56.5)(53.75,58.75)
\multiput(51.68,56.43)(.4,.45){6}{{\rule{.4pt}{.4pt}}}
%\end
%\dottedline(79.75,48.75)(81.75,51)
\multiput(79.68,48.68)(.4,.45){6}{{\rule{.4pt}{.4pt}}}
%\end
%\dottedline(81.75,47.75)(84,50.5)
\multiput(81.68,47.68)(.45,.55){6}{{\rule{.4pt}{.4pt}}}
%\end
%\dottedline(83.75,46.5)(87.25,50.5)
\multiput(83.68,46.43)(.5,.57143){8}{{\rule{.4pt}{.4pt}}}
%\end
%\dottedline(85.75,45.75)(90.75,50.5)
\multiput(85.68,45.68)(.625,.59375){9}{{\rule{.4pt}{.4pt}}}
%\end
%\dottedline(88.5,44.5)(93.5,49.75)
\multiput(88.43,44.43)(.55556,.58333){10}{{\rule{.4pt}{.4pt}}}
%\end
%\dottedline(90.25,43.5)(93.25,46.75)
\multiput(90.18,43.43)(.5,.54167){7}{{\rule{.4pt}{.4pt}}}
%\end
\end{picture}

Since both $e^{-1}$ and $g_1$ are directed downward
and they start from the same vertex, there must be a cell $\pi$ in $\Gamma$ corresponding to the relation
$\alpha p=1$ and having a common edge with the path $eg_1.$ 
Similarly, there is a cell $\pi'$, corresponding to the relation $1=\alpha p$ and having
an edge from the subpath $g_lf.$ Observe that if $\pi$ belongs to $O,$ then $\pi'$ lies
outside this region, and vice versa, since $p$ is a part of the boundary of $O.$ Let us
assume that $\pi$ does not belong to $O.$ However the $\alpha p$-cell $\pi$ must belong to one of the
$\alpha q$-lenses enclosed in $\Gamma.$ This contradicts the definition of the region $O$ since
$\pi$ has an edge belonging to the  boundary of $O.$ Thus Case 2 is impossible,
and Statement (a) is proved. 

{\bf (b,c)} The same proof as for (a) since two sides of the triangle are simply unlabeled now,
and in Case (c), only $a$-band corresponding to the letters from $Y_r$ can start on $\cal C.$  

{\bf (d)} The proof is similar to that in (a) but there might happen that a segment of the boundary
of the closed region $O$ connects the top and the bottom of $\Delta.$ As in Case 2 above, this
segment must be vertical. But this would imply that $\Delta$ is a divisible trapezium, a contradiction. 
\endproof

Let a minimal trapezium $\Delta$ have no $A$-triangles and have an $\alpha q$-cup (or a cap) $\Gamma$ satisfying the assumptions of Lemma \ref{type2} (c). We also assume that $\Delta$ has no cups (or caps) of smaller height than $\Gamma$ (i.e.,
with maximal $q$-band shorter than $\cal C$). 
We define the {\it base label} $b(\Gamma)$ of $\Gamma$ as follows. If the type of the $q$-band
$\cal C$ of $\Gamma$ is $0$ (and so all cells of $\cal C$ are $p$-cells), then $b(\Gamma)$ is just the word $\alpha W p$ we read 
on the top (or on the bottom) of $\Delta$ between the $\alpha$-edge and the $q$-edge of $\Gamma.$
If the type of $\cal C$ is $1,$ then $\cal C$ has one $q_1\omega$-cell, and so one
maximal $\omega$-band $\cal D$ starts on $\cal C$ from the right and, by Lemma \ref{omega} (a), ends on the top %\footnote{Ssylka?}
(or on the bottom) of $\Delta.$ Then $b(\Gamma)$ is the word 
$\alpha U q V\omega$ we read between the ends of $\cal B$ and $\cal D$ on the top/bottom
of $\Delta.$ 

\begin{lemma} \label{cup} Under the above restrictions, (1) if $\cal C$ has type $0,$ then $W$ is a word in the alphabet $A_l;$ 

(2) if the type of $\cal C$ is $1,$ then the base label $b(\Gamma)$ is a reachable configuration
of the machine $M_5.$
\end{lemma}
\proof  If the type of $\cal C$ is $0,$ then every maximal band starting on the segment labeled by $W$ ends on a $p$--cell of $\cal C.$ This implies the first statement of the lemma.
To proof the second one, we consider the derivation bands of $\Delta$ crossing $\cal D.$ They form a derivation
trapezium $\Delta'.$ Let $\Delta''$ be a subtrapezium of  $\Delta'$ bounded from the left
by the left side of $\cal B$ (which is vertical) and bounded from the right by
the right side of $\cal D$ (which is vertical too). 

The bottom of $\Delta''$ 
has label of the form $\alpha W q_1\omega,$ where $W$ is a word in $A_l$ since the underlying
part of the cup $\Gamma$ has the $q$-band of type $0.$ 
The part of $\Delta''$ between $\cal B$ and $\cal C$ has no auxiliary cells (corresponding
to Relations (\ref{aux})) because $\Gamma$ surrounds no $q$-cells. There are no $\alpha q$-cups/caps between $\cal C$ and $\cal D$
by the minimality of the height of $\Gamma.$ Also there are no lenses there
by Lemma \ref{capcup} (d). Hence the part of $\Delta''$ between $\cal C$ and $\cal D$
has no auxiliary cells too. Therefore $\Delta''$ is a machine trapezium with bottom
label $\alpha W q_1\omega,$ and so its top label $b(\Gamma)$ is reachable by $M_5,$
as required.
\endproof

\section{Indivisible trapezia and completion of proofs}
%\footnote{Ne zabyt', chto divisible trapezium mozhet delitsya liniej
%homotopnoj side.} 

\subsection{Upper bounds for spaces of derivations}

{\it We will assume in Lemmas \ref{ao} - \ref{p} that $\Delta$ is an indivisible minimal trapezium
without caps, cups, and $A$-triangles. Let $\Delta$  correspond to a derivation $D:w_0\to\dots\to w_t$ over $H,$
and $\bf x, y$ are the top and the bottom of $\Delta,$ resp. }

\begin{lemma}\label{ao} $\Delta$ has at most one  $\alpha$-band (at most one $\omega$-band) 
connecting $\bf x$ and $\bf y.$ If such a band exists, its left side
(resp., right side)
coincides with the left (resp., right) side of $\Delta.$
\end{lemma}
\proof The letter $\alpha$ (resp., $\omega$) can occur in a defining relation $u=v$
of $H$ only as the left-most (the right-most) letter of $u$ or $v.$ It follows that
the left side of an $\alpha$-band (the right side of an $\omega$-band) connecting $\bf x$ and $\bf y$
is a vertical line. Since $\Delta$ is indivisible, this line must be equal to the
left (to the right) side of $\Delta,$ and the statement of the lemma follows.
\endproof

\begin{lemma}\label{noq} If $\Delta$ has no through $q$-bands, then
$space_H(w_0, w_t)\le S_5(\max(|w_0|_a,|w_t|_a)+3.$
\end{lemma}
\proof We may assume that $space(D)>1.$ Then there are no $\omega$-bands connecting
$\bf x$ and $\bf y$ since either the left or the right side of such a
band would make $\Delta$ divisible. Similarly $\Delta$ has no through $\alpha$-bands
since any $\alpha q$-cell of $\Delta$ must belong to a lens. 

Therefore the top and bottom labels are words in $A\cup Y_l\cup Y_r.$ By Remark \ref{parts} (a),
%\footnote{Skazat' bolee javno ob etom v tom Zamechanii ?}, 
every
$a$-edge of the outer boundary of a thick lens is labeled by a letter from $A.$
%\footnote{Imeetsya v vidu linza s predelannoj $\omega$-ruchkoj sprava.} 
Hence
every maximal $a$-band starting on $\bf y$ with an edge labeled by a letter
from $Y_l\cup Y_r$ consists of trivial cells only and ends on $\bf x.$ This
makes $\Delta$ divisible, a contradiction. 

Thus the top and bottom labels are words in the alphabet $A.$ By Lemma \ref{inj}
we have $w_0=_S w_t,$ and Remark \ref{sp} completes the proof.
\endproof 

\begin{lemma}\label{1q} $\Delta$ has at most one through $q$-band.
\end{lemma}
 \proof
 Assume that there are two through $q$-bands ${\cal C}_1$ and ${\cal C}_2,$
 where ${\cal C}_2$ is from the right of ${\cal C}_1,$ and there are no
 through $q$-bands between them. Note that only maximal $a$-bands corresponding to
 the letters from $Y_l$ can start on the left side of ${\cal C}_2.$ These
 $a$-band can end either on ${\cal C}_2$ or on the top/bottom of $\Delta.$
 Hence there is a vertical path $\bf z$ connecting $\bf x$ and $\bf y$ whose edges
 belong either to the left sides of some of these $a$-bands or to the left side of ${\cal C}_2.$
 It follows that the derivation trapezium $\Delta$ is divisible by this vertical line,
 a contradiction. The lemma is proved.
 \endproof
 
 %TeXCAD (http://texcad.sf.net/) Picture. File: [semi9.pic]. Options on following lines.
%\grade{\on}
%\emlines{\off}
%\epic{\off}
%\beziermacro{\on}
%\reduce{\on}
%\snapping{\off}
%\pvinsert{% Your \input, \def, etc. here}
%\quality{8.000}
%\graddiff{0.005}
%\snapasp{1}
%\zoom{4.0000}
\unitlength 1mm % = 2.845pt
\linethickness{0.4pt}
\ifx\plotpoint\undefined\newsavebox{\plotpoint}\fi % GNUPLOT compatibility
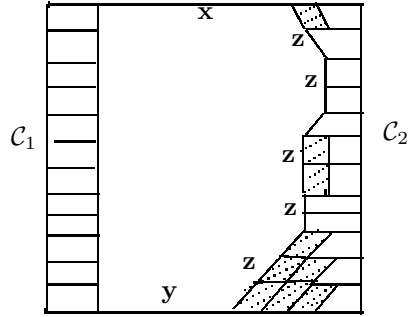
\begin{picture}(105,45.25)(15,17)
\put(102,60){\line(0,-1){41.25}}
%\emline(96,60.25)(97.75,56.75)
\multiput(96,60.25)(.03365385,-.06730769){52}{\line(0,-1){.06730769}}
%\end
\put(97.75,56.75){\line(1,0){4.25}}
\put(94.75,56.75){\line(1,0){3.5}}
%\emline(94.75,56.5)(97.5,52.75)
\multiput(94.75,56.5)(.03353659,-.04573171){82}{\line(0,-1){.04573171}}
%\end
\put(97.5,52.75){\line(1,0){4.5}}
\put(97.25,52.75){\line(0,-1){7.25}}
\put(97.25,45.5){\line(1,0){5}}
\put(97.25,49){\line(1,0){4.75}}
%\emline(97,45.5)(94.5,42.5)
\multiput(97,45.5)(-.03333333,-.04){75}{\line(0,-1){.04}}
%\end
\put(94.5,42.5){\line(1,0){7.25}}
\put(97.75,42.5){\line(0,-1){8}}
\put(97.75,34.5){\line(1,0){4.25}}
\put(94.5,38.75){\line(1,0){7.25}}
\put(94.25,34.5){\line(1,0){4.25}}
\put(94.5,34.75){\line(0,-1){5}}
\put(94.5,29.75){\line(1,0){7.5}}
%\emline(98,29.5)(94.75,26.25)
\multiput(98,29.5)(-.033505155,-.033505155){97}{\line(0,-1){.033505155}}
%\end
\put(94.75,26.25){\line(1,0){7.5}}
\put(94.25,32.25){\line(1,0){7.75}}
%\emline(98.75,26.25)(95.25,22.75)
\multiput(98.75,26.25)(-.033653846,-.033653846){104}{\line(0,-1){.033653846}}
%\end
\put(95.25,22.75){\line(1,0){6.5}}
%\emline(98.75,22.5)(95.75,19)
\multiput(98.75,22.5)(-.03370787,-.03932584){89}{\line(0,-1){.03932584}}
%\end
\put(95.75,19){\line(1,0){6.5}}
\put(60.5,60){\line(1,0){41.75}}
%\emline(92.75,60)(94.5,56.75)
\multiput(92.75,60)(.03365385,-.0625){52}{\line(0,-1){.0625}}
%\end
\put(94.25,42.5){\line(0,-1){7.5}}
%\emline(95.75,23)(92.25,19.25)
\multiput(95.75,23)(-.033653846,-.036057692){104}{\line(0,-1){.036057692}}
%\end
\put(60,19.5){\line(0,-1){.25}}
\put(60,19){\line(1,0){42}}
%\emline(95,26.5)(88.75,19.25)
\multiput(95,26.5)(-.033602151,-.038978495){186}{\line(0,-1){.038978495}}
%\end
%\emline(94.5,30)(85,19.5)
\multiput(94.5,30)(-.0336879433,-.0372340426){282}{\line(0,-1){.0372340426}}
%\end
%\emline(91.25,26.5)(95.75,26.25)
\multiput(91.25,26.5)(.5625,-.03125){8}{\line(1,0){.5625}}
%\end
%\emline(87.75,23)(96.25,23.25)
\multiput(87.75,23)(1.0625,.03125){8}{\line(1,0){1.0625}}
%\end
%\dottedline(93.5,58.5)(95.75,59.75)
\multiput(93.43,58.43)(.5625,.3125){5}{{\rule{.4pt}{.4pt}}}
%\end
%\dottedline(94.5,57.5)(96.5,58.5)
\multiput(94.43,57.43)(.5,.25){5}{{\rule{.4pt}{.4pt}}}
%\end
%\dottedline(96,57)(97.25,57.5)
\multiput(95.93,56.93)(.625,.25){3}{{\rule{.4pt}{.4pt}}}
%\end
%\dottedline(94,40.5)(97.25,42.25)
\multiput(93.93,40.43)(.65,.35){6}{{\rule{.4pt}{.4pt}}}
%\end
%\dottedline(94.5,38.75)(97.75,40.75)
\multiput(94.43,38.68)(.65,.4){6}{{\rule{.4pt}{.4pt}}}
%\end
%\dottedline(95.25,37)(97.5,38.5)
\multiput(95.18,36.93)(.5625,.375){5}{{\rule{.4pt}{.4pt}}}
%\end
%\dottedline(94.5,35)(97.5,36.75)
\multiput(94.43,34.93)(.6,.35){6}{{\rule{.4pt}{.4pt}}}
%\end
%\dottedline(97.25,34.75)(97.25,35.25)
\multiput(97.18,34.68)(0,.25){3}{{\rule{.4pt}{.4pt}}}
%\end
%\dottedline(92.5,27.5)(96.25,30)
\multiput(92.43,27.43)(.75,.5){6}{{\rule{.4pt}{.4pt}}}
%\end
%\dottedline(90.25,24.75)(96.75,28.5)
\multiput(90.18,24.68)(.8125,.46875){9}{{\rule{.4pt}{.4pt}}}
%\end
%\dottedline(86.5,22)(95.25,26.25)
\multiput(86.43,21.93)(.79545,.38636){12}{{\rule{.4pt}{.4pt}}}
%\end
%\dottedline(86,20.25)(96.75,26.25)
\multiput(85.93,20.18)(.82692,.46154){14}{{\rule{.4pt}{.4pt}}}
%\end
%\dottedline(87.5,19.75)(98,25.5)
\multiput(87.43,19.68)(.80769,.44231){14}{{\rule{.4pt}{.4pt}}}
%\end
%\dottedline(89.5,19.5)(96.25,23.75)
\multiput(89.43,19.43)(.75,.47222){10}{{\rule{.4pt}{.4pt}}}
%\end
%\dottedline(92.25,19.75)(96.75,22.75)
\multiput(92.18,19.68)(.64286,.42857){8}{{\rule{.4pt}{.4pt}}}
%\end
%\dottedline(94.25,19.75)(98,22.5)
\multiput(94.18,19.68)(.75,.55){6}{{\rule{.4pt}{.4pt}}}
%\end
%\dottedline(95.5,20)(96.75,20.5)
\multiput(95.43,19.93)(.625,.25){3}{{\rule{.4pt}{.4pt}}}
%\end
%\dottedline(60.5,60.25)(60.5,59.75)
\multiput(60.43,60.18)(0,-.25){3}{{\rule{.4pt}{.4pt}}}
%\end
\put(60.25,60){\line(0,-1){40.75}}
\put(67,60){\line(0,-1){40.75}}
\put(55.5,41.25){${\cal C}_1$}
\put(105,42){${\cal C}_2$}
\put(80.25,58){$\bf x$}
\put(75.5,21){$\bf y$}
\put(92.75,54.5){$\bf z$}
\put(94.5,49){$\bf z$}
\put(91.5,39){$\bf z$}
\put(91.75,32){$\bf z$}
\put(86.25,24.75){$\bf z$ }
\put(60.25,22.75){\line(1,0){6.75}}
\put(60.25,25.75){\line(1,0){6.75}}
\put(60.5,29.25){\line(1,0){6.75}}
\put(60.25,32){\line(1,0){6.75}}
\put(60.25,34.75){\line(1,0){7}}
\put(60.25,38.25){\line(1,0){6.25}}
\put(61.25,41.75){\line(1,0){5.5}}
\put(60.5,45.25){\line(1,0){6.5}}
\put(60.5,49){\line(1,0){6.25}}
\put(60.25,52.25){\line(1,0){6.5}}
\put(60.5,56.5){\line(1,0){6.5}}
\end{picture}

 \begin{lemma}\label{vs} Assume that $\Delta$ has one through $q$-band $\cal C.$ Then 
 
 (1) each cell from the left of $\cal C$ is a trivial $a$-cell corresponding to
a letter from $Y_l$ or it is an $\alpha$-cell;  

(2) if $\Delta$ 
 has no through $\omega$-bands, 
 and  $\cal C$ has no cells having an edge labeled by $p,$  
 then every $q$-cell of $\Delta$ belongs to $\cal C$ and all $a$-edges from
 the right of $\cal C$ are labeled by letters from $Y_r.$ 
 %(resp., from $A$).
  
 \end{lemma}
 \proof We will prove Statement (2) of the lemma. By the assumptions, 
no $\omega$-bands can start/end on $\cal C.$  
Let ${\cal A}_1,\dots$ be all the 
maximal $a$-bands starting on $\cal C$ from
the right. Since no cell of $\cal C$ has a $p$-edge, all these bands correspond
to  letters from $ Y_r.$ Hence each of ${\cal A}_i$-s must end either on
$\cal C$ or on $\bf x,$ or on $\bf y$ (but cannot end on the outer boundary of a thick lens by Remark \ref{parts} (a)).
By Lemmas \ref{qa} (a) and \ref{capcup} (b), there are neither lenses nor $\omega$-cells between any two of these $a$-bands or between some ${\cal A}_i$ and
$\cal C.$ 

There is a vertical line composed of the side edges of these $a$-bands and
of $\cal C$. Since $\Delta$ is indivisible, this line coincides with the right side
of $\Delta,$ and so every cell in and from the right of $\cal C$ corresponds to a machine
relation or is trivial; the trivial $a$-cells are labeled by letters from $Y_r.$

Similarly, each cell from the left of $\cal C$ is a trivial $a$-cell corresponding to
a letter from $Y_l$ or it is an $\alpha$-cell.
\endproof
 
 \begin{lemma}\label{aq} 
 Under the assumptions of Lemma \ref{vs} (2),  $space_H(w_0,w_t)\le
 c_3\max(|w_0|, |w_t|)+c_4,$ where the constants $c_3,c_4$ do not depend on the derivation $D.$
 \end{lemma} 
\proof It follows from Lemma \ref{vs} that all the cells of $\Delta$ correspond to the
machine relations (\ref{machine}), and the derivation $D$ is a peeled
machine derivation,
where the letter $\omega$ is not involved in the commands
of the corresponding computation $\cal C.$ Also there is a reduced computation $C': w_0\to\dots\to w_{t'}=w_t.$ 
Then $t'$ is bounded by a linear function of $\min(|w_0|, |w_t|)$
by Lemma \ref{M} (e). Since the set of defining relations of $H$ is finite, this
implies that the space of $C',$ 
is bounded by a linear function of $\max(|w_0|, |w_t|).$ \endproof
 
\begin{lemma}\label{aqo} Assume that $\Delta$ has one through $q$-band $\cal C,$
 has a through $\omega$-band $\cal D$, 
 and $\cal C$ has no cells having an edge labeled by $p$. Then 
 $$space_H(w_0, w_t)\le\max(c_3|w_0|+c_4,\; c_3|w_t|+c_4,\; S'_5(max(|w_0|_a,|w_t|_a)))$$
 \end{lemma} 
\proof By Lemma \ref{ao} a through $\alpha$-band $\cal B$ (if any exists)  consists of the left-most
cells of $\Delta,$ and
by Lemma \ref{capcup} (d), there are no lenses between $\cal C$ and $\cal D.$
Therefore, as in the proof of Lemma \ref{aq}, we obtain that $\Delta$ is a machine or
a peeled machine trapezium (depending on the presence of a through $\alpha$-band in it).
If it is peeled machine, then $space_H(w_0,w_t)\le 
 c_3\max(|w_0|, |w_t|)+c_4,$ as in Lemma \ref{aq}. If $\Delta$ is a machine trapezium,
 then the statement of the lemma follows from the definition of the function $S'_5(n).$
\endproof

\begin{lemma}\label{p} Assume that $\Delta$ has one through $q$-band $\cal C,$
 and $\cal C$ has an edge labeled by $p$. Then \\ $space_H(w_0,w_t)\le S'_5(\max(|w_0|, |w_t|)+5.$
 \end{lemma} 
 
 \proof By Lemma \ref{vs} (1), each cell from the left of $\cal C$ is a trivial $a$-cell corresponding to
a letter from $Y_l$ or it is an $\alpha$-cell. 

Assume that an $\omega$-band $\cal D$ starts and ends on $\cal C$ (from the right). 
%Let us remove the first and the last cells from $\cal D$ (which are $q_1\omega$-cells) and denote
 % by $\cal D'$ the remaining $\omega$-band. 
% Let $\bf z$ be the (vertical) right side of $\cal D$ extended by the right sides of two %$q_1\omega$-cells at the ends of $\cal D.$
 Let $\Gamma$ be the subtrapezium of $\Delta$ bounded from the right by the right side $\bf z$ of $\cal D$ and bounded from
 the left by a part of the left side of $\Delta.$ By Lemma \ref{qa} (b), $\Gamma$ is a (peeled) augmented machine
 trapezium corresponding to a computation of $M_5$. Applying the type separation trick, we may
 assume that  all cells from the right of $\bf z$ are trivial, and the computation of $M_5$ is reduced.
 Moreover, it is non-empty since otherwise the type of $\cal C$ could be decreased after  removing of two derivation bands containing the $q_1\omega$-cells of $\Gamma.$
 
 If $\Gamma$ has no cell having both an $\alpha$-edge and a $q$-edge, then we have a contradiction with Lemma \ref{M} (f). Suppose there is such a cell $\pi,$ and let us choose it 
to be the closest one to the bottom of $\Gamma.$ Then by Lemma \ref{M} (f), the part of $\cal C$ between
the bottom of $\Gamma$ and $\pi$ has no edges labeled by letters from $Y_l\backslash A_l.$
Hence only $a$-bands corresponding to letters from $A_l$ can start on the bottom of $\Gamma$ and
end on $\cal C$ from the left. Therefore the $a$-letters from the  bottom label of $\Gamma$ belong to  $A_l$
contrary to Lemma \ref{type2} (b). 
 So we may assume that no $\omega$-band starts and ends on $\cal C.$
 
 %TeXCAD (http://texcad.sf.net/) Picture. File: [semi10.pic]. Options on following lines.
%\grade{\on}
%\emlines{\off}
%\epic{\off}
%\beziermacro{\on}
%\reduce{\on}
%\snapping{\off}
%\pvinsert{% Your \input, \def, etc. here}
%\quality{8.000}
%\graddiff{0.005}
%\snapasp{1}
%\zoom{4.0000}
\unitlength 1mm % = 2.845pt
\linethickness{0.4pt}
\ifx\plotpoint\undefined\newsavebox{\plotpoint}\fi % GNUPLOT compatibility
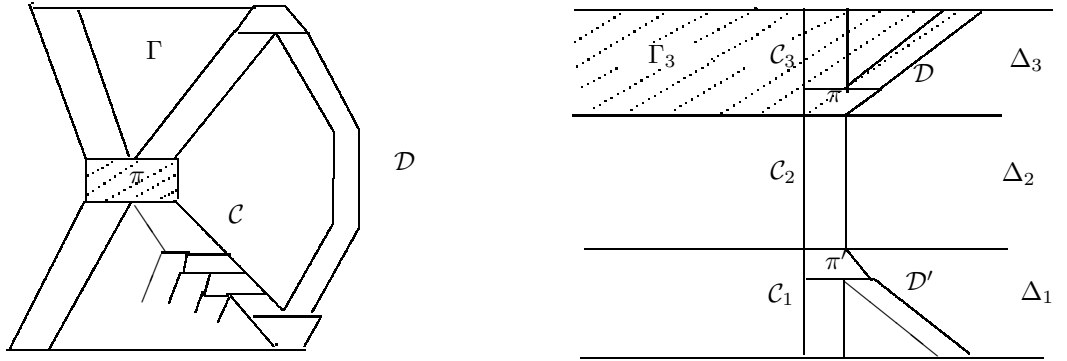
\begin{picture}(177.5,60)(35,20)
\put(38.5,28.75){\line(1,0){39.75}}
\put(49,54.25){\line(0,-1){5.75}}
\put(49,48.5){\line(1,0){12}}
%\emline(61,48.5)(60.75,48.75)
\multiput(61,48.5)(-.03125,.03125){8}{\line(0,1){.03125}}
%\end
\put(49,54.25){\line(1,0){12.25}}
\put(61.25,54.25){\line(0,-1){5.25}}
\put(74.25,29.5){\line(1,0){.25}}
\put(71.25,33.25){\line(1,0){9}}
%\emline(80.25,33.25)(78,29.25)
\multiput(80.25,33.25)(-.03358209,-.05970149){67}{\line(0,-1){.05970149}}
%\end
%\emline(60.75,48.75)(75.25,34)
\multiput(60.75,48.75)(.0337209302,-.0343023256){430}{\line(0,-1){.0343023256}}
%\end
%\emline(48.75,48.5)(38.75,28.75)
\multiput(48.75,48.5)(-.0336700337,-.0664983165){297}{\line(0,-1){.0664983165}}
%\end
%\emline(55,48.5)(44,28.75)
\multiput(55,48.5)(-.0336391437,-.0603975535){327}{\line(0,-1){.0603975535}}
%\end
%\emline(68.25,36)(74,29.25)
\multiput(68.25,36)(.033625731,-.039473684){171}{\line(0,-1){.039473684}}
%\end
\put(68,36){\line(-1,0){3.25}}
%\emline(64.75,36)(65.5,39)
\multiput(64.75,36)(.0326087,.1304348){23}{\line(0,1){.1304348}}
%\end
\put(61.75,39){\line(1,0){3.75}}
%\emline(62,39.25)(62.5,41.75)
\multiput(62,39.25)(.0333333,.1666667){15}{\line(0,1){.1666667}}
%\end
\put(59,41.75){\line(1,0){3.75}}
\put(55.5,48){\line(2,-3){4}}
\put(59,41.5){\line(-2,-5){2.5}}
%\emline(61.5,39)(60,35)
\multiput(61.5,39)(-.03333333,-.08888889){45}{\line(0,-1){.08888889}}
%\end
%\emline(64.5,36.5)(63.5,33.75)
\multiput(64.5,36.5)(-.0333333,-.0916667){30}{\line(0,-1){.0916667}}
%\end
%\emline(68,35.75)(66.5,32.5)
\multiput(68,35.75)(-.03333333,-.07222222){45}{\line(0,-1){.07222222}}
%\end
\put(54.75,51){$\pi$}
%\emline(60.75,54.25)(74.25,71)
\multiput(60.75,54.25)(.0336658354,.0417705736){401}{\line(0,1){.0417705736}}
%\end
\put(74.25,71){\line(0,1){0}}
\put(69.25,71){\line(1,0){9.5}}
%\emline(78.5,71)(75.5,74.5)
\multiput(78.5,71)(-.03370787,.03932584){89}{\line(0,1){.03932584}}
%\end
\put(75.5,74.5){\line(-1,0){4.25}}
%\emline(71.25,74.5)(68.75,71)
\multiput(71.25,74.5)(-.03333333,-.04666667){75}{\line(0,-1){.04666667}}
%\end
%\emline(68.75,71.25)(55.5,54.25)
\multiput(68.75,71.25)(-.0337150127,-.0432569975){393}{\line(0,-1){.0432569975}}
%\end
\put(42,75){\line(-1,0){.25}}
\put(42,74.25){\line(1,0){29.5}}
%\emline(41.5,74.75)(49,54.25)
\multiput(41.5,74.75)(.033632287,-.091928251){223}{\line(0,-1){.091928251}}
%\end
%\emline(48,74)(54.75,54.5)
\multiput(48,74)(.03358209,-.097014925){201}{\line(0,-1){.097014925}}
%\end
%\emline(74.25,71)(82,57.75)
\multiput(74.25,71)(.033695652,-.057608696){230}{\line(0,-1){.057608696}}
%\end
%\emline(82,57.75)(81.75,45.5)
\multiput(82,57.75)(-.03125,-1.53125){8}{\line(0,-1){1.53125}}
%\end
%\emline(81.75,45.5)(75.25,34)
\multiput(81.75,45.5)(-.033678756,-.059585492){193}{\line(0,-1){.059585492}}
%\end
%\emline(78.5,71)(85.25,58.25)
\multiput(78.5,71)(.03358209,-.063432836){201}{\line(0,-1){.063432836}}
%\end
\put(85.25,58.25){\line(0,-1){13.25}}
%\emline(85.25,45)(79,33.75)
\multiput(85.25,45)(-.033602151,-.060483871){186}{\line(0,-1){.060483871}}
%\end
\put(90,52.75){$\cal D$}
\put(68,45.75){$\cal C$}
\put(57,67.5){$\Gamma$}
\put(62.25,41.5){\line(1,0){6}}
\put(65,39){\line(1,0){5}}
\put(67.75,36.25){\line(1,0){5.25}}
%\dottedline(49.25,52)(52.5,53.75)
\multiput(49.18,51.93)(.65,.35){6}{{\rule{.4pt}{.4pt}}}
%\end
%\dottedline(49.25,49.5)(56.5,54.25)
\multiput(49.18,49.43)(.725,.475){11}{{\rule{.4pt}{.4pt}}}
%\end
%\dottedline(51.75,48.75)(61.25,54)
\multiput(51.68,48.68)(.79167,.4375){13}{{\rule{.4pt}{.4pt}}}
%\end
%\dottedline(55.75,49.25)(60.75,52)
\multiput(55.68,49.18)(.71429,.39286){8}{{\rule{.4pt}{.4pt}}}
%\end
%\dottedline(59,49)(61,50)
\multiput(58.93,48.93)(.5,.25){5}{{\rule{.4pt}{.4pt}}}
%\end
%\dottedline(138,74.25)(138.25,74)
\multiput(137.93,74.18)(.125,-.125){3}{{\rule{.4pt}{.4pt}}}
%\end
\put(144.5,74.25){\line(0,-1){46}}
\put(144.5,60){\line(0,1){.25}}
\put(113.75,60){\line(1,0){57}}
\put(115.25,42.25){\line(1,0){56.25}}
\put(150.25,74.25){\line(0,-1){11.25}}
\put(144.5,63.5){\line(1,0){10.25}}
%\emline(154.75,63.5)(150.25,60.25)
\multiput(154.75,63.5)(-.046391753,-.033505155){97}{\line(-1,0){.046391753}}
%\end
\put(150,60.25){\line(0,-1){18.25}}
%\emline(150.25,42)(153.5,38)
\multiput(150.25,42)(.033505155,-.041237113){97}{\line(0,-1){.041237113}}
%\end
\put(153.5,38.25){\line(-1,0){8.75}}
\put(149.75,38){\line(0,-1){10}}
%\emline(150.25,64)(163,74)
\multiput(150.25,64)(.0429292929,.0336700337){297}{\line(1,0){.0429292929}}
%\end
%\emline(154.5,63.5)(168,73.75)
\multiput(154.5,63.5)(.0444078947,.0337171053){304}{\line(1,0){.0444078947}}
%\end
\put(149.75,38){\line(5,-4){12.5}}
%\emline(153.75,38.25)(166.5,28.25)
\multiput(153.75,38.25)(.0429292929,-.0336700337){297}{\line(1,0){.0429292929}}
%\end
\put(114,74){\line(1,0){63.5}}
\put(114.75,27.75){\line(1,0){61.5}}
\put(147.5,61.75){$\pi$}
\put(147.25,39.5){$\pi'$}
\put(159,64.5){$\cal D$}
\put(158,36.75){$\cal D'$}
\put(139.75,35.25){${\cal C}_1$}
\put(140,51.5){${\cal C}_2$}
\put(140,67){${\cal C}_3$}
\put(173.25,35.5){$\Delta_1$}
\put(170.75,51.25){$\Delta_2$}
\put(171.75,66.5){$\Delta_3$}
\put(123.75,67){$\Gamma_3$}
%\dottedline(114,69.75)(120,73.5)
\multiput(113.93,69.68)(.75,.46875){9}{{\rule{.4pt}{.4pt}}}
%\end
%\dottedline(114.25,66.5)(125.75,73.5)
\multiput(114.18,66.43)(.76667,.46667){16}{{\rule{.4pt}{.4pt}}}
%\end
%\dottedline(115,63.75)(119.5,66.25)
\multiput(114.93,63.68)(.75,.41667){7}{{\rule{.4pt}{.4pt}}}
%\end
%\dottedline(124.5,69.75)(131.5,73.5)
\multiput(124.43,69.68)(.77778,.41667){10}{{\rule{.4pt}{.4pt}}}
%\end
%\dottedline(116.5,61)(125.25,66)
\multiput(116.43,60.93)(.79545,.45455){12}{{\rule{.4pt}{.4pt}}}
%\end
%\dottedline(130,69.75)(137,73.5)
\multiput(129.93,69.68)(.77778,.41667){10}{{\rule{.4pt}{.4pt}}}
%\end
%\dottedline(121.5,60.75)(132.5,67.25)
\multiput(121.43,60.68)(.78571,.46429){15}{{\rule{.4pt}{.4pt}}}
%\end
%\dottedline(137.5,70.25)(143,74)
\multiput(137.43,70.18)(.78571,.53571){8}{{\rule{.4pt}{.4pt}}}
%\end
%\dottedline(127.25,60.75)(136.5,66.25)
\multiput(127.18,60.68)(.77083,.45833){13}{{\rule{.4pt}{.4pt}}}
%\end
%\dottedline(143,70.25)(149.75,73.5)
\multiput(142.93,70.18)(.84375,.40625){9}{{\rule{.4pt}{.4pt}}}
%\end
%\dottedline(133.5,60.25)(154.5,73.75)
\multiput(133.43,60.18)(.80769,.51923){27}{{\rule{.4pt}{.4pt}}}
%\end
%\dottedline(139,60.5)(159,73.75)
\multiput(138.93,60.43)(.76923,.50962){27}{{\rule{.4pt}{.4pt}}}
%\end
%\dottedline(145.25,60.5)(165.25,74)
\multiput(145.18,60.43)(.76923,.51923){27}{{\rule{.4pt}{.4pt}}}
%\end
%\dottedline(132.25,67.5)(138.25,70.75)
\multiput(132.18,67.43)(.75,.40625){9}{{\rule{.4pt}{.4pt}}}
%\end
\end{picture}

  There is at most one $\omega$-band $\cal D$ starting with a cell $\pi$ of $\cal C$ and ending on the top of $\Delta$ and at most one $\omega$-band $\cal D'$ starting on the bottom of $\Delta$ and ending with a cell $\pi'$
 of $\cal C.$ (We take into account that the $p$-cells corresponding to the relation $p=q_1\omega$ and to $q_1\omega=p$ alternate in $\cal C.$) 
 Since $\cal C$ has a $p$-cell, $\Delta$ must consists of
 the following pieces enumerated from the bottom to the top (some of them may be absent):
 the subtrapezium $\Delta_1$ crossed by the 
 the part ${\cal C}_1$ of $\cal C,$ connecting the bottom of $\Delta$ and $\pi',$
 %$\omega$-band ${\cal D}_1$ which ends with $\pi'$,
 %the transition band $\Delta_2$ of $\Delta$ containing the cell $\pi'$ of $\cal C,$ 
 the subtrapezium
 $\Delta_2$ crossed by the part ${\cal C}_2$ of $\cal C,$ containing only $p$-cells corresponding
 to the auxiliary relations involving letters from $A,$  
 %the transition band $\Delta_4$ containing the cell $\pi,$
 and the subtrapezium $\Delta_3$  crossed by the part ${\cal C}_3$ of $\cal C,$
 starting with $\pi$ and ending on the top of $\Delta.$
 %the $\omega$-band ${\cal D}_2$ starting with $\pi.$
 The band $\cal C$ cannot share a cell with a through $\omega$-band,
 and so there are no through $\omega$-bands in $\Delta$ since $\Delta$ is indivisible. 
 
 Let 
 %${\cal C}_2$ be the part of $\cal C$ crossing $\Delta_5$ and 
 $\Gamma_3$ is the part of $\Delta_3$
 bounded from the right by the right side   ${\cal D}.$ By Lemma \ref{capcup} (c), $\Gamma_3$ is
 a (peeled) augmented machine trapezium. 
 %Adding the part of the transition band with the $q_1\omega$-cell,
 %we have the augmented trapezium $\bar \Gamma_2$
 
 If ${\cal C}_3$ has no $\alpha$-cells, then
 every $a$-band starting on ${\cal C}_3$ from the left, corresponds to the $a$-letter from $A_l$
 by Lemma \ref{M} (f). If ${\cal C}_3$ has an $\alpha$-cell, then one can choose such a cell to
 be the closest to the bottom of $\Gamma_3$ and, as above, obtain a contradiction with Lemma \ref{type2} (b). 
 Similarly, every $a$-band starting on ${\cal C}_1$  from the left, corresponds to the $a$-letter from $A_l,$ and therefore the same property holds for the whole $\cal C.$ 
 
 It follows that there are no $a$-edges corresponding to the letters of $Y_l\backslash A_l$ from
 the left of $\cal C$ since the maximal $a$-band through such edges would make the trapezium $\Delta$
 divisible. Moreover, we see that  
 if $\Delta$ has an $\alpha$-band $\cal B$, then no $\omega$-band starts/ends on $\cal C.$
 Now we consider two cases. 
 
  {\bf Case 1}: $\Delta$ has no through $\alpha$-bands. 
 
 One can apply the time separation trick to $\Delta_3$ and its 'left half' $\Gamma_3.$ Therefore
 one can assume that
 %the machine parts of 
 the lenses
 from the right of $\cal C$ (if they exist) lie in $\Delta_2.$ Furthermore, if $\Delta_2$ corresponds
 to a subderivation $w_i\to\dots \to w_j$ of $D,$ then $|w_0|= |w_1|=\dots =|w_i|$ and
 $|w_j|=\dots = |w_t|$ by Lemma \ref{M} (f). Hence it suffices to estimate $space_H(w_i, w_j).$
 
 We have $w_i=UpV,$ where $U$ is a word in $A_l$ and $V$ is a word in $A.$ Indeed, 
 every maximal $a$-band $\cal A$ of $\Delta_2$ disposed from the right of $\cal C$ and corresponding to a non-$A$ $a$-letter cannot end either on the $p$-cell of $\cal C$ or on an outer boundary of a thick lens of $\Delta_2.$ So both sides of it divide $\Delta_2,$ and the sides of the maximal extension of $\cal A$ in the whole $\Delta$ divide $\Delta,$ a contradiction. Similar form has $w_j=U'pV'.$ Hence
 $\psi(w_i)=\psi(w_j)$ by Lemma \ref{oneq}. 
 
 To complete the proof, we first use the relations $a_lp=pa$ to replace $w_i$ by $w'_i = pU_AV,$ where
 $U_A$ is the copy of $U$ in the alphabet $A.$ The derivation $w_i\to\dots\to w'_i$ has space $|w_i|=|w'_i|.$
 Similarly we obtain $w'_j.$ By Remark \ref{sp}, $space_H(U_AV,U'_AV')\le S'_5(\max(|U_AV|,|U'_AV'|)+3.$
 Hence $$space_H(w_i,w_j)= space_H(w'_i,w'_j)\le S'_5(\max(|w_0|,|w_t|)+4$$
 
 {\bf Case 2}:  $\Delta$ has a through $\alpha$-band $\cal B$ (as in Lemma \ref{ao}).
 Since there are no $q_1\omega$-cells in $\cal C,$ 
 every cell of this band is a cell corresponding to a relation $pa=a_lp$ or a trivial
 $p$-cell.
 In particular,
 $\cal C$ does not share cells with $\omega$-bands, and therefore by Lemma \ref{ao}, 
 $\Delta$ has no through $\omega$-bands since both sides  of such a band would be
 vertical, but $\Delta$ is indivisible. 
 
 Since every $a$-band starting/ending on $\cal C$ from the right or on the outer boundary
 of a thick lens, correspond to a letter from $A,$ the top and bottom labels of $\Delta$
 from the right of $\cal C$ are the words in the alphabet $A$ (again, because $\Delta$
 indivisible). 
 So $w_0=\alpha UpV,$ $w_t=\alpha U'pV',$ where $V, V'$ are words in $A$ (and $U,U'$ are words
 in $A_l$).

 Hence, by Lemma \ref{oneq}, $\psi(w_0)=_S\psi(w_t).$
 The relations $pa=a_lp$ preserve the value of $\psi$ and does not change the length.
 So we may assume that $w_0=\alpha pV, w_t=\alpha pV',$ where $V$ and $V'$ are words in $A.$
 By Remark \ref{sp}, the words $V$ and $V'$ can be connected by an $H$-derivation of space
 at most $S'_5(\max(|V|,|V'|)+3.$ Therefore $space_H(w_0,w_t)\le S'_5(\max(|w_0|_a, |w_t|_a)+5,$
 as required.
 
 \endproof
 
 Summarizing, we obtain
 \begin{lemma}\label{sum} Assume that $\Delta$ is an indivisible minimal trapezium
without caps, cups, and $A$-triangles, and $\Delta$  corresponds to a derivation $D:w_0\to\dots\to w_t$ over $H.$ Then $$space_H(w_0, w_t)\le max(c_3|w_0|+c_4, c_3|w_t|+c_4, S'_5(max(|w_0|_a,|w_t|_a))+5)$$
 \end{lemma} 
 
 \proof By Lemmas \ref{noq} and \ref{1q}, we may assume that $\Delta$ has exactly one
 through $q$-band $\cal C.$ If $\cal C$ has no $p$-cells, then the statement of the 
 lemma follows from Lemmas \ref{aq} and \ref{aqo}. Otherwise it
 follows from Lemma \ref{p}.
 \endproof
 
 Now we want to eliminate the restrictions of Lemma \ref{sum}
 imposed on $\Delta.$
 
 \begin{lemma}\label{final} The space function of $H$ is bounded
 from above by a function 
 equivalent to the function $S'_5(n).$
 \end{lemma} 
 \proof
 We define the function $f(n)= S'_5(n+c)+c_3n+c_4+5$ for $n\ge 1,$ 
 where $c$ is the constant from Lemma \ref{M} (d), $c_3, c_4$ are from Lemma \ref{aq} (and \ref{sum}),
 and define $f(0)=0$ $(=S'_5(0)).$
 Obviously, $f(n)\sim S'_5(n).$ We can use the inequality $f(n-k)+k\le f(n)$ for $0\le k\le n.$ 
Indeed the function $S'(n)$ is
non-descending, and one can select $c_3\ge 1.$
 
 Now we modify the length of a word $w:$ by definition $||w||$ is the number
 of letters, where every $\alpha$- or $q$-letter is counted with weight
 $c,$  and other letters are counted with weight $1.$ 
 
 Note that $|w|\le ||w||\le c|w|$ for every word in the generators of $H.$ Therefore to prove the lemma, it suffices to prove the inequality $space_H(w, w')\le f(||w||+||w'||)$ for any pair of equal in $H$ words $w$ and $w'.$ This will be proved by induction on $\Sigma=||w||+||w'||$ with trivial base $\Sigma =0.$ So we will assume that $\Sigma>0$ and consider a derivation $D:$ $w=w_0\to\dots\to w_t=w'.$ 
 Let us denote by $\Delta$ the corresponding minimal trapezium. Of course, one may assume that
 the unique vertical line connecting the endpoints of the left side of $\Delta$ is the left side
 itself since otherwise one can replace $\Delta$ by a subtrapezium. Similar assumption is taken for the
 right side of $\Delta.$
 
First assume  that the trapezium $\Delta$ is divisible and use
the notation of Remark \ref{indiv}. Then (see formula (\ref{div})) there is a derivation
$$D': w=w_0(1)w_0(2)\to\dots\to w_t(1)w_0(2)\to\dots\to w_t(1)w_t(2)=w',$$
where  $\max(||w_0(1)||+||w_t(1)||, ||w_0(2)||+||w_t(2)||)< ||w||+||w'||, $ and so by the inductive hypothesis, the first half (the second half) of the
derivation $D'$ can be chosen with space at most $f(||w_0(1)||+||w_t(1)||) +|w_0(2)|$ (resp., at most $f(||w_0(2)||+||w_t(2)||) +|w_t(1)|.$) Hence 
$$space(D')\le \max(f(||w||+||w'||-||w_0(2)||)+||w_0(2)||,
f(||w||+||w'||-||w_t(1)||)+||w_t(1)||)$$ $$
\le f(||w||+||w'||)$$ 
Thus we may further assume that the derivation trapezium  $\Delta$
is indivisible. 

Now assume that there is an $A$-triangle in $\Delta.$ It means that the bottom
(or the top) label of $\Delta$ is of the form $w=\bar w u\bar{\bar w},$ where $u$ is a non-empty
word in the alphabet $A$ and $u=_S 1$ by Lemma \ref{inj}. By Remark \ref{sp}, $space_H (u,1)\le S'_5(|u|)+3.$
Therefore there is a derivation $w\to\dots\to \bar w \bar{\bar w}$ over $H$ of space at most
$$S'_5(|u|)+3+|\bar w|+|\bar{\bar w}|\le f(|u|) +|\bar w|+|\bar{\bar w}|=
f(|w|- |\bar w|-|\bar{\bar w}|)+|\bar w|+|\bar{\bar w}| \le f(|w|)\le f(||w||)$$
By the inductive hypothesis, there is a derivation $\bar w \bar{\bar w}\to\dots\to w'$
of space $\le f(||\bar w \bar{\bar w}||+||w'||) \le f(||w||+||w'||)$, hence
$space_H(w, w')\le f(||w||+||w'||).$ Thus we may further assume that $\Delta$ has no
$A$-triangles.

Assume that $\Delta$ has a cup (or cap). Let $\Gamma$ be a cup of minimal height.
By Lemma \ref{capcup} (a), there are no lenses enclosed in $\Gamma.$ Therefore the
type of the maximal $q$-band ${\cal C}_{\Gamma}$ of $\Gamma$ is $0$ or $1$ by Lemma
\ref{type2} (c).

If the type of $\cal C$ is $0,$ then by Lemma \ref{cup}, the word $w'$ is of the
form $\bar w \alpha W p \bar{\bar w},$ where $W$ is a word in $A_l.$ Hence there is
a derivation $w'\to \dots\to \bar w \alpha  p W' \bar{\bar w}\to \bar w W' \bar{\bar w}$, where $W'$ is the
copy of $W$ in the alphabet $A,$ and this derivation has space $|w'|\le ||w'||.$ Since 
$||\bar w W' \bar{\bar w}||=||w'||-2c,$ we obtain by the inductive hypothesis, that
$$space_H(w, \bar w W' \bar{\bar w})\le f(||w||+||w'||-2c)\le f(||w||+||w'||)$$ Hence
$space_H (w,w')\le f(||w||+||w'||),$ as desired.

If the type of $\cal C$ is $1,$ then, by Lemma \ref{cup}, the word $B=b(\Gamma)$ is reachable by the machine $M_5.$ By Lemma \ref{M} (d), there is a computation $B\to\dots\to B'$ of $M_5,$
where $B'$ is an input configuration of $M_5$ and $|B'|_a \le |B|_a+c,$ and so $||B'|| \le ||B||+c$. Denote by $C'$
the corresponding machine derivation over $H.$ Furthermore
applying the auxiliary relations in the standard way, one can extend $C',$ 
remove the letters $\alpha,$ $q_1$ and $\omega$ from $B'$ and obtain a word $B''$ with $||B''||< ||B'||-2c < ||B||-c.$
Hence we have a derivation $w'=\bar w B \bar{\bar w}\to\dots \to \bar w B'' \bar{\bar w}$   
of space at most $S'_5(|B|+c)+|w'|-|B|.$ By the inductive hypothesis, there
is a derivation $w\to \dots\to \bar w B'' \bar{\bar w}$ of space at most $f(||w||+||w'||-c).$
Therefore $space_H(w,w')\le f(||w||+||w'||),$ as required. 

Thus we may assume that $\Delta$ satisfies the assumptions of Lemma \ref{sum},
and therefore $space_H(w,w')\le f(||w||+||w'||)$ by Lemma \ref{sum} and the definition of $f(n).$
The Lemma is proved.

\subsection{Lower bounds and  completion of proofs}

We define a monoid $H'$  as follows. 
The set of generators of $H'$ is $A_{H'}= Y_l \sqcup Y_r\sqcup Q\sqcup\{\alpha, \omega\},$
i.e., $A_{H'}=A_H\backslash\{A\cup\{p\}\}$).
The set of defining relations of $H'$ consists of only machine relations of $H:$ 
$R_{H'}=\{V'=V \;\;for\; every\; command\;\; V\to V'\; of\; M_5\}.$
\begin{lemma}\label{H'} The space function of $H'$ is 
equivalent to $S'_5(n).$
\end{lemma}
\proof

The space function of $H'$ is bounded from above by a function
equivalent to $S'_5(n).$ This statement is a very easy version of Lemma \ref{final},
since the analogs of Lemmas \ref{ao} - \ref{aqo} and \ref{sum} become 
trivial when we have no letters from $A,$ no $p,$ no defining relations of $H'$ with $1$ in the left/right sides, and consequently, no triangles, lenses, caps and cups.

By the definition of the set of relations $R_{H'},$ every computation of $M_5$ can
be considered as a derivation over $H',$ and vice versa, every derivation $w\to\dots,$ where $w$ is
a configuration of $M_5,$ is a computation.
Let a derivation $C$ over $H'$ connects $w=\alpha U q V \omega$ and  $\alpha U' q' V' \omega,$
where $q,q'\in Q,$ the words $U,U'$ are words in the alphabet $Y_l$ and $V,V'$ are words in $Y_r.$ Then $C$
is a machine derivation of $M_5$ since all defining relations of $H'$ are machine relations. 
Hence if for two configurations $w$ and $w'$ of $M_5,$
%$w=\alpha U q_1 V \omega$ and  $w'\alpha U' q_1 V' \omega,$ 
we have
$space_{M_5}(w,w')= s$ for some $s,$ then $w$ and $w'$ cannot be connected by
a derivation over $H'$ with space $\le s.$ Hence the space function of $H'$ is at least
$S'_5(n),$ and the statement of the lemma follows.
\endproof

{\bf Proof of Theorem \ref{main}.} Let a monoid $H''$ be a copy of $H'$ given by a finite presentation with
a set of generators $A_{H''}$ disjoint with $A_H.$ We define the monoid $P$ announced in Theorem \ref{main}, as the free
product $H\star H''$ and consider its space function $s(n)$ with respect to the presentation
$\la A_H\cup A_{H''}\mid R_H\cup R_{H''}\ra.$

On the one hand, any derivation over $P$ projects on a derivation $C$ over $H''.$ (One just deletes the letters
from $A_H$ in any word from $C.$) Therefore the space function $s(n)$
%\footnote{Space funktsija ne opredelena vo vvedenii (tol'ko v abstracte) !} 
of $P$ is greater than or equal to the space
function of $H''.$ Hence $s(n)\succeq S'_5(n)$ by Lemma \ref{H'}. On the other hand
two words  $w$ and $w$ over $A_H\cup A_{H''}$ are equal in $P$ iff their corresponding
$A_H$- and $A_{H''}$-syllables are equal in $H$ and in $H',$ respectively. Since
the derivations between equal words can be define syllable-by-syllable, we see
that, up to equivalence, $s(n)$ does not exceed the maximum of the space functions
of $H$ and of $H'$. Therefore $s(n)\preceq S'_5(n)$ by Lemmas \ref{final} and \ref{H'}.

Our estimates show that $s(n)\sim S'_5(n).$ Recall that $S'_5(n)\sim S_0(n)$ by
Lemma \ref{M} (b). Hence $s(n)\sim S_0(n)$, and by Lemma \ref{inj}, the theorem is proved.$\Box$
\medskip   

{\bf Proof of Corollary \ref{pspace}}. If the function $S_0(n)$ is bounded by a
polynomial, then so is $S'_5(n)$ by Lemma \ref{M} (b). By Lemma \ref{inj}, $S$
is a subsemigroup (submonoid) of the monoid $H,$ and the space function of $H$
is polynomial by Lemma \ref{final}.

Conversely, assume that a finitely generated semigroup (monoid) $S$ is
embedded in a finitely presented semigroup (monoid) $H$ with polynomial space
function $s(n)$. Then the word problem in $S$ is solvable by an NTM with space
$\preceq s(n):$ This machine takes any word $w$ in the generators of $S,$ rewrites it
in the generators  of $H,$ and if $w=_H 1,$ it produces a derivation $w\to\dots\to 1.$
(Recall that an NTM may guess and verify. The head of this NTM can move along a word 
and can replace a subword $u$ by $v$ if $u=v$ or $v=u$ is one of defining relations of $H.$ See details
in \cite{Bir}.)  
By the remarkable theorem of Savitch (see \cite{DK}, Corollary 1.31), if an NTM has
polynomial space complexity, then there exists a DTM solving the same algorithmic problem  with a polynomial space as well.
Therefore the corollary is proved.
$\Box$
\medskip

{\bf Proof of Corollary \ref{realiz}.} The word problem in a 1-element group $S$ is
linear space decidable. But one can force to solve this problem with space
complexity $f(n)$ of given deterministic machine $M.$ For this goal, the machine $M_0$
from Subsection \ref{inout} should do the following extra work. Given an input word $uv'$
of length $n,$ then in the beginning, $M_0$  let machine $M$ to use extra tapes and to accept or to reject in consecutive order all words $w$ of length $\le n$ in the tape alphabet
of $M$. Clearly the space function of such a machine $M_0$ will be equivalent to $f(n).$ Then
we apply Theorem \ref{main} to complete the proof. $\Box$
\medskip

{\bf Proof of Corollary \ref{complete}.} There exists a finitely presented semigroup $S$ 
(even a group, see \cite{T} or \cite{V}) with polynomial space (PSPACE) complete word problem 
(see \cite{DK} for the definition).
By Theorem \ref{main}, $S$ is a submonoid of a finitely presented monoid $P$ with
polynomial space function, and so the word problem in $P$ is at least PSPACE hard.
On the other hand, there is a polynomial $f(n)$ such that two words $w$ and $w'$ are equal
in $P$ iff there exists a derivation $w\to\dots\to w'$ of space $\le f(\max(|w|,|w'|)).$ It follows
(as in the proof of Corollary \ref{pspace}) that there is an NTM  of space complexity $\preceq f(n)$ 
which solves the word problem in $P,$ 
and so there is a DTM  solving the same problem in polynomial time. Thus the corollary is proved. $\Box$
\medskip

{\bf Proof of Corollary \ref{alpha}.} As we mentioned in Introduction, our notion
of space function for semigroups differs from that used in \cite{O} for groups. Nevertheless
the proof of Corollary 1.7 \cite{O} can be literally repeated here to deduce the proof of Corollary \ref{alpha} from Corollary \ref{realiz}. $\Box$

\end{document}